\documentclass[11pt]{amsart}
\usepackage{amssymb}
\usepackage{graphicx}
\usepackage{xcolor} 
\usepackage{tensor}
\usepackage{fullpage} 
\usepackage{amsmath}
\usepackage{amsthm}
\usepackage{mathrsfs} 
\usepackage{verbatim}
\usepackage[colorlinks=true,citecolor=blue,linkcolor=blue]{hyperref}
\usepackage{esint}
\usepackage{enumitem}
\usepackage{todonotes}
\setlist[enumerate]{leftmargin=1.5em}
\setlist[itemize]{leftmargin=1.5em}

\providecommand{\MR}{\relax\ifhmode\unskip\space\fi MR }

\usepackage{xcolor}
\definecolor{purple}{rgb}{0.5, 0, 1}
\definecolor{orange}{rgb}{1,.5,0}

\providecommand{\MR}{\relax\ifhmode\unskip\space\fi MR }

\usepackage{seqsplit} 
\definecolor{green}{rgb}{0,0.8,0} 

\newtheorem{theorem}{Theorem}[section]

\newtheorem{lemma}[theorem]{Lemma}
\newtheorem{proposition}[theorem]{Proposition}

\theoremstyle{definition}
\newtheorem{definition}[theorem]{Definition}

\theoremstyle{definition}
\newtheorem{remark}[theorem]{Remark}

\numberwithin{equation}{section}
\newcommand{\nrm}[1]{\Vert#1\Vert}

\newcommand{\set}[1]{\{#1\}}
\newcommand{\tld}[1]{\widetilde{#1}}

\newcommand{\nnrm}[1]{{\vert\kern-0.25ex\vert\kern-0.25ex\vert #1 
		\vert\kern-0.25ex\vert\kern-0.25ex\vert}}

\newcommand{\dist}{\mathrm{dist}\,}

\newcommand{\supp}{{\mathrm{supp}}\,}

\newcommand{\0}{\emptyset}

\newcommand{\alp}{\alpha}

\newcommand{\dlt}{\delta}

\newcommand{\veps}{\varepsilon}

\newcommand{\tht}{\theta}

\newcommand{\omg}{\omega}
\newcommand{\Omg}{\Omega}

\newcommand{\bfe}{{\bf e}}

\newcommand{\bfu}{{\bf u}}
\newcommand{\bfv}{{\bf v}}

\newcommand{\bfx}{{\bf x}}
\newcommand{\bfy}{{\bf y}}

\newcommand{\bfI}{{\bf I}}
\newcommand{\bfJ}{{\bf J}}

\newcommand{\bbR}{\mathbb R}

\newcommand{\calB}{\mathcal B}

\newcommand{\calG}{\mathcal G}

\newcommand{\calI}{\mathcal I}
\newcommand{\calJ}{\mathcal J}
\newcommand{\calK}{\mathcal K}

\newcommand{\calM}{\mathcal M}

\newcommand{\calS}{\mathcal S}

\newcommand{\varep}{\varepsilon}
\newcommand{\dd}{\mathrm{d}}

\newcommand{\bfone}{\mathbf1}
\newcommand{\ift}{\infty}
\newcommand{\q}{\mbox{ }}
\newcommand{\qd}{\quad}

\begin{document}
	
	\title{Existence and stability of Sadovskii vortices:\\
    from patch to smooth vortices}
	
	\author{Ken Abe}
	\address{Department of Mathematics, Graduate School of Science, Osaka Metropolitan University, 3-3-138 Sugimoto, Sumiyoshi-ku Osaka, 558-8585, Japan.}
	\email{kabe@omu.ac.jp}
	
	\author{Kyudong Choi}
	\address{Department of Mathematical Sciences, Ulsan National Institute of Science and Technology, 50 UNIST-gil, Eonyang-eup, Ulju-gun, Ulsan 44919, Republic of Korea.}
	\email{kchoi@unist.ac.kr}
	
	\author{In-Jee Jeong}
	\address{School of Mathematics, Korea Institute for Advanced Study, 85 Heogi-ro, Seoul 02455, Republic of Korea.}
	\email{ijeong@kias.re.kr} 
	
	\author{Young-Jin Sim}
	\address{School of Mathematics, Korea Institute for Advanced Study, 85 Heogi-ro, Seoul 02455, Republic of Korea.}
	\email{yjsim@kias.re.kr}
	
	\author{Kwan Woo}
	\address{Departement Mathematik und Informatik, Universit\"at Basel, CH-4051 Basel, Switzerland.}
	\email{kwan.woo@unibas.ch}
	
	\date{\today}
	
	\renewcommand{\thefootnote}{\fnsymbol{footnote}}

	\renewcommand{\thefootnote}{\arabic{footnote}}


	\begin{abstract}  
    We establish a scaling-invariant variational framework for steadily translating dipoles of the two-dimensional incompressible Euler equations. Specifically, we consider the maximization of the kinetic energy subject to constraints on the impulse and the $L^p$-norm ($1<p\leq \infty$) of vorticity without imposing any restriction on the total mass. In contrast to variational constructions with a mass constraint, the shape of the vortices in our framework is determined directly by the scaling-invariant structure of the energy.

    We prove that every vortex arising from this variational principle necessarily touches the symmetry axis and is therefore what is known in the physics literature as a \emph{Sadovskii vortex}, a configuration known to emerge as the endpoint of steady vortex-dipole branches. The construction is based on a sharp energy inequality under the two constraints of fixed impulse and $L^p$-norm, combined with an application of the concentration--compactness principle. This yields a family of axis-touching solutions parameterized by $p \in (1,\infty]$, interpolating between the classical Sadovskii vortex patch (corresponding to $p=\infty$) and vortices of arbitrarily high regularity as $p \to 1$. In particular, the classical Chaplygin--Lamb dipole (corresponding to $p=2$) is recovered within this family.

    Previous variational constructions treated only particular cases, notably $p=2$ and $p=\infty$, and were formulated under an additional mass constraint. Under this mass constraint, the natural scaling leads to the restriction $p>4/3$. By removing the mass constraint, we obtain a unified scaling-invariant variational principle valid for all $1<p\le\infty$. The endpoint $p=1$ nonetheless remains a natural limiting threshold for the present framework.

    As a consequence of the variational structure, we establish {an orbital} stability result, demonstrating that the axis-touching geometry persists under small perturbations. Finally, we derive a quantitative bound on the horizontal center of mass of perturbed solutions, showing that they propagate at nearly the same speed as the underlying Sadovskii vortex.
	\end{abstract}
    
	\maketitle

	\tableofcontents 
	
	\section{Introduction} 
	
	\subsection{Sadovskii Vortices}
	
	The motion of incompressible fluids is often best understood in terms of rotational motion rather than the velocity field itself.
	In two-dimensional incompressible and inviscid flows, this rotational character is described by the vorticity, which is transported by the flow.
	This viewpoint is particularly useful for understanding the formation and evolution of coherent structures. 
	
	One of the classical configurations is a counter-rotating vortex pair, which appears in various physical settings such as aircraft wakes and vortex ring interactions.
	A particularly striking configuration arises when a counter-rotating pair of vortices touch each other across the symmetry axis while translating steadily--we term such a configuration as a \textit{Sadovskii vortex}.
	Numerical and experimental studies have repeatedly identified this touching state as the endpoint of vortex-pair branches, indicating its special role in the structure of steady Euler flows \cite{FW2012, LDW2016}.
	As the translation speed increases for symmetric vortex pairs with fixed impulse, the branch of steady dipole solutions terminates precisely at this touching configuration.
	
	In the special case when the vorticity is uniform on its support (i.e. an odd symmetric pair of vortex patches), such configurations were first numerically computed by Childress \cite{Child1966}, Goldshtik \cite{Gold}, and Sadovskii \cite{Sad1971}. Their primary motivation came from the Prandtl--Batchelor theory \cite{Pra,Bat}, which predicts the emergence of such patch configurations in the inviscid limit; see extensive numerical studies in \cite{McL,Chernyshenko1993,Chernyshenko.1988}. Surprisingly, Sadovskii vortex patches appear as asymptotic cross-sectional profiles in axisymmetric three-dimensional flows, for instance in head-on collisions of vortex rings \cite{Childress.2007,Childress.2008,CG2018, Car2002, SLF2008, CGV2016, CJ2021, GMT2023, LJ2025}.
	In such settings, stability of Sadovskii patches may control long-term vorticity growth, further motivating their mathematical study \cite{CJS2025,HT2025}. 
	
	In the current work, we use the term Sadovskii vortex to describe any \textit{touching} vortex dipoles in $\bbR^{2}$ which are steadily traveling (Definition \ref{def: Sadovskii vortex} below), including the Sadovskii patch and the Chaplygin--Lamb dipole \cite{Cha1903,Lamb1906} as special cases. A motivation for this generalization comes from the study of viscous vortex dipoles, where the occurrence of such objects has been repeatedly reported, even when initially the dipoles are well-separated from each other \cite{MV1994,Nielsen,Trieling,Sipp,DeRo2009,DoGa}. While earlier studies suggested that the appearing object is well described by the Chaplygin--Lamb dipole (e.g. \cite{Nielsen,Sipp}), a more recent study \cite{DeRo2009} shows a deviation of the vorticity profile from that of the Chaplygin--Lamb dipole.

	These considerations motivate a systematic study of dynamical properties of \textit{general} Sadovskii vortex dipoles. However, even the problem of construction has been open, except for the very special cases of the Chaplygin--Lamb dipole and the Sadovskii vortex patch (solved only very recently in \cite{CJS2025,HT2025}). The main technical difficulty in their construction comes from the defining property that these vortices are touching each other. 
	
	In this work, we develop a new variational framework which provides the existence and stability of a family of Sadovskii vortices. Our key finding is that these structures are obtained by solving the kinetic energy maximization problem \textit{without} a constraint on the total circulation (equivalently, the $L^1$-mass under the sign condition); more specifically, Sadovskii vortices are characterized by the energy maximizers, with constraints on the impulse and $L^p$-norm of the vorticity for $p>1$. In other words, given these two constraints, the Sadovskii vortices have the largest total circulation, at least among energy-maximizing vortex configurations. 
    This ``largest circulation'' property suggests a structural mechanism underlying their appearance in various settings, including collisions of anti-parallel vortex rings and the evolution of viscous dipoles.

	To precisely define the concept of Sadovskii vortices used in the paper, we begin by recalling the incompressible Euler equations in vorticity form in $\bbR^2$: 
	\begin{equation}\label{eq: Euler eq.}
		\left\{
		\begin{aligned}
			&\partial_{t} \omega + \bfv \cdot \nabla \omega = 0,\\
			&\bfv=\nabla^\perp(-\Delta_{\bbR^2})^{-1}\omega,\\
			&\omega|_{t=0} = \omega_{0},
		\end{aligned}
		\right.
	\end{equation} where $\nabla^\perp = (\partial_{x_2}, -\partial_{x_1})$, and $\omg_0:\bbR^2\to\bbR$ is the  given initial data. The \textit{stream function} $\psi=\psi[\omg]:=(-\Delta_{\bbR^2})^{-1}\omega$ allows us to recover the \textit{velocity} $\bfv$ from the \textit{vorticity} $\omega$ through the Biot--Savart law:
	$$
	\bfv(\bfx)  = \nabla^\perp \psi(\bfx) = -\frac{1}{2\pi} \int_{\bbR^2}\frac{(\bfx-\bfy)^\perp}{|\bfx-\bfy|^2}\omega(\bfy)  \dd \bfy, \quad \mbox{ where }\,\, \bfx^\perp=(x_2,-x_1).
	$$
	In this paper, we establish the existence and stability of the \textit{Sadovskii vortex}, which is a profile of a traveling wave solution in the form of an odd-symmetric touching dipole:

	\begin{definition}\label{def: Sadovskii vortex}
		A function $\overline{\omega} \in L^\infty(\bbR^2)$ is called a \textit{Sadovskii vortex} if it is \textbf{compactly supported} and satisfies the following properties:
		\begin{itemize}[itemsep=1em, topsep=1em]
			\item[(i)] \textbf{Traveling wave.} 
			There exists a constant $W > 0$ such that  
			\[
			\omega(t, \bfx) := \overline{\omega}(\bfx - (W, 0)t)
			\]
			solves the Euler equations \eqref{eq: Euler eq.}.  
			We refer to $W$ as the traveling speed of $\overline{\omega}$.
			
			\item[(ii)] \textbf{Odd symmetric dipole.}
			The vorticity $\overline{\omega}$ is odd in the $x_2$-variable:
			\[
			\overline{\omega}(x_1, x_2) = -\overline{\omega}(x_1, -x_2) \geq 0 
			\quad \mbox{for }\, x_1 \in \bbR, \,\, x_2 > 0.
			\]
			
			\item[(iii)] \textbf{Touching.} The dipole exhibits full contact; there exist $r > 0$ and $c \in \bbR$ such that  
			\begin{equation}\label{touching_scale}
			\left\{ \bfx \in \bbR^2 : |\bfx - (c, 0)| < r \right\} \subset \supp \overline{\omega}.
			\end{equation}
		\end{itemize}
	\end{definition}
	For convenience, we shall occasionally use the term \emph{Sadovskii vortex} also for its restriction to the upper half-plane $\mathbb{R}^2_+$, i.e. $\omega := \overline{\omega}|_{\bbR^2_+}$, whenever no confusion arises.
	We remark that our definition of a Sadovskii vortex deviates from its meaning in the applied literature; we defer this discussion to Section \ref{subsec:refs} below.
	
	\medskip

	For any general odd-symmetric vortex dipole, the vorticity in the upper half-plane generates the motion of its symmetric counterpart in the lower half-plane through the Biot--Savart law, and vice versa. As the two components approach each other, this mutual influence becomes increasingly singular. A Sadovskii vortex, in particular, represents the limiting case in which the two counter-rotating vortices are in full contact along the symmetry axis. This contact intensifies the nonlocal interaction between two blobs to its extreme, yet the configuration retains symmetry and translates steadily without deformation.
	Such a phenomenon reveals a remarkably delicate balance between nonlinearity and symmetry, raising questions about their existence, regularity, stability, and structural rigidity.

	\subsection{Main Results}
	
We construct a one-parameter family of steadily translating
solutions to the two-dimensional Euler equations, referred to as the \textit{Sadovskii vortex family}, which contains, as special cases,  classical examples such as the Chaplygin--Lamb dipole and the Sadovskii vortex patch. 
This family is characterized by an energy-maximization principle 
without a mass constraint, in contrast with previous variational approaches.
By removing the structural rigidity imposed by mass constraints, our formulation reveals a genuinely scaling-invariant variational structure capable of treating the axis-touching regime.     
As a consequence, it yields a continuous transition between smooth and patch-type vortices, parameterized by the exponent $p\in(1,\infty]$. 
The existence and properties of this family are summarized in the following theorem.

\medskip
\noindent More detailed statements are deferred to Theorem~\ref{thm_main_1} (also see Theorem~\ref{thm_stability of patch}).

	\begin{theorem}[Theorem~\ref{thm_main_1}: Sadovskii vortex family]\label{thm_intro_main}
		For each $p \in (1,\infty]$, there exists a nonempty set $\mathcal{S}(p)$ of vorticities defined on $\bbR^2_+$ satisfying the following:
		\begin{itemize}[itemsep=0.8em, topsep=0.8em]
			
			\item[(i)] \textbf{Sadovskii vortex.}  	 The odd extension of every element in $\mathcal{S}(p)$ is a Sadovskii vortex with traveling speed $W=W(p)$ and touching scale $r=r(p)$ satisfying \eqref{touching_scale} in Definition~\ref{def: Sadovskii vortex}.

			\item[(ii)] \textbf{Regularity transition: from patch to smooth vortices.}  
			Every element of $\mathcal{S}(\infty)$ is a vortex patch, i.e. for each $\omg \in \mathcal{S}(\infty)$, there is $A \subset \mathbb{R}^{2}_{+}$ such that $\omg = \mathbf{1}_{A}$.
			As $p \searrow 1$, the regularity of the corresponding Sadovskii vortices increases. In particular, for every $m \in \mathbb{N}$, there exist Sadovskii vortices of class $C^m(\mathbb{R}^2)$.

			\item[(iii)] \textbf{{Orbital} stability.}  
			For any $1 < p < \infty$, the set $\mathcal{S}(p)$ is {orbitally} stable with respect to the norm  $\|\cdot\|_{L^p} + \|x_2(\cdot)\|_{L^1}$.
            {
            For $p = \infty$, stability in the norm $\|x_2(\cdot)\|_{L^1}$ holds under the smallness assumption for the excess $L^q$-norm $1 < q < \infty$ on the initial data, which is specified in Theorem~\ref{thm_main_1}-(iii).
            }
		\end{itemize}
		More precisely, each $\omega \in \mathcal{S}(p)$ maximizes the kinetic energy under unit impulse and unit $L^p$ norm constraints, subject to positivity on $\bbR^2_+$ and odd-symmetry, with no restriction on the total mass.
	\end{theorem}

The mechanism behind Theorem~\ref{thm_intro_main} is a scale-invariant variational structure that couples the impulse and the $L^{p}$-norm through the sharp energy inequality 
\[
(\mbox{the kinetic energy induced by }\omega) \lesssim_{p} \|x_2\omega\|_{L^1}^{\alpha}\|\omega\|_{L^p}^{\beta},
\]
established in this work (see \eqref{eq_0704_02} in Lemma~\ref{lem_energy}). 
The two exponents $\alpha$ and $\beta$ depend only on the choice of $p\in(1,\infty]$ and are uniquely determined by the requirement that the estimate be invariant under the natural two-parameter $(\varepsilon,\lambda)$ scaling 
\begin{equation}\label{eq:two_sca}
   \omega(x) \,\mapsto\,\varepsilon\, \omega(\lambda x),  
\end{equation}
which reflects a basic symmetry of the Euler equations. 
Consequently, the resulting homogeneity balance enforces compactness of maximizing sequences, ruling out both concentration and vertical escape to infinity.
In particular, the axis-touching configuration emerges as a necessary extremal geometry selected by the variational structure, rather than being imposed \emph{a priori}. 
The detailed analytic ingredients underlying this mechanism are presented in Subsection~\ref{sec_key}. Beyond the present setting, this scale-invariant variational perspective is expected to extend beyond planar dipoles to related models, including axisymmetric Euler flows in $\mathbb{R}^3$ and the two-dimensional surface quasi-geostrophic equation.

As $p \searrow 1$, the variational problem enters a highly regular regime.
In this regime, the construction of smooth dipoles whose vorticity remains in contact with the symmetry axis is particularly delicate, since the free boundary meets the axis and the configuration lies at the borderline between regular and singular geometries. Obtaining high regularity in this touching regime therefore requires a mechanism that remains stable under such geometric degeneracy. Our variational framework achieves precisely this, producing smooth axis-touching dipoles within the Sadovskii family. We note that smooth dipoles whose components remain separated from the symmetry axis have also been constructed in the literature, e.g. \cite[Theorem~11]{SS2010}. 
Such configurations avoid the geometric degeneracy inherent in the touching case considered here.

While the regime $p \searrow 1$ leads to increasingly smooth configurations, the endpoint $p=\infty$ corresponds to the Sadovskii vortex \textit{patch}, the most singular member of the family.
At this endpoint, the stability notion stated in Theorem~\ref{thm_intro_main}-(\textit{iii}) is weaker than that in the regime $p<\infty$ since the impulse norm $\|x_2(\cdot)\|_{L^1}$ alone does not sufficiently capture the dynamical behavior of perturbations.
In this regard, the following theorem provides stronger control over the dynamics.
Here we present a {normalized} version of Theorem~\ref{thm_stability of patch} and Lemma~\ref{lem_shift estimate} in Section~\ref{sec:patch_stability}.

\medskip
\noindent The set $\calS'(\infty)$ is a subset of $\calS(\infty)$ consisting of maximizers attaining the minimal $L^1$-mass (see \eqref{def_min_mass} for details).

	\begin{theorem}[Stability of the Sadovskii vortex patch]\label{thm_intro_patch}
		There exists a nonempty subcollection $\mathcal{S}'(\infty) \subset \mathcal{S}(\infty)$ which is stable in the following sense: 
		if we set
		\[
		\mathscr{A}:=\{A\subset\mathbb{R}^2_+ : \mathbf{1}_A\in\mathcal{S}'(\infty)\},
		\]
		then the following holds.
		
		\begin{itemize}[itemsep=1em, topsep=1em]
			
			\item[(i)] \textbf{Stability with total circulation {(Theorem~\ref{thm_stability of patch})}.}
			For every $\varepsilon>0$, there exists $\delta=\delta(\varepsilon)>0$ such that, for any bounded initial set $\Omega_0\subset\mathbb{R}^2_+$ satisfying
			\[
			\inf_{A\in\mathscr{A}} \int_{\Omega_0\triangle A}(x_2+1)\,\mathrm{d}\bfx < \delta,
			\]
			the corresponding solution 
			\(\mathbf{1}_{\Omega(t)}-\mathbf{1}_{\Omega(t)_-}\) of the Euler equations \eqref{eq: Euler eq.} with initial data 
			\(\mathbf{1}_{\Omega_0}-\mathbf{1}_{\Omega_{0-}}\) satisfies
			\[
			\inf_{A\in\mathscr{A}} \int_{\Omega(t)\triangle A} (x_2+1)\,\mathrm{d}\bfx < \varepsilon
			\quad \textrm{ for all }\, t\in\bbR,
			\]
			where
              $A_-:=\set{\bfx\in\bbR^2:\,(x_1,-x_2)\in A}$ denotes the mirror image set
            of $A$, and
            $A_1\triangle A_2$ denotes the symmetric difference between $A_1$ and $A_2$.

            \medskip
            {
            \noindent More generally, for each $q \in (1, \infty)$ and $\varepsilon >  0$, there exists $\delta =  \delta(q, \varepsilon) > 0$ such that, for any non-negative $\theta_0 \in L^\infty_c(\bbR^2_+)$ satisfying
            \[
            \inf_{\omega\in\calS'(\infty)} \left\{ \|\theta_0-\omega\|_{L^1\cap L^q} + \|x_2(\theta_0 - \omega)\|_1 \right\} <  \delta,
            \]
            the unique global solution of \eqref{eq: Euler eq.} with initial data given by the odd extension of $\theta_0$, restricted to $\bbR^2_+$ and denoted by $\theta(t)$, satisfies
            \[
            \inf_{\omega\in\calS'(\infty)} \left\{ \|\theta(t)-\omega\|_{L^1\cap L^q} + \|x_2(\theta(t)-\omega)\|_1 \right\} <  \varepsilon \quad\mbox{for all }\,\, t \in \bbR.
            \]
            }
			
			\item[(ii)]\textbf{Horizontal displacement {(Lemma~\ref{lem_shift estimate})}.}
			Let $P(t)$ denote the horizontal center-of-mass displacement of the patch solution $\bfone_{\Omega(t)}$ in (i), measured by
			\[
			P(t) := \fint_{\Omega(t)} x_1 \,\dd \bfx - \fint_{\Omega_0} x_1 \, \dd \bfx.
			\]
			There exists $\varepsilon_0 > 0$ such that, by assuming $\varepsilon \in (0,\varepsilon_0)$ in (i), $P(t)$ satisfies the quantitative estimate
			\begin{equation}\label{est_shift_intro}
				|P(t)-W_\ift t| \lesssim \varepsilon |t|
				\quad \textrm{ for all }\, t\in\bbR,
			\end{equation}
			where $W_\ift$ is the {common} traveling speed of vortices in $\calS(\ift)${, as shown in Lemma~\ref{lem_pohozaev}}.
		\end{itemize}
	\end{theorem}
	
	The proof of Theorem~\ref{thm_intro_patch}-(\textit{i}) above follows directly from a more comprehensive statement, Theorem~\ref{thm_stability of patch}, where the class of perturbed initial data is generic $L^\ift_c(\bbR^2)$-functions that are nonnegative in $\bbR^2_+$.
	Theorem~\ref{thm_intro_patch}-(\textit{ii}) is a simple consequence of Lemma~\ref{lem_shift estimate}.

	\begin{remark}
		The uniqueness of the maximizer (up to an $x_1$ translation) is open, although numerical simulations suggest uniqueness of the Sadovskii patch (cf. \cite{HT2025}).
		Hence, the {orbital} stability of $\calS(p)$, $p \in (1,\infty]$, does not exclude the possibility that the Sadovskii vortex closest to the perturbed solution may vary over time.
		Nevertheless, the overall \textit{Sadovskii-like} nature of the flow, in particular its touching behavior, is preserved; for example with $p=\infty$, the set $\Omg(t)$ from any  solution  in Theorem~\ref{thm_intro_patch}-(\textit{i}) that is initially perturbed from $\mathscr{A}$ stays \emph{almost in contact} with its mirror set $\Omg(t)_-$ uniformly in time: 
		There exists a shift function $\tau(t)$ such that 
		$$
		\Big|\{
		\bfx\in\bbR^2:|\bfx-(\tau(t),0)|<r_\infty\}
		\,\Big\backslash\,\big(\Omg(t)\cup\Omg(t)_-\big)\Big|_{\bbR^2}\,\lesssim\,\varepsilon \quad \mbox{ for all }\,t\in\bbR,
		$$
		where  $r_\infty>0$ is the touching scale (radius)    satisfying \eqref{touching_scale} (for $\tau(t)$, see Theorem \ref{thm_stability of patch_shift} and Proposition \ref{lem_shift estimate_1} for details).
	\end{remark}

\begin{remark}[Steiner symmetry]\label{rmk_Steiner}
Each maximizer $\omega \in \calS(p)$, $p \in (1,\infty]$, is Steiner symmetric with respect to the $x_1$-variable: for every $x_2>0$, the map $x_1 \mapsto \omega(x_1,x_2)$ is even and nonincreasing in $|x_1|$, up to translation. 
This follows from rearrangement arguments. 
In contrast to earlier approaches such as \cite{FT1981, Bur1988, BNL2013, AC2022}, however, our proofs do not rely on this symmetry, reflecting the intrinsic nature of the variational structure.
\end{remark}	
	\subsection{Reduction to the Half-Plane and the Vorticity--Stream Function Relation} \label{subsec_sadovskii}
	
	Since Definition~\ref{def: Sadovskii vortex} requires a Sadovskii vortex to be odd symmetric in $x_{2}$, for its existence, we may study the Euler equations restricted to the upper half-plane $\bbR^2_+$:
	\begin{equation}\label{eq: half Euler eq.}
		\left\{
		\begin{aligned}
			& \partial_{t}\omega + \bfv \cdot \nabla\omg = 0 \quad \mbox{ in } \,  [0,\infty)\times\bbR^2_+,\\
			& \bfv =  \nabla^\perp(-\Delta_{\bbR^2_+})^{-1}\omega \quad \mbox{ in }\, [0,\infty)\times\bbR^2_+,
		\end{aligned}
		\right.
	\end{equation} where the stream function $\calG[\omg]:=(-\Delta_{\bbR^2_+})^{-1}\omg$ is given by
	$$
	\calG[\omg](\bfx)=\frac{1}{4\pi}\int_{\bbR^2_+}\log\left(1+\frac{4x_2y_2}{|\bfx-\bfy|^2}\right)\omg(\bfy) \, \dd\bfy,
	$$
	and the velocity $\bfv=\nabla^\perp\calG[\omg]$ satisfies $v_{2} = 0$ on the boundary $\partial\bbR^2_+ = \{ x_{2} = 0 \}$.
	
	As the Sadovskii vortex $\overline{\omg}$ gives rise to a traveling wave solution of the form 
	$$
	\omg(t,\cdot)=\overline{\omg}(\cdot-(W,0)t),
	$$
	the equations \eqref{eq: half Euler eq.} lead to the steady Euler equations for the traveling profile $\overline{\omg}$ with the moving frame $(W,0)$: 
	\begin{equation}\label{eq: steady Euler eq.}
		\left\{
		\begin{aligned}
			&(\overline{\bfv}-(W,0)) \cdot \nabla \overline{\omg} = 0 \quad \mbox{ in } \, \bbR^2_+,\\
			&\overline{\bfv} = \nabla^\perp \calG \left[\overline{\omg}\right] \quad \mbox{ in } \, \overline{\bbR^2_+}:=\{\bfx\in\bbR^2: x_2\geq0\}.
		\end{aligned}
		\right.
	\end{equation} To find solutions to this steady equation, one often prescribes a vorticity-stream relation
	\begin{equation}\label{eq_intro_vorticity}
		\overline{\omg} = f \left(\calG\big[\overline{\omg}\big]-Wx_2-\gamma \right)\quad \mbox{ in } \, \bbR^2_+
	\end{equation}
	for some non-decreasing \textit{vorticity function} $f:\bbR\to\bbR$ satisfying $f\equiv0$ in $(-\infty,0]$ and for a \textit{flux constant} $\gamma\geq0$.
	Two examples of the vorticity function are $f(s):=\mathbf{1}_{\set{s>0}}$ and $f(s):=s_+=\max\set{s,0}$.
	If $\gamma=0$,  the former corresponds to the Sadovskii vortex patch and the latter to the Chaplygin--Lamb dipole. 
	In this work, we construct steady solutions \emph{without} flux constant, that is, $\gamma = 0$ and solutions in the form 
	\begin{equation}\label{eq: traveling wave form}
		\overline{\omg}  =  \lambda \left(\calG\left[ \overline{\omg}\right] - Wx_2 \right)_+^{1/(p-1)}\qd\mbox{ in }\q\bbR^2_+
	\end{equation} for some positive constant $\lambda$.
	It turns out that having zero flux constant is linked with the \textit{touching behavior}, a distinctive feature of the Sadovskii vortices. That is, from \eqref{eq: traveling wave form}, we are able to show that 
	$$
	\{ \bfx \in \bbR^2_+: |\bfx-(c,0)| < r\}  \subset \left\{\bfx \in \bbR^2_+: \calG\left[ \overline{\omega}\right] - Wx_2>0\right\} \subset \supp \bar{\omega}
	$$
	holds for some $r > 0$ and for some $c\in\mathbb{R}$.

	\subsection{Key Ideas} \label{sec_key}
	In this section, we present the main ideas of the work. 
	
	\medskip
	
	\noindent \textbf{Energy maximization under two constraints}.
	The present study is based on a \textit{variational framework}, which selects certain steady solutions as energy maximizers under appropriate constraints.
	Here, the energy is viewed as a functional associated with a vorticity and is defined by
	\begin{equation*}
		\begin{split}
			E(\omg):= \frac{1}{2}\iint_{\bbR^2_+ \times \bbR^2_+}G(\bfx,\bfy)\,\omg(\bfx)\omg(\bfy)\,\dd\bfx\dd\bfy = \frac{1}{8\pi} \iint_{\bbR^2_+ \times \bbR^2_+}\log\left(1+\frac{4x_2y_2}{|\bfx-\bfy|^2}\right)\,\omg(\bfx)\omg(\bfy)\,\dd\bfx\dd\bfy. 
		\end{split}
	\end{equation*} 
	To obtain physically meaningful steady flows, one must impose suitable constraints, and for this constrained variational problem to be relevant for dynamical stability, the constraints should be conserved quantities of the Euler equations. From incompressibility, all the $L^{p}$ norms of $\omg$ are preserved in time, and since the log kernel in $E(\omg)$ does not decay at infinity, it has been a standard practice to constrain the \emph{total mass} (or \emph{total circulation}), i.e. $\|\omega\|_1$. The $L^{1}$ norm by itself is insufficient for at least two reasons, namely \textit{concentration} and \textit{drift to infinity in the vertical direction}. First, if $\omega|_{\bbR^2_+} \sim \varepsilon^{-2}\mathbf{1}_{B_{\varepsilon}(\tilde \bfx)} \to \delta_{\tilde \bfx}$ for some point $\tilde \bfx\in\mathbb{R}^2_+$, then $E(\omg)\sim |\log \varepsilon|$, so any sequence concentrating to a Dirac mass produces an arbitrarily large energy. While this is easily ruled out by placing an $L^{p}$ constraint with any $p > 1$, the vorticity may escape vertically to infinity, keeping all the $L^{p}$ norms constant. This again makes the energy divergent, since the interaction kernel $G(\bfx,\bfy)$ blows up as $x_2, y_2 \to \infty$. The role of the impulse constraint
	\[
	\mu(\omg) := \int_{\bbR^2_+} x_2\,\omega(\bfx)\,\mathrm{d}\bfx   < \infty
	\]
	is to prevent this vanishing. However, note that although this quantity is conserved in time for solutions of the Euler equations, unlike the $L^{p}$ norms, it is coercive only under the sign constraint $\omg\ge0$, which is indispensable in the entire process.
	
	\medskip

	\noindent\textbf{Existence of maximizer without mass constraint.}
	To explain the geometric consequence of the mass constraint, we consider the simple setting where the impulse $\mu(\omg)$ and the maximum $\nrm{\omg}_{\infty}$ are both capped to $1$.
	Then, it is
    by now standard in the variational theory of vortex patches (for instance, see the recent work \cite{CQZZ2025})
    that for $0< m \ll 1$, the energy maximizer under these constraints is given by the patch $\mathbf{1}_{A}$ where $A$ is approximately a ball centered near $(0,m^{-1})$ (up to shifts in $x_{1}$) with radius $\pi^{-1/2}m^{1/2}$.
	As one increases $m$, the maximizing patch becomes larger, closer to the $x_{1}$-axis, and looks less circular.
	Then, it is natural to ask what happens to the maximizers as $m\to\infty$. It was shown in \cite{CJS2025} that after a certain mass threshold, the maximizer does not saturate the mass constraint and becomes the Sadovskii vortex patch. This approach consists of imposing an additional mass constraint and subsequently reducing the impulse parameter to produce axis-touching configurations. 
	A similar observation was made in an earlier work \cite{AC2022} for the Chaplygin--Lamb dipole, where one just replaces $L^\infty$ by $L^2$ constraint.
	Moreover, the following principle was noted in these works: \textit{not} saturating the mass bound for the energy maximizers is equivalent to the touching property. 
	\medskip
	
	\noindent 
	These two ``case studies'' suggest that for any $p > 1$, directly maximizing the energy with just $\nrm{\omg}_{p}$ and $\mu(\omg)$ constraints may give  rise to a Sadovskii vortex.
    However, extending this strategy to the full range $p \in (1,\infty]$ presents serious challenges due to the possible slow decay of $\omega$ at infinity when $p$ is small.
	Regarding this difficulty, one of the key tools we introduce is the sharp energy inequality, which bounds the kinetic energy solely in terms of two constraints: for any $p \in (1,\infty]$, we prove 
	\begin{equation}
		\label{eq_intro_energy}
		E(\omega) \le C \|x_2\omega\|_{L^1}^{\alpha}\|\omega \|_{L^p}^{\beta},\quad \alpha=\frac{4p-4}{3p-2},\quad \beta=\frac{2p}{3p-2}
	\end{equation}
	for an optimal constant $C=C(p)$; the inequality is sharp, and equality is attained by the Sadovskii vortices (Lemma~\ref{lem_energy} and Remark~\ref{r: maxmusigma}).
    The two exponents $\alpha, \beta$ are not arbitrary: the balance of homogeneities under the natural two-parameter scaling \eqref{eq:two_sca} uniquely determines them whenever an inequality of this type holds. Consequently, the estimate naturally fits into a scale-invariant variational framework that plays a central role in our analysis.
    \medskip
	
	\noindent 
	To prove the existence of a maximizer, we apply Lions' concentration--compactness principle \cite{Lions1984} (Lemma~\ref{lem:Lions lemma}) to the sequence $\{x_2\omega_n \}$, where $\{\omega_n\}$ is an energy maximizing sequence. 
	The principle asserts that essentially only three behaviors are possible: \textit{(i) compactness {up to translations}}, \textit{(ii) vanishing}, and \textit{(iii) dichotomy}.
	Among these, since $\{\omega_n\}$ is maximizing, \textit{(ii) vanishing} is excluded by the energy decay, while \textit{(iii) dichotomy} is ruled out by the convexity structure  $$ { \alpha + \frac{\beta}{p} > 1}$$ of the exponents appearing in the sharp inequality \eqref{eq_intro_energy} (see Lemma \ref{lem_superadd} and Theorem \ref{thm_existence_new} for details).
	Consequently, \textit{(i) compactness {up to translations}} occurs, which yields the existence of a maximizer (Theorem~\ref{thm_existence_new}) and ensures the \textit{{orbital} stability} of the set $\calS(p)$ of maximizers  (Theorems~\ref{thm_main_1} and~\ref{thm_stability of patch}).
	In this procedure, the main difficulty lies in the absence of a mass constraint,
    and we had to develop several 
    inequalities for the stream function and velocity under just impulse and $L^p$ bounds.
	
	\medskip
	
	\noindent \textbf{Comparison with the penalized energy method.}
	We emphasize a fundamental difference between our formulation and the relevant \textit{penalized} energy approach.
	In the penalized framework appearing, for instance, in \cite{AC2022}, touching dipoles are obtained by solving
	\[
	\calJ_{\mu} = \calJ_{\mu}(p) = \sup_{\omega \in \calM_{\mu}(p)}E_p(\omega),
	\]
	where
	\[
	E_p(\omega) := E(\omega) - \frac{1}{p}\|\omega\|_{L^p}^p, \qquad \calM_{\mu} = \mathcal{M}_{\mu}(p):= \{ 0 \le \omega \in L^p(\bbR^2_+): \|x_2\omega\|_{L^1} = \mu\ \}.
	\]
	Such penalizations are naturally introduced when the target vorticity-stream function relation \eqref{eq_intro_vorticity} is prescribed, so that the resulting maximizer yields the desired profile, for example, the Chaplygin--Lamb dipole \cite{AC2022}.
	Here, the competition between the kinetic energy $E(\omega)$ and the penalization term $\|\omega\|_{L^p}^p$ raises both coercivity and (strict) superadditivity issues.
	That is, the following may \emph{fail}:
	\[
	\calJ_{\mu} > 0 \quad \textrm{and} \quad \calJ_{\mu} > \calJ_{\mu - \alpha} + \calJ_{\alpha} \quad \textrm{ for }\, 0 < \alpha < \mu.
	\]
	In fact, scaling analysis shows that $\calJ_{\mu}>0$ only holds for $p>4/3$ and the superadditivity is known only for $p = 2$.
	In particular, in the case $p = 4/3$, the penalized energy method does not distinguish the zero solution from nontrivial maximizers.
	Losing these structures is a serious obstacle in proving the existence of maximizers via the Lions' compactness principle.
	We refer to \cite{ACJ2025}   for a recent reformulation for the case $p=2$; the Chaplygin--Lamb dipole. See also the penalized energy framework based on the functional $E(\omega)-\mu(\omega)$ developed in \cite{Bur2005, BNL2013}, which provides another mechanism for producing maximizers.

	\medskip
	
	\noindent In our work, by directly imposing $\|x_2 \omega\|_{L^1}$ and $\|\omega\|_{L^p}$ constraints, which align with the sharp energy inequality \eqref{eq_intro_energy}, we overcome these obstructions, particularly for $p \le 4/3$.
	This allows one to construct Sadovskii vortices of arbitrarily high regularity, surpassing the $C^{2, 1}$ regularity corresponding to $p = 4/3$. 
	
	\medskip 
	
\begin{figure}[htbp]
		\centering
		\includegraphics[width=0.9\linewidth]{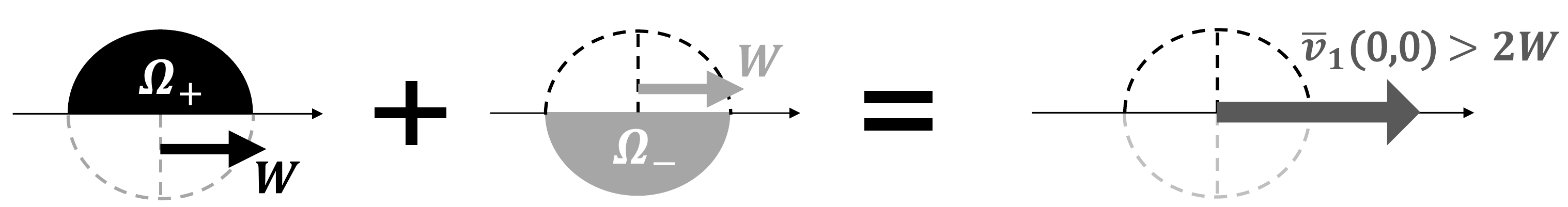}
		\caption{For an odd-symmetric dipolar vortex $\mathbf{1}_{\Omg_+}-\mathbf{1}_{\Omg_-}$, the upper region $\Omg_+$ and the lower region $\Omg_-$ interact to produce the steady flow $(W,0)$ that carries the dipole.
			Their mutual interaction ($\Omg_+$ pushes $\Omg_-$ and vice versa) overlaps at $(0,0)$, enhancing the central velocity so that $\overline{v}_1(0,0) > 2W$.
		} 
		\label{fig: center speed}
	\end{figure}

	\noindent\textbf{Euler--Lagrange equation.}
	Given the existence of a maximizer $\omega$, we derive in Section~\ref{sec:EL} the associated \textit{Euler--Lagrange equation}.
	This functional identity serves as a central tool for understanding the structure of the maximizer and, in particular, for identifying the Sadovskii vortex configuration.
	Formally, introducing Lagrange multipliers $W$ and $\kappa$ for the respective constraints, the extreme value condition
	\[
	\big(E(\theta)-W\|x_2\theta\|_{L^1}-\kappa\|\theta\|_{L^p}^{p}\big)'\Big|_{\theta = \omega}=0
	\]
	suggests the relation
	\begin{equation}\label{eq_intro_EL}
		\kappa\,\omega^{p-1}=\big(\mathcal{G}[\omega]-Wx_2\big)_+.
	\end{equation}
	This \textit{Euler--Lagrange equation} captures the balance between the nonlocal self-interaction of the maximizer $\omega$, its stream function $\mathcal{G}[\omega]$, and the traveling speed $W$.
	In particular, we get information on the location of the support:
	\begin{equation}\label{eq_intro_supp}
		\{ \omega > 0 \} = \{ \mathcal{G}[\omega]-Wx_2 > 0 \}=\{\mathcal{G}[\omega]/x_2 > W\}.
	\end{equation}
	Combined with properties of $\mathcal{G}[\omega]$, the relation \eqref{eq_intro_EL} successively yields the sharp regularity of $\omega$ (Lemma~\ref{lem_regularity_1} and Remark~\ref{rmk_optimal_regularity}),
	the Pohozaev-type identity (Lemma~\ref{lem_pohozaev}), 
	the compact-support property (Lemma~\ref{lem_cpt_support}), 
	and the uniform mass bound (Remark~\ref{lem_unif}) in order.  
	In particular, the information \eqref{eq_intro_supp} serves as the key ingredient in verifying the \textbf{touching behavior} (Lemma~\ref{prop_touching}) of the Sadovskii configuration.  See Figure~\ref{fig: center speed} for an illustration.\\
	
	\medskip
	
	\noindent\textbf{Further stability: the minimal $L^1$-mass condition.}
	For $p \in (1,\infty)$, the {orbital} stability follows from the strong convergence of maximizing sequences in the norm $\|\cdot\|_{L^p} + \|x_2(\cdot)\|_{L^1}$.
	In contrast, when $p=\infty$, the compactness becomes weaker: we only have strong convergence in the impulse norm $\nrm{x_2(\cdot)}_{L^1}$ along with weak-$^*$ convergence in $L^\infty$, which yields only a restricted form of stability for Sadovskii patches.
	
	\medskip
	
	\noindent To recover a meaningful notion of stability, we exploit a special feature of the case $p=\infty$: every maximizer is a \emph{patch}, i.e. $\omega = \bfone_{\Omega}$.
	Since uniqueness of maximizers is not guaranteed, we focus on the subclass $\calS'(\infty)\subset\calS(\infty)$ consisting of patches $\omega'$ satisfying the \emph{minimal-mass condition} 
    among all maximizers sharing the same impulse:
    $$\|\omega'\|_{L^1} =\inf_{\omega\in \calS(\infty)}\{\|\omega\|_{L^1}\}.$$
	This restriction is introduced to retain compactness: the concentration--compactness argument ensures compactness only in the impulse $\|x_2(\cdot)\|_{L^1}$ and does not exclude a possible dichotomy in $\|\cdot\|_{L^1}$.
	The minimal-mass condition, together with weak-$^*$ convergence to $\omega = \bfone_{\Omega}$ in $L^\infty$,  prevents such splitting and thereby restores compactness in $\|\cdot\|_{L^1}$.

	\subsection{Related Topics} \label{subsec:refs}
	This work is closely related to previous studies on traveling vortex solutions and their properties.
	We briefly highlight a few particularly relevant directions below.\\

	\noindent\textbf{On the terminology \textit{Sadovskii vortex}.}
	In  his 1971 paper \cite{Sad1971}, Sadovskii considered the problem of finding a steady streamline connecting two points on $\partial\bbR^2_+$ separating a constant vorticity region from irrotational flow in $\bbR^2_+$.
	He actually computed a family of such solutions, parameterized by the jump size of Bernoulli's constant at the streamline, obtaining the (pure) patch as the limiting case when the jump vanishes.
	We note that a few years earlier than Sadovskii, Childress \cite{Child1966} studied similar objects, for some simple approximations of the Euler equations.
	Since then, similar configurations have appeared under the terminology ``Sadovskii vortex'' in numerous studies~\cite{Pierrehumbert1980, ST1982, MST1988, Chernyshenko1993, AZ2021, CJS2025, HT2025}, concerning various types of boundaries and external flows.
	Often this terminology refers to objects with accompanying vortex sheets (\cite{Freilich2016, FL2017}), but in the present work, we adopt the term only for vortices without vortex sheets\footnote{A rigorous construction of Sadovskii vortices surrounded by a vortex sheet  appears to remain open.}, and in a slightly generalized sense: an odd-symmetric pair of touching vortices whose upper-half component is a nonnegative, Steiner-symmetric function; see Remark~\ref{rmk_Steiner} and the introduction in \cite{CJS2025_atmos}.
	Steiner symmetry is natural for energy-maximizing vortices, and indeed many known traveling-wave solutions exhibit this symmetry.
	The possible existence of traveling waves without it remains an open question.
	We refer to \cite[Section 1.1]{CJS2025} for a more detailed historical account on Sadovskii vortices.
	
	\medskip
	
	\noindent\textbf{Variational construction and stability of stationary solutions.}
	A powerful framework for constructing steady states is to view them as critical points of the kinetic energy under physically motivated constraints.
	This idea goes back to Kelvin and was formalized in the Kelvin--Benjamin principle \cite{Ben1976} and Arnold’s theory~\cite{Arn1966}.
	In the axisymmetric setting, Friedman--Turkington~\cite{FT1981} applied this principle to establish the existence of vortex rings by maximizing kinetic energy under total circulation and impulse constraints.
	Alternative formulations based on rearrangement classes were later developed by Turkington~\cite{Tur1983} and extended by Burton~\cite{Bur1987, Bur1988}, who also studied the uniqueness of circular vortex pairs~\cite{Bur1996}.
	The stability of such variationally constructed vortices has been extensively investigated within Arnold’s framework, which provides 
        stability for energy-critical steady flows (see~\cite{GS2024} for an overview).
	Examples include the orbital stability of symmetric vortex pairs~\cite{BNL2013},   stability and variational characterization of the Chaplygin--Lamb dipole~\cite{AC2022}, and subsequent refinements by~\cite{Wang2024, ACJ2025}.
	Related results for the Hill’s spherical vortex were obtained in~\cite{Choi2024}, together with the uniqueness characterization of~\cite{AF1986}, and a variational construction of the Sadovskii patch was proposed in~\cite{CJS2025}.
	For stronger notions such as asymptotic stability, we refer to~\cite{BM2015, BCM2019, IJ2022, IJ2023, MZ23}.

	\medskip

	\noindent\textbf{Comparison with Huang--Tong \cite{HT2025}.}
	Besides our variational characterization of Sadovskii vortices (see also \cite{CJS2025, ACJ2025}), an independent work by Huang and Tong~\cite{HT2025} established the existence of a touching vortex \emph{patch} by a completely different approach.
	Their method is based on a fixed-point argument for the upper boundary curve of the patch, which provides detailed \emph{geometric information} on the touching interface and a constructive description of its boundary.
	On the other hand, our framework focuses on a unified treatment of both patch and non-patch regimes, offering a variational and energetic perspective that naturally leads to a \textit{stability} result.

	\medskip

	\noindent\textbf{Desingularization of point vortices.}
	When two vortex blobs are well-separated, smooth traveling dipole solutions can be constructed by desingularizing point vortex configurations.
	These methods replace the singular dipole $\delta_{(0,1)}-\delta_{(0,-1)}$ with smooth vorticity supported on compact sets.
	Several mathematical strategies have been developed for this purpose.
	The \emph{stream-function method} prescribes vorticity as a nonlinear function of the stream function, leading to semilinear elliptic equations analyzed perturbatively or variationally~\cite{Nor1975, Yang1991, SS2010}.
	The \emph{bifurcation approach} generates curves of traveling or rotating dipoles from radially symmetric states~\cite{Wan1986, Wan1988, HMW2020}, while \emph{KAM theory} yields quasi-periodic vortex structures with strong stability~\cite{HHR2024}.
	Most constructions rely on reflection or rotational symmetry to simplify nonlinear interactions.
	Recent developments have extended these ideas to more general settings~\cite{HM2017, CWZ2020, CLZ2021, CLW2024, CJY2024, DPMP, BCM2025, DHLM2025}, illustrating how coherent structures emerge from singular models.

	\medskip
	
	\noindent\textbf{Non-uniqueness of planar flows via dipole structures.}
	The Chaplygin--Lamb dipole has recently been used as a building block in revealing striking non-uniqueness phenomena.
	In~\cite{BCK2024}, localized perturbations inspired by its structure were used to construct infinitely many non-conservative weak solutions with integrable vorticity.
	This approach makes use of the Chaplygin--Lamb dipole’s robust geometric and analytic features, demonstrating that explicit stationary solutions can serve as prototypes in the analytical study of weak Euler flows.

    {
    \medskip

    \noindent \textbf{General vorticity functions.}
    A natural extension of the present framework is to replace the $L^p$ constraint by a suitable convex Casimir constraint, leading to more general vorticity functions $f$ in \eqref{eq_intro_vorticity}.
    In particular, for a compactly supported traveling profile satisfying \eqref{eq_intro_vorticity} with $\gamma = 0$, the contact argument (Lemma~\ref{prop_touching}) does not rely on the power-law form of $f$.
    The compactness argument must be adapted, since the strict superadditivity established here relies on the homogeneity of the $L^p$ constraint.

    \medskip

    \noindent \textbf{Transition to separated dipoles.}
    In this work, varying the impulse $\mu$ only rescales the set of maximizers, so the resulting dipoles remain touching for every $\mu > 0$.
    By contrast, the variational family considered in~\cite{AC2022} maximizes an \emph{enstrophy-penalized} energy under an impulse constraint $\mu$ and an \textit{extra} mass bound.
    With the  mass bound and penalization coefficient fixed, the results of~\cite{AC2022} imply that the maximizers are touching Lamb dipoles for sufficiently small $\mu$ and separated dipoles for sufficiently large $\mu$.
    }

	\subsection{Open Problems}\label{sec: open problem}
	
	We discuss some open questions arising from our study of Sadovskii vortices.
	
	\medskip

	\noindent\textbf{Uniqueness.}
	While the existence and qualitative properties of Sadovskii vortices are established here, the question of \emph{uniqueness} remains open even under symmetry.
	Both~\cite{CJS2025} and~\cite{HT2025} constructed symmetric traveling vortex patches that touch the symmetry axis by distinct methods, yet it is unknown whether they represent the same solution.
	Uniqueness is currently known only for the special case $p=2$, corresponding to the Chaplygin--Lamb dipole~\cite{Bur1996, Bur2005}, where it was proved through an isometric transformation of the stream function and the moving-plane method.
	At present, no analogous geometric or analytic principle is known for general $p\ne 2$, and multiple symmetric maximizers with identical physical parameters cannot be ruled out.
	
	\medskip
	
    \noindent\textbf{Geometric structure.}
    The precise \emph{geometry} of the vortex support (core) remains only partially understood—how the boundary meets the symmetry axis, whether the core is convex, and how its curvature behaves. In the patch setting, the work \cite{HT2025}, based on a fixed-point argument, provides detailed geometric information on the touching interface together with a constructive description of the boundary.  Beyond the patch setting, a comparable geometric characterization remains open.
    
		\medskip
	
    \noindent\textbf{Infinitely differentiable (touching) dipole.}  
    While the present work establishes the existence and stability of Sadovskii vortices with arbitrary but finite regularity, the existence of a $C^{\infty}$ (or analytic) touching dipole remains open. It is natural to ask whether such a highly smooth touching configuration can exist within the Euler flow. The present variational approach naturally stops short of this endpoint, as the underlying energy inequality \eqref{eq_intro_energy} degenerates at the threshold $p=1$ (see Remark~\ref{rem_optimality_p}).  This may suggest a limitation of the functional setting itself, as the borderline $L^1$ regime is closely related to   measure-valued vorticities (see, for instance, general results on measure-valued Euler solutions such as \cite{Delort}; cf. \cite{GallagherGallay, Abe21_measure} in the viscous setting). Progress in this direction would likely require a more refined understanding of the energy kernel and its variational implications.

	\subsection{Organization of the Paper}
	
	The remainder of the paper is organized as follows.\\

In Section~\ref{sec:prel}, we introduce the variational framework and establish the basic estimates for the stream function and kinetic energy.
Section~\ref{sec:EL} is devoted to the derivation of the Euler--Lagrange equation.
Qualitative properties of maximizers, including regularity and the touching structure, are analyzed in Section~\ref{sec:maximizer_properties}.
The compactness of maximizing sequences is proved in Section~\ref{sec:existence}, which forms the core of the existence theory.
The proofs of the main results are presented in Section~\ref{sec:main}.
Finally, Appendix~\ref{app_shift} contains the quantitative estimate for the horizontal translation of perturbed solutions.

	\section{Variational Principle}\label{sec:prel}
	
	We show that the kinetic energy in a half-plane $\mathbb{R}^{2}_{+}= \left\{ \mathbf{x} = (x_1, x_2) \in \mathbb{R}^2 : x_2 > 0 \right\}$ is finite for $L^{p}$-integrable vorticity $\omega$ for $1<p \leq \infty$ with the finite impulse $x_2\omega\in L^{1}(\mathbb{R}^{2}_{+})$ and formulate the maximization of the kinetic energy \eqref{eq_maximization} under two constraints $$||x_2\omega||_{1}:=||x_2\omega||_{L^1(\mathbb{R}^2_+)}=\mu \quad\mbox{and}\quad||\omega||_{p}:=||\omega||_{L^p(\mathbb{R}^2_+)}=\sigma.$$
	We show that the maximum of the kinetic energy $\calI_{\mu,\sigma}(p)$ is positive, increasing, and strictly superadditive for $(\mu,\sigma)$. 
	
	\subsection{Stream Function Estimates}
	
	We define the stream function $\mathcal{G}[\omega]$ and the kinetic energy $E(\omega)$ for vorticity $\omega$ in a half plane by using the Green's function for the Laplacian subject to the Dirichlet boundary condition:  
	\begin{align}
		G(\bfx, \bfy)=- \frac{1}{2\pi} \left( \log|\bfx-\bfy|  -  \log|\bfx-\bfy^{*}| \right) = \frac{1}{4\pi} \log \left( 1 + \frac{4x_2y_2}{|\bfx - \bfy|^2} \right),  \label{eq: Green}
	\end{align}
	where $\mathbf{y}^* = (y_1, -y_2)$.
	By
	\[
	0<\log{(1+t)}\lesssim_{\alpha} t^{\alpha} \quad \textrm{ for }  \, t>0,
	\]
	where $\alpha\in (0,1]$, the Green's function satisfies the pointwise estimate
	\begin{equation}
		G(\bfx, \bfy) \lesssim \frac{x_2^\alpha y_2^{\alpha}}{|\bfx - \bfy|^{2\alpha}},\quad \bfx,\bfy\in \mathbb{R}^{2}_{+}.  \label{eq: GE}
	\end{equation} 
	We set the stream function and the kinetic energy by 
	\begin{equation}\label{eq:G}
		\begin{split}
			\calG[\omega](\bfx)= \int_{\bbR^2_+} G(\bfx, \bfy)\,\omega(\bfy)\,\dd\bfy,
		\end{split}
	\end{equation} 
	\begin{equation}
		E(\omg)= \frac12 \int_{ \mathbb{R}^{2}_{+}} 	\calG[\omega](\bfx) \,  \omg(\bfx) \,\dd\bfx  =  \frac{1}{2}\iint_{\bbR^2_+ \times \bbR^2_+}G(\bfx,\bfy)\,\omg(\bfx)\omg(\bfy)\,\dd\bfx\,\dd\bfy.  \label{eq: KE}
	\end{equation}
	
	\begin{proposition}\label{lem0801_01}
		Let $1\leq q <\infty$.
		For every  $\alpha>0$ and $\bfx=(x_1,x_2)\in\mathbb{R}^2_+$, 
		\begin{align}
			\int_{y_2>0,\;|\bfx-\bfy|<\alpha}  G \left(\bfx, \bfy\right)^q \dd \bfy  &\lesssim_q   x_2 \alp,   				\label{eq240803_2} \\
			\int_{0 < y_2 \le \alpha} G \left(\bfx, \bfy\right)^q  \dd \bfy   &\lesssim_q  x_2^{\frac{1}{2}}\alpha^{\frac{3}{2}}. \label{eq241003_1}
		\end{align}
	\end{proposition}
	
	\begin{proof}
		By changing the variables,
		\begin{align*}
			\int_{y_2>0,\;|\bfx-\bfy|<\alpha} G(\bfx,\bfy)^q \,\dd\bfy
			\lesssim &x_2^2 \int_{|\bfy|<\frac{\alpha}{x_2}} \left| \log\left(1+\frac{4(1-y_2)}{|\bfy|^2}\right)\right|^q \dd\bfy \\ \lesssim &x_2^2 \int_0^{\frac{\alpha}{x_2}} r \left| \log \left(1+\frac{4(r+1)}{r^2}\right)\right|^q \dd r.
		\end{align*}
		The integrand is bounded for all $r>0$ and \eqref{eq240803_2} follows. We estimate 
		\[
		\int_{0<y_2\le\alpha} G(\bfx,\bfy)^q \,\dd\bfy 
		\;\lesssim\; \int_{0<y_2\le\alpha} \left|\log \left(1+\frac{4\alpha x_2}{|\bfy|^2}\right)\right|^q \dd\bfy.
		\]
		By changing the variable, the right-hand side is bounded by
		\[
		\alpha^{\frac{3}{2}}x_2^{\frac{1}{2}} \int_0^\infty \left|\log\!\left(1+\frac{1}{t^2}\right)\right|^q \dd t \lesssim \alpha^{\frac{3}{2}}x_2^{\frac{1}{2}},
		\]
		and \eqref{eq241003_1} follows. 
	\end{proof}

	\begin{proposition}\label{prop: IE}
		Let $1<r\leq p\leq \infty$. For $\omega\in L^{p}(\mathbb{R}^{2}_{+})$ satisfying $x_2\omega\in L^{1}(\mathbb{R}^{2}_{+})$,  
		\begin{align}
			||x_2^{\alpha}\omega||_{r}\leq ||x_2\omega||_{1}^{\alpha}||\omega||_{p}^{1-\alpha},\quad \alpha=\frac{p-r}{r(p-1)}.  \label{eq: IE}
		\end{align}
	\end{proposition}
	
	\begin{proof}
		For $q=1/(\alpha r)$ and $1/q+1/q'=1$, we apply H\"older's inequality to estimate 
		\begin{align*}
			\int_{\mathbb{R}^{2}_{+}}|x_2^{\alpha}\omega|^{r}\dd\bfx=\int_{\mathbb{R}^{2}_{+}}|x_2\omega|^{\alpha r}|\omega|^{(1-\alpha) r}\dd\bfx
			\leq ||x_2\omega||_{1}^{\frac{1}{q}}||\omega||_{p}^{(1-\alpha)r}.
		\end{align*}
		By taking the $1/r$-th power of both sides, \eqref{eq: IE} follows.
	\end{proof}

	\begin{lemma}\label{lem_stream_bdd_2}
		Let $1<p\leq \infty$. Let $\omega\in L^p(\mathbb{R}^2_+)$ satisfy $\omega\geq 0$ and $x_2\omega\in L^1(\mathbb{R}^2_+)$. Then for every $\bfx=(x_1,x_2)\in\mathbb{R}^2_+$,
		\begin{equation}\label{eq_0719_01}
			0< \calG[\omega](\bfx) \lesssim_p \, x_2^{\frac{p-1}{2p-1}} \, \|x_2\omega\|_1^{\frac{p-1}{2p-1}} \, \|\omega\|_p^{\frac{p}{2p-1}},
		\end{equation}
		\begin{equation}\label{eq250718_2}
			\|\calG[\omega]\|_{\infty}\lesssim_p \|x_2\omega\|_{1}^{\frac{2p-2}{3p-2}}\|\omega\|_{p}^{\frac{p}{3p-2}}.
		\end{equation}
	\end{lemma}
	
	\begin{proof}
		For $L>0$, we set 
		\[
		\calG[\omega](\bfx) 
		= \int_{\mathbb{R}^2_+} G(\bfx,\bfy)\,\omega(\bfy)\,\dd\bfy
		= \int_{|\bfx-\bfy|>L} + \int_{|\bfx-\bfy|\le L} 
		=: I_1+I_2.
		\]
		By using \eqref{eq: GE} for $\alpha=1$,	$I_1 \lesssim x_2L^{-2}\|y_2\omega\|_{1}$. We take $\alpha \in (0, \frac{p-1}{2p-1})$ and set
		\[
		q = \frac{p}{(p-1)(1-\alpha)}, \quad q'= \frac{q}{q-1} = \frac{p}{1-\alpha+p\alpha},
		\]
		so that $\alpha q < 1$. By the pointwise estimate \eqref{eq: GE} and H\"older's inequality, 
		\[
		I_2 \lesssim x_2^\alpha
		\left(\int_{|\bfy|\le L}\frac{1}{|\bfy|^{2\alpha q}}\,\dd\bfy\right)^{\frac{1}{q}}
		\|y_2^\alpha \omega\|_{q'} \lesssim x_2^\alpha L^{-2\alpha+\frac{2}{q}}\|y_2^\alpha \omega\|_{q'}.
		\]
		By using \eqref{eq: IE} of Proposition \ref{prop: IE} for $r=q'$, $\|y_2^\alpha \omega\|_{q'} \lesssim_p \|y_2\omega\|_{1}^{\alpha}\|\omega\|_{p}^{1-\alpha}$ and
        \[
        I_2 \lesssim_p L^{-2\alpha+2/q}x_2^\alpha \|y_2\omega\|_{1}^{\alpha}\|\omega\|_{p}^{1-\alpha}.
        \]
        By combining the bounds for $I$ and $II$, 
		\[
		\calG[\omega](\bfx)\lesssim_p  L^{-2} x_2 \|y_2\omega\|_{1}
		+L^{-2\alpha+\frac{2}{q}}x_2^\alpha \|y_2\omega\|_{1}^{\alpha}\|\omega\|_{p}^{1-\alpha}.
		\]
		By taking $L>0$ so that the two terms are equal, \eqref{eq_0719_01} follows.
		
		We set  
		\[
		\calG[\omega](\bfx)=\int_{|\bfx-\bfy|\ge \frac{x_2}{2}}G(\bfx,\bfy)\omega(\bfy)\,\dd\bfy
		+\int_{|\bfx-\bfy|<\frac{x_2}{2}}G(\bfx,\bfy)\omega(\bfy)\,\dd\bfy=:J_1+J_2.
		\]
		By \eqref{eq: GE} for $\alpha=1$,
		\[
		J_1 \lesssim \frac{1}{x_2}\int_{|\bfx - \bfy| \ge \frac{x_2}{2}}  y_2\omega(\bfy)\,\dd\bfy
		\le \frac{1}{x_2}\|x_2\omega\|_{1}.
		\]
		By combining this estimate with $J_1<\mathcal{G}[\omega]$ and \eqref{eq_0719_01}, 
		\begin{align*}
			J_1 = J_1^{\frac{p-1}{3p-2}} J_1^{\frac{2p-1}{3p-2}} 
			&\lesssim_p  \left( \frac{1}{x_2}\|x_2\omega\|_{1} \right)^{\frac{p-1}{3p-2}} \left( x_2^{\frac{p-1}{2p-1}} \|x_2 \omega\|_1^{\frac{p-1}{2p-1}} \|\omega\|_p^{\frac{p}{2p-1}}\right)^{\frac{2p-1}{3p-2}} =  \|x_2 \omega\|_1^{\frac{2p-2}{3p-2}} \|\omega\|_p^{\frac{p}{3p-2}}.
		\end{align*}
		For $J_2$, we fix constants
		\[
		q = \frac{3p-2}{2p-1} \in (1, p), \quad \theta = \frac{2p-2}{3p-2} \in (0, 1)
		\]
		so that $1/q = \theta + (1-\theta)/p$ and $2/q' = \theta$. By applying H\"older's inequality and \eqref{eq240803_2},
		\begin{align*}
			J_2 \le & \, \left( \int_{|\bfx - \bfy| < \frac{x_2}{2}} G(\bfx, \bfy)^{q'}  \dd \bfy \right)^{1/q'}\left( \int_{|\bfx - \bfy| < \frac{x_2}{2}} \omega( \bfy)^{q}  \dd \bfy \right)^{1/q} \\
			\lesssim_{p} & \,x_2^{\frac{2}{q'}} \|\omega\|_{L^1(|\bfx - \bfy| < x_2/2)}^{\theta} \| \omega\|_{L^p(|\bfx - \bfy| < x_2/2)}^{1-\theta}.
		\end{align*}
		For $|\bfx - \bfy|< x_2/2$, $1< 2x_2^{-1}y_2$ and 
		\[
		\|\omega\|_{L^1(|\bfx - \bfy| < x_2/2)} \le \frac{2}{x_2} \int_{\bbR^2_+} y_2\omega(\bfy)\,\dd \bfy.
		\]
		This implies 
		\[
		J_2 \lesssim_p x_2^{\frac{2}{q'} - \theta} \| x_2 \omega\|_1^{\theta} \| \omega \|_p^{1-\theta} = \|x_2 \omega\|_1^{\frac{2p-2}{3p-2}} \|\omega\|_p^{\frac{p}{3p-2}},
		\]
		and \eqref{eq250718_2} follows.
	\end{proof}

	\subsection{Kinetic Energy Estimates}
	
	\begin{lemma}\label{lem_energy}
		Let $1<p \leq \infty$. For $\omega,\eta \in L^p(\bbR^2_+)$ satisfying $x_2\omega, x_2\eta \in L^1(\bbR^2_+)$, 
		\begin{align}
			&\iint_{\bbR^2_+ \times \bbR^2_+} G(\bfx,\bfy)\,|\omega(\bfx)\,\eta(\bfy)|\,\dd\bfx\,\dd\bfy
			\, 
			\lesssim_p \, 
			\|x_2\omega\|_{1}^{\frac{2p-2}{3p-2}}\,
			\|\omega\|_{p}^{\frac{p}{3p-2}}\,
			\|x_2\eta\|_{1}^{\frac{2p-2}{3p-2}}\,
			\|\eta\|_{p}^{\frac{p}{3p-2}},   \label{eq_0704_01}\\
			&E(\omega) 
			\lesssim_p 
			\|x_2\omega\|_{1}^{\frac{4p-4}{3p-2}}\,
			\|\omega\|_{p}^{\frac{2p}{3p-2}},		\label{eq_0704_02}  \\
			&|E(\omega)-E(\eta)| \;
			\lesssim_p\;
			\|x_2(\omega+\eta)\|_{1}^{\frac{2p-2}{3p-2}}\,
			\|\omega+\eta\|_{p}^{\frac{p}{3p-2}}\,
			\|x_2(\omega-\eta)\|_{1}^{\frac{2p-2}{3p-2}}\,
			\|\omega-\eta\|_{p}^{\frac{p}{3p-2}}.  \label{eq_0704_03}
		\end{align}
	\end{lemma}
	
	\begin{proof}
		The energy bound \eqref{eq_0704_02} follows from \eqref{eq_0704_01}. Since 
		\[
		2\left( E \left( \omega \right)-E\left( \eta \right) \right)
		=\iint_{\bbR^2_+\times\bbR^2_+} G\left( \bfx,\bfy \right)\left(\omega+\eta\right)\left( \bfx \right) \left(\omega-\eta \right)\left(\bfy \right)\,\dd\bfx\dd\bfy,
		\]
		\eqref{eq_0704_03} also follows from \eqref{eq_0704_01}. For 
		\[
		\alpha \in \left( \max \left\{ \frac{p-2}{2(p-1)}, 0 \right\} , \frac{2(p-1)}{3p-2} \right),
		\]
		we set 
		\[
		q = \frac{p}{(p-1)(1-\alpha)}, \quad q'= \frac{q}{q-1}= \frac{p}{1-\alpha+p\alpha}.
		\]
		so that $\alpha q < 2 \le q$ and $\alpha q'\le 1$. For $L>0$, we set 
		\[
		\iint_{\bbR^2_+\times\bbR^2_+} G(\bfx,\bfy)\,|\omega(\bfx)\,\eta(\bfy)|\,\dd\bfx\dd\bfy
		= \iint_{|\bfx-\bfy|>L} + \iint_{|\bfx-\bfy|\le L} =: I+II.
		\]
		By using \eqref{eq: GE}, we estimate 
		\[
		I \lesssim L^{-2}\|x_2\omega\|_{1}\|x_2\eta\|_{1}.
		\]
		We set
		\[
		h(\bfx)=\int_{\{|\bfx-\bfy|<L\}}\frac{y_2^\alpha}{|\bfx-\bfy|^{2\alpha}}\,\eta(\bfy)\,\dd\bfy.
		\]
		By using \eqref{eq: GE}, 
		\[
		II\le \int_{\bbR^2_+} h(\bfx)\,x_2^\alpha \omega(\bfx)\,\dd\bfx
		\le \|h\|_{q}\,\|x_2^\alpha\omega\|_{q'}.
		\]
		Let $r:=q/2$ and note that $r\ge1$ and $1+1/q=1/r+1/q'$.
		By Young's inequality,
		\[
		\|h\|_{q}\le \left(\int_{\{|\bfx|<L\}}|\bfx|^{-2\alpha r}\,\dd\bfx\right)^{1/r}\,\|y_2^\alpha\eta\|_{q'}
		\lesssim L^{-2\alpha+\frac{2}{r}}\,\|y_2^\alpha\eta\|_{q'}
		= L^{-2\alpha+\frac{4}{q}}\,\|y_2^\alpha\eta\|_{q'},
		\]
		since $-2\alpha r + 1  > -1$ is equivalent to $\alpha < \frac{2(p-1)}{3p-2}$.
		Moreover, by H\"older's inequality with $\alpha q'\le1$ and $(q'-\alpha q')/(1-\alpha q')=p$,
		\[
		\|x_2^\alpha\omega\|_{q'} \lesssim_p \|x_2\omega\|_{1}^\alpha\,\|\omega\|_{p}^{1-\alpha},
		\qquad
		\|y_2^\alpha\eta\|_{q'} \lesssim_p \|x_2\eta\|_{1}^\alpha\,\|\eta\|_{p}^{1-\alpha}.
		\]
		Hence
		\[
		II \lesssim_p L^{-2\alpha+\frac{4}{q}}\,
		\|x_2\omega\|_{1}^\alpha\|\omega\|_{p}^{1-\alpha}\,
		\|x_2\eta\|_{1}^\alpha\|\eta\|_{p}^{1-\alpha}.
		\]
		Combining the bounds we have
		\[
		I + II \lesssim_p 
		L^{-2}\|x_2\omega\|_{1}\|x_2\eta\|_{1}
		+L^{-2\alpha+\frac{4}{q}}
		\|x_2\omega\|_{1}^{\alpha}\|\omega\|_{p}^{1-\alpha}
		\|x_2\eta\|_{1}^{\alpha}\|\eta\|_{p}^{1-\alpha}
		\]
		and optimizing in $L>0$ yields \eqref{eq_0704_01}. 
	\end{proof}
	\begin{remark}\label{rem_optimality_p}
		The condition $p>1$ for the kinetic energy estimate \eqref{eq_0704_02} in the preceding lemma   
		is essentially optimal.
		Indeed, the endpoint case $p=1$ already fails at a very basic heuristic level.
		Formally consider a dipole consisting of two point vortices of opposite signs,
		\[
		\omega = \delta_{(0,1)} - \delta_{(0,-1)}.
		\]
		This dipole may be viewed as having   $L^1$ norm and   impulse of order one.
		However, the associated velocity field exhibits a local singularity of order
		\[
		|\bfu( \bfx)| \sim |\bfx - \bfx_0|^{-1}
		\qquad \text{near } \bfx_0=(0,\pm 1).
		\]
		As a consequence,   the total kinetic energy diverges logarithmically.		
		This observation indicates that no bound of the kinetic energy in terms of the
		$L^1$-norm and the impulse alone can hold.
		Moreover, this heuristic observation can be made fully rigorous by replacing
		the point vortices above with suitable compactly supported smooth approximations.
	\end{remark}
	
	\begin{proposition}
		\label{lem_stream_bdd}
		Let $1<p\leq \infty$. For $\omega\in L^p(\bbR^2_+)$ satisfying $x_2\omega\in L^1(\bbR^2_+)$,  
		\begin{align}
			\sup_{|\bfx| \ge R}\calG[\omega](\bfx)&\to 0 \quad \text{as } \,\, R\to\infty, \label{eq_vanishing} \\
			E(\omega)&=\frac{1}{2}\left\|\nabla \calG[\omega] \right\|_{2}^{2}.  \label{eq_L2energy}
		\end{align}
	\end{proposition}
	
	\begin{proof}
		We take a compactly supported approximation $\{\omega_n\}$ such that $|\omega_n| \le  |\omega|$, and $\omega_n$ converges to $\omega$ with the norms $||x_2(\cdot)||_1$ and $||\cdot||_p$. For an arbitrary $\varepsilon>0$, we take large $N$ such that 
		\[
		\| x_2 (\omega - \omega_N)\|_1^{\frac{2p-2}{3p-2}} \cdot \|\omega - \omega_N\|_{p}^{\frac{p}{3p-2}} \le \varepsilon.
		\]
		Since $\omega_N$ is compactly supported, there exists $R = R(\varepsilon) > 0$ such that
		\[
		\frac{1}{\operatorname{dist} (\bfx, \supp \omega_N)} \, \le \varepsilon \quad \textrm{for} \,\, |\bfx| > R(\varepsilon).
		\]
		By applying \eqref{eq250718_2} to $\calG[\omega - \omega_N]$,
		\begin{align*}
			\left| \calG[\omega](\bfx) \right| &{\leq} \left| \calG[\omega - \omega_N](\bfx) \right| + \left| \calG[\omega_N](\bfx) \right| 
			\lesssim_{p, \sigma} \, \|x_2(\omega - \omega_N)\|_{1}^{\frac{2p-2}{3p-2}}\|\omega - \omega_N\|_{p}^{\frac{p}{3p-2}}  + \left| \calG[\omega_N](\bfx) \right| \\
			&\lesssim \, \varepsilon  + \int_{\bbR^2_+} \frac{x_2y_2}{|\bfx - \bfy|^2}  \left| \omega_N(\bfy) \right| \dd\bfy,
		\end{align*}
		Thus, $\sup_{|x|\geq R(\varepsilon)} | \calG[\omega](\bfx) | \lesssim  \varepsilon$ and \eqref{eq_vanishing} follows.

		The kinetic energy is finite by \eqref{eq_0704_02}. Let $\varphi\in C_c^\infty(\bbR^2)$ be a radial cut-off function with $\varphi=1$ on $B_1$ and $\varphi=0$ outside $B_2$, and set $\varphi_R(\bfx)=\varphi(\bfx/R)$. By testing $ \varphi_R\,\calG[\omega]$ against $-\Delta \calG[\omega] = \omega$ and integrating by parts, 
		\begin{align*}
			\int_{\bbR^2_+}  \omega(\bfx) \calG[\omega](\bfx) \varphi_R(\bfx) \dd  \bfx    &= \int_{\bbR^2_+}   \left| \nabla\calG[\omega](\bfx) \right|^2 \varphi_R(\bfx) \dd  \bfx + \int_{\bbR^2_+} \frac{1}{2} \nabla \left|\calG[\omega](\bfx)\right|^2 \cdot \nabla \varphi_R(\bfx) \dd \bfx \\
			&=\int_{\bbR^2_+}   \left| \nabla\calG[\omega](\bfx) \right|^2 \varphi_R(\bfx) \dd  \bfx - \frac{1}{2}\int_{\bbR^2_+} \left|\calG[\omega](\bfx) \right|^2 \Delta \varphi_R(\bfx) \,\dd \bfx.
		\end{align*}
		Since $\calG[\omega]$ is bounded and vanishes as $|\bfx| \to \infty$, letting $R \to \infty$ implies \eqref{eq_L2energy}.
	\end{proof}

	\subsection{Maximization of the  Kinetic Energy}
	
	Let $1<p\leq \infty$ and $0 < \mu, \sigma < \infty$. We define the admissible set 
	\begin{equation}
		\begin{aligned}
			\calK_{\mu, \sigma}(p)   = 
			\left\{ \omega \ge 0: \|x_2\omega\|_1 \le \mu, \,   \|\omega\|_p \le \sigma \right\}	
		\end{aligned}
		\label{eq: AS}
	\end{equation}
	and consider the maximization problem
	\begin{equation}\label{eq_maximization}
		\calI_{\mu, \sigma}(p)= \sup_{\omega \in \calK_{\mu, \sigma}(p)} E(\omega).    
	\end{equation}
	The constant $\calI_{\mu, \sigma}(p)$ is positive and finite by \eqref{eq_0704_02}.
	By the scaling 
	\begin{equation}\label{two_scaling}
		\overline{\omega}(\bfx) = \frac{\tau^3}{\mu} \omega(\tau \bfx),\quad \tau=\left(\frac{\mu}{\sigma}\right)^{\frac{p}{3p-2}},
	\end{equation}
	the maximum satisfies 
	\begin{align}
		{\calI}_{\mu,\sigma}(p)=\calI_{1,1}(p)\mu^{\alpha}\sigma^{\beta},\quad 0< \alpha=\frac{4p-4}{3p-2}\leq \frac{4}{3},\quad \frac{2}{3}\leq \beta=\frac{2p}{3p-2}<2.  \label{eq: scaling}
	\end{align}
	We show the strict superadditivity of $\calI_{\mu,\sigma}(p)$ for two parameters, cf. \cite[IV]{Lions84b}.
	
	\begin{lemma}\label{lem_superadd}
		Let $1<p\leq \infty$ and $0<\mu_i,\sigma_i<\infty$ for $i=1,2$. Then, 
		\begin{equation}
			\begin{aligned}
				\calI_{\mu_1,\sigma_1}(p)+\calI_{\mu_2,\sigma_2}(p)&<\calI_{\mu_1+\mu_2,(\sigma_1^{p}+\sigma_2^{p})^{1/p} }(p),\quad 1<p<\infty,    \\
				\calI_{\mu_1,\sigma_1}(\infty)+\calI_{\mu_2,\sigma_2}(\infty)&<\calI_{\mu_1+\mu_2,\sigma_1\vee \sigma_2 }(\infty),\quad p=\infty,
			\end{aligned}
			\label{eq: SSA}
		\end{equation}
		where $\sigma_1 \vee \sigma_2 = \max \{ \sigma_1, \sigma_2\}$.
	\end{lemma}
	
	\begin{proof}
		For $l>1$ and $0<t<1$, 
		\begin{align*}
			t^{l}+(1-t)^{l}<1.
		\end{align*}
		For $p=\infty$, taking $l=\alpha=4/3$ and $t=\mu_1/(\mu_1+\mu_2)$ imply $\mu^{\alpha}_{1}+\mu^{\alpha}_{2}<(\mu_1+\mu_2)^{\alpha}$ and \eqref{eq: SSA} follows.
		
		We show \eqref{eq: SSA} for $1<p<\infty$. For $0<t,s<1$, we apply Young's inequality for $1<q<\infty$ with the conjugate exponent $q'$ to estimate 
		\begin{align*}
			t^{\alpha}s^{\frac{\beta}{p}}=t^{\alpha}s^{\frac{(2-\alpha)}{p}}\leq \frac{1}{q}t^{\alpha q}+\frac{1}{q'}s^{\frac{(2-\alpha)q'}{p}}.
		\end{align*}
		Observe that $0<\alpha<4/3$ and 
		\begin{align*}
			\alpha=\frac{4p-4}{3p-2}>\frac{3p-4}{3p-2}=1-\frac{2-\alpha}{p}.
		\end{align*}
		We choose $q$ satisfying $\alpha>1/q>1-(2-\alpha)/p$ so that $\alpha q>1$ and $(2-\alpha)q'/p>1$. Then, 
		\begin{align*}
			t^{\alpha}s^{\frac{\beta}{p}}+(1-t)^{\alpha}(1-s)^{\frac{\beta}{p}}\leq \frac{1}{q}t^{\alpha q}+\frac{1}{q'}s^{\frac{(2-\alpha)q'}{p}}+\frac{1}{q}(1-t)^{\alpha q}+\frac{1}{q'}(1-s)^{\frac{(2-\alpha)q'}{p}}<1.
		\end{align*}
		Taking $t=\mu_1/(\mu_1+\mu_2)$ and $s=\sigma_1^{p}/(\sigma_1^{p}+\sigma_2^{p})$ imply 
		\begin{align*}
			\mu_1^{\alpha}\sigma_1^{\beta}+\mu_2^{\alpha}\sigma_2^{\beta}
			<(\mu_1+\mu_2)^{\alpha}(\sigma_1^{p}+\sigma_2^{p})^{\frac{\beta}{p}},
		\end{align*}
		and \eqref{eq: SSA} follows.
	\end{proof}

	\begin{remark}\label{r: maxmusigma}
		The maximizers $\omega\in \calK_{\mu,\sigma}(p)$ to \eqref{eq_maximization} satisfy $||x_2\omega||_{1}=\mu$ and $||\omega||_{p}=\sigma$. Indeed, suppose that $\omega\in K_{\mu,\sigma}(p)$ is a maximizer satisfying $\mu_0=||x_2\omega||_{1}\leq \mu$ and $\sigma_0=||\omega||_{p}\leq\sigma$.
		If $\mu_0<\mu$ or $\sigma_{0}<\sigma$, we obtain a contradiction $\calI_{\mu,\sigma}(p)=E(\omega)\leq \calI_{\mu_0,\sigma_0}(p)<\calI_{\mu,\sigma}(p)$ since $\calI_{\mu,\sigma}(p)$ is increasing. Thus, $\mu_0=\mu$ and $\sigma_0=\sigma$.
		
		Similarly, non-negative sequences $\{\omega_n\}$ satisfying
		\begin{align*}
			\calI_{\mu,\sigma}(p)\leq \liminf_{n\to\infty}E(\omega_n),\quad \limsup_{n\to\infty}||x_2\omega_n||_{1}\leq \mu,\quad \limsup_{n\to\infty}||\omega_n||_{p}\leq \sigma,
		\end{align*}
		satisfy the convergence 
		\begin{align*}
			\limsup_{n\to\infty}||x_2\omega_n||_{1}= \mu,\quad \limsup_{n\to\infty}||\omega_n||_{p}= \sigma.
		\end{align*}
		Indeed, if $\mu_n=||x_2\omega_n||_{1}$ and $\sigma_n=||\omega_n||_{p}$ satisfy $\mu_n\to \mu_0< \mu$ or $\sigma_n\to \sigma_0< \sigma$, we obtain a contradiction $\calI_{\mu,\sigma}(p)\leq \liminf_{n\to\infty}E(\omega_n)\leq \liminf_{n\to\infty} \calI_{\mu_n,\sigma_n}(p)=\calI_{\mu_0,\sigma_0}(p)<\calI_{\mu,\sigma}(p)$. Thus, $\mu_n\to \mu$ and $\sigma_n\to \sigma$. 
	\end{remark}

	\section{Euler--Lagrange Equation}\label{sec:EL}

	Let $\calS_{\mu,\sigma}(p)$ be the set of maximizers to \eqref{eq_maximization}, that is,
	\[
	\calS_{\mu, \sigma}(p) := \{ \omega \in \calK_{\mu, \sigma}(p): E(\omega) = \calI_{\mu, \sigma}(p) \}.
	\]
	We show that the Euler--Lagrange equation for $\omega\in \calS_{\mu,\sigma}(p)$ is of the form 
	
	\begin{equation}
		\left\{
		\begin{aligned}
			\kappa \, \omega(\bfx)^{p-1}&=(\mathcal{G}[\omega](\bfx)-Wx_2)_{+},\quad 1<p<\infty,\\
			\omega(\bfx)&=\sigma \mathbf{1}_{\{\bfy\in\mathbb{R}^2_+\,:\,\mathcal{G}[\omega](\bfy)-Wy_2>0 \}}(\mathbf{x}),\quad p=\infty, 
		\end{aligned}
		\right.
		\label{eq: EL}
	\end{equation}
	for $s_+= \max\{s, 0\}$ and the characteristic function   
	\[
	\mathbf{1}_A(\mathbf{x}) :=  
	\begin{cases}  
		1, & \text{if } \mathbf{x} \in A, \\  
		0, & \text{if } \mathbf{x} \in \bbR^2_+\setminus A,
	\end{cases}  
	\]
	for a measurable set $A \subset \mathbb{R}^2_+$, with the Lagrange multipliers $W_p(\omega)\geq 0$ and $\kappa_{p}(\omega)>0$ associated with the constraints $||x_2\omega||_1=\mu$ and $||\omega||_{p}=\sigma$. We show the Euler--Lagrange equation \eqref{eq: EL} for the case $\mu=\sigma=1$.

	The functional 
	\[
	E: \calK= \{\omega \in L^p(\bbR^2_+): \; x_2\omega \in L^1(\bbR^2_+)\}\to [0,\infty)
	\]
	is continuous by \eqref{eq_0704_03}.
	Moreover, $E$ is Fr\'echet differentiable on $K$, and its derivative at $\omega\in \calK$ in the direction $\eta\in \calK$ is
	\begin{equation}
		\label{eq_derivative}
		\langle E'[\omega],\eta\rangle= \frac{\dd}{\dd\tau}E(\omega+\tau\eta)\Big|_{\tau=0}
		= \int_{\bbR^2_+}\calG[\omega](\bfx)\,\eta(\bfx)\,\dd\bfx.
	\end{equation}
	In the sequel, $|A|$ denotes the Lebesgue measure of a measurable set $A \subset \bbR^2$ and $L^\infty_c(\bbR^2_+)$ denotes functions in $L^\infty(\bbR^2_+)$ with compact support in $\bbR^2_+$.

	\subsection{Case \texorpdfstring{$1<p<\infty$}{1<p<infinity}}
		
		We differentiate the functional $E$ on $ \calK_{1,1}(p)$ with the constraints $\omega\geq 0$, $||x_2\omega||_{1}=1$, and $||\omega||_{p}=1$. For a maximizer $\omega\in \calS_{1,1}(p)$ and an arbitrary function $h$, we consider the function of the form $\omega(\tau)=\omega+\tau \eta +s(\tau)\eta_0$ belonging to $\calK_{1,1}(p)$ for small $\tau \geq 0$ and differentiate $E(\omega(\tau))$ at $\tau =0$. The function $h$ is modified by $\eta$ so that its impulse vanishes. We choose the function $s(\tau)$ so that $||\omega(\tau)||_{p}=1$ by using the implicit function theorem.

		\begin{proposition}\label{prop:h_aux}
			Let $1<p<\infty$ and $\omega\in \calS_{1,1}(p)$. Let $\delta_0>0$ be a constant such that $|\{\omega>\delta_0\}|>0$. Let $h_{*}\in L^{\infty}_{c}(\mathbb{R}^{2}_{+})$ be a function such that $\textrm{spt}\ h_{*}\subset \{\omega>\delta_0\}$ and $\int_{\mathbb{R}^{2}_{+}}x_2 h_{*}\dd\bfx=1$. For $h\in L^{\infty}_{c}(\mathbb{R}^{2}_{+})$, set 
			\begin{align}
				\eta=h-\left(\int_{\mathbb{R}^{2}_{+}}x_2 h\dd\bfx \right)h_{*}.   \label{eq: eta}
			\end{align}
			Then, $\eta\in L^{\infty}_c(\mathbb{R}^{2}_{+})$ satisfies 
			\begin{equation}
				\begin{aligned}
					\int_{\mathbb{R}^{2}_{+}}x_2 \eta \dd\bfx=0.   
				\end{aligned}
				\label{eq: etacond}
			\end{equation}
		\end{proposition}
		
		\begin{proof}
			The conclusion follows immediately from the definition of $\eta$.
		\end{proof}

		\begin{proposition}\label{p: existenceh_0}
        {
        Under the assumption in Proposition~\ref{prop:h_aux}, after decreasing $\delta_0$ if necessary while keeping $h_*$ fixed, there exists} a non-negative $h_0\in L^{\infty}_{c}(\mathbb{R}^{2}_{+})$ such that $\eta_0$ defined by \eqref{eq: eta} {with $h = h_0$} satisfies  
			\begin{align}
			 { h_0 = 0 \quad \textrm{on }\,\, \{0 \le \omega \le \delta_0 \} \quad \textrm{and} \quad} 	\int_{\mathbb{R}^{2}_{+}}\omega^{p-1}\eta_0 \dd\bfx\neq 0.   \label{eq: nonzero}
			\end{align} 
                        
		\end{proposition}

{
\begin{proof}
    We fix $\omega \in \calS_{1, 1}(p)$.
    Suppose that
    \[
    g := \omega^{p-1} - \left(\int_{\mathbb{R}^2_+} \omega^{p-1}h_*\dd\bfx \right) x_2 = 0 \quad \textrm{a.e. on }\,\, \{ \omega > 0 \}.
    \]
    Integrating $g$ against $\omega$ and using $\|\omega\|_p = \|x_2\omega\|_1 = 1$, we see $\int_{\mathbb{R}^2_+} \omega^{p-1}h_*\dd\bfx = 1$ and $\omega^{p-1} = x_2$ a.e. on $\{ \omega > 0 \}$.

    Fix $\delta\in(0,\delta_0)$ and take any $h \in L^\infty_c (\mathbb{R}^2_+)$ such that $h = 0$ on $\{ 0 \le \omega \le \delta \}$.
    Define $\eta$ by \eqref{eq: eta} with this $h$.
    For $\varepsilon>0$, set
    \[
    \omega_\tau = \omega + \tau(\eta - \varepsilon h_*).
    \]
    For sufficiently small $\tau > 0$, $\omega_{\tau} \ge 0$.
    Since
    \[
    \int_{\mathbb{R}^2_+}\omega^{p-1}\eta \dd\bfx = \int_{\mathbb{R}^2_+} x_2 \eta \dd\bfx = 0, \quad \int_{\mathbb{R}^2_+}\omega^{p-1}h_* \dd\bfx = 1,
    \]
    it satisfies
    \[
    \int_{\mathbb{R}^2_+}x_2\omega_\tau \dd \bfx = 1 - \varepsilon \tau, \quad \int_{\mathbb{R}^2_+}\omega_\tau^p \dd\bfx = 1 - p \varepsilon \tau + o(\tau) < 1.
    \]
    Hence $\omega_\tau\in\calK_{1,1}(p)$, and then maximality of $\omega$ yields
    \[
    0 \ge \int_{\mathbb{R}^2_+} \mathcal{G}[\omega] (\eta - \varepsilon h_*) \dd\bfx.
    \]
    Letting $\varepsilon \to 0$ and applying the same argument to $-h$, we obtain
    \[
    \int_{\mathbb{R}^2_+}(\mathcal{G}[\omega] - C_* x_2) h \dd \bfx = 0, \quad C_*=\int_{\mathbb{R}^2_+}\mathcal{G}[\omega] h_*\dd\bfx.
    \]
    Since $h$ and $\delta$ are arbitrary, $\mathcal{G}[\omega] = C_* x_2$ a.e. on $\{ \omega > 0 \}$.
    Since $\mathcal{G}[\omega] \in W^{2,p}_{\mathrm{loc}}(\mathbb{R}^2_+)$ by elliptic regularity, we obtain $\omega = -\Delta \mathcal{G}[\omega] = -\Delta (C_*x_2) = 0$ a.e. on $\{\omega > 0\}$, a contradiction.

    Thus $g$ is nonzero on a subset of $\{\omega > 0\}$ of positive measure.
    We choose $\delta_1 \in (0, \delta_0)$ and a compact set $K \subset \{\omega > \delta_1 \}$ with $|K| > 0$ on which $g$ is either strictly positive or strictly negative.
    For $h_0 = \mathbf{1}_K$, the corresponding $\eta_0$ defined by \eqref{eq: eta} with $h = h_0$ satisfies
    \begin{align*}
    \int_{\mathbb{R}^2_+}\omega^{p-1}\eta_0 \dd \bfx &= \int_{K} \omega^{p-1}  \dd \bfx - \left( \int_K x_2 \dd \bfx\right) \int_{\mathbb{R}^2_+} \omega^{p-1}h_* \dd \bfx \\
    &= \int_K g \, \dd \bfx \neq 0.
    \end{align*}
    By replacing $\delta_0$ by $\delta_1$ while $h_*$ is fixed, we complete the proof.
\end{proof}

}

		\begin{proposition}\label{p: IFT}
			Let $h_0$ and $\eta_0$ be as in Proposition \ref{p: existenceh_0}. Let $\delta \in (0,\delta_0)$.  Let $h\in L^{\infty}_{c}(\mathbb{R}^{2}_{+})$ be a function such that 
			\begin{align}
				h\geq 0\ \textrm{on}\ \{0\leq \omega\leq \delta\}.    \label{eq: h}
			\end{align}
			Then $\eta$ in \eqref{eq: eta} satisfies \eqref{eq: etacond} and 
			\begin{align}
				\eta=h\quad \textrm{on}\ \{0\leq \omega\leq \delta\}.    \label{eq: etan}
			\end{align}
			Set 
			\begin{align}
				j(\tau,s)=\int_{\mathbb{R}^{2}_{+}}|\omega+\tau \eta+s\eta_0|^{p}\dd\bfx.  \label{eq: j}
			\end{align}
			Then there exists a $C^{1}$-function $s(\tau)$ defined for small $|\tau|$ such that $s(0)=0$,  $j(\tau,s(\tau))=1$, and  
			
			\begin{align}
				s'(0)=-\frac{\int_{\mathbb{R}^{2}_{+}}\omega^{p-1}\eta \dd\bfx}{\int_{\mathbb{R}^{2}_{+}}\omega^{p-1}\eta_0\dd\bfx}.    \label{eq: s0}
			\end{align}
		\end{proposition}
		
		\begin{proof}
			It is easy to see that $\eta$ satisfies \eqref{eq: etan}. By differentiating \eqref{eq: j} for $s$,
			\begin{align*}
				\partial_sj(\tau,s)=p\int_{\mathbb{R}^{2}_{+}}|\omega+\tau \eta+s\eta_0|^{p-2}(\omega+\tau \eta+s\eta_0)\eta_0\dd\bfx.
			\end{align*}
			By \eqref{eq: nonzero}, 
			\begin{align*}
				\partial_sj(0,0)=p\int_{\mathbb{R}^{2}_{+}}\omega^{p-1}\eta_0\dd\bfx\neq 0.
			\end{align*}
			By $j(0,0)=1$ and the implicit function theorem, there exists a $C^{1}$-function $s=s(\tau)$ such that $j(\tau,s(\tau))=1$ for small $|\tau|$. By differentiating $j(\tau,s(\tau))=1$ for $\tau$, \eqref{eq: s0} follows.
		\end{proof}

		\begin{proposition}
			The function 
			\begin{align}
				\omega(\tau)=\omega+\tau \eta+s(\tau)\eta_0   \label{eq: adjusted}
			\end{align}
			satisfies{, for all sufficiently small $\tau \ge 0$,}
			\begin{equation}
				\begin{aligned}
					\int_{\mathbb{R}^{2}_{+}}\omega(\tau)^{p}\dd\bfx=1,\quad
					\int_{\mathbb{R}^{2}_{+}}x_2 \omega(\tau)\dd\bfx=1,\quad
					\omega(\tau)\geq 0.   
				\end{aligned}
				\label{eq: constraints}
			\end{equation}
		\end{proposition}

\begin{proof}
    {On a.e. $\{0 \le \omega \le \delta\}$, $h_* = h_0 = 0$ by construction.
    This with \eqref{eq: h} implies that, for $\tau \ge 0$, $\omega(\tau) = \omega + \tau h \ge 0$ on $\{0 \le \omega \le \delta\}$.}
    On $\{\omega > \delta\}$, $\omega(\tau) \ge \delta - \tau \|\eta\|_\infty - |s(\tau)| \|\eta_0\|_\infty > 0$ for sufficiently small $\tau \ge0$, since $s(0)=0$.
    Thus $\omega(\tau)$ is non-negative.
    The first condition now follows from Proposition \ref{p: IFT}.
    Since $\eta$ and $\eta_0$ satisfy \eqref{eq: etacond}, the second condition follows.
\end{proof}

		\begin{lemma}\label{prop_representation}
			For $\omega \in \calS_{1,1}(p)$, 
			there exist unique constants $\kappa=\kappa_p(\omega)>0$ and $W=W_p(\omega)\ge0$ such that \eqref{eq: EL} holds.
		\end{lemma}
		
		\begin{proof}
			We set $\psi=\mathcal{G}[\omega]$. By computation, 
			\begin{align*}
				E(\omega(\tau))=E(\omega+\tau \eta+s(\tau)\eta_0)
				&=\frac{1}{2}\int_{\mathbb{R}^{2}_{+}}|\nabla \mathcal{G}[\omega+\tau \eta+s(\tau)\eta_0]|^{2}\dd\bfx \\
				&=E(\omega)+\tau \int_{\mathbb{R}^{2}_{+}}\psi\left(\eta+\frac{s}{\tau}\eta_0\right)\dd\bfx+E(\tau \eta+s(\tau)\eta_0).
			\end{align*}
			Since $\omega\in \calS_{1,1}(p)$ is a maximizer, by using \eqref{eq: s0}, 
			\begin{align*}
				0\geq \frac{d}{d \tau}E(\omega(\tau))\Bigg|_{\tau=0}
				=\int_{\mathbb{R}^{2}_{+}}\psi\left(\eta+s'(0)\eta_0\right)\dd\bfx
				&=\int_{\mathbb{R}^{2}_{+}}\psi \eta \dd\bfx+s'(0)\int_{\mathbb{R}^{2}_{+}}\psi \eta_0 \dd\bfx \\
				&=\int_{\mathbb{R}^{2}_{+}}\psi \eta \dd\bfx-\frac{\int_{\mathbb{R}^{2}_{+}}\psi \eta_0 \dd\bfx}{\int_{\mathbb{R}^{2}_{+}}\omega^{p-1}\eta_0\dd\bfx}\int_{\mathbb{R}^{2}_{+}}\omega^{p-1}\eta \dd\bfx\\
				&=\int_{\mathbb{R}^{2}_{+}}(\psi-\kappa \omega^{p-1})\eta \dd\bfx,
			\end{align*}
			for the constant $\kappa=(\int_{\mathbb{R}^{2}_{+}}\psi \eta_0 \dd\bfx)/(\int_{\mathbb{R}^{2}_{+}}\omega^{p-1}\eta_0\dd\bfx)$. We set $W=\int_{\mathbb{R}^{2}_{+}}(\psi-\kappa \omega^{p-1}) h_{*}\dd\bfx$ and observe that 
			\begin{align*}
				0\geq \int_{\mathbb{R}^{2}_{+}}(\psi-\kappa \omega^{p-1})\left(h-\left(\int_{\mathbb{R}^{2}_{+}}x_2h\dd\bfx  \right)h_*\right) \dd\bfx
				=\int_{\mathbb{R}^{2}_{+}}(\psi-Wx_2-\kappa \omega^{p-1})h\dd\bfx.
			\end{align*}
			Since $h$ is an arbitrary function satisfying \eqref{eq: h}, 
			\begin{equation*}
				\left\{
				\begin{aligned}
					&\psi-Wx_2-\kappa \omega^{p-1}=0\quad \textrm{on}\ \{\omega>\delta\},\\
					&\psi-Wx_2-\kappa \omega^{p-1}\leq 0\quad \textrm{on}\ \{0\leq \omega\leq \delta\}.
				\end{aligned}
				\right.
			\end{equation*}
			By letting $\delta\to 0$, 
			\begin{equation*}
				\left\{
				\begin{aligned}
					&\psi-Wx_2-\kappa \omega^{p-1}=0\quad \textrm{on}\ \{\omega>0\},\\
					&\psi-Wx_2\leq 0\quad \textrm{on}\ \{\omega=0\}.
				\end{aligned}
				\right.
			\end{equation*}
            {
            We first show that $\kappa>0$.
            Set
            \[
a = \int_{\mathbb{R}^2_+} \omega^{p-1} \eta_0 \dd\bfx \neq 0, \quad v = -\frac{\eta_0}{a}.
\]
Since $\eta_0 = 0$ on $\{0 \le \omega \le \delta_0\}$, the function $\omega + tv$ is non-negative for sufficiently small
$t > 0$.
Moreover,
\[
    \int_{\mathbb{R}^2_+}x_2(\omega + tv) \dd\bfx = 1, \quad \int_{\mathbb{R}^2_+}(\omega+tv)^p\dd \bfx = 1 - pt + o(t) < 1.
\]
Thus $\omega + tv \in \calK_{1,1}(p)$.
By maximality and the definition $\kappa=(\int_{\mathbb{R}^{2}_{+}}\psi \eta_0 \dd\bfx)/(\int_{\mathbb{R}^{2}_{+}}\omega^{p-1}\eta_0\dd\bfx)$,
\[
    0\ge E(\omega + tv) - E(\omega) = t \int_{\mathbb{R}^2_+}\psi v \, \dd \bfx + t^2 E(v) = -\kappa t + t^2E(v).
\]
Since $v \neq 0$, by \eqref{eq_L2energy} we see $E(v) = \|\nabla\mathcal{G}[v]\|_2^2/2 > 0$.
Therefore $\kappa > 0$.

The preceding relations also give
\[
    \kappa\omega^{p-1} \ge \psi - Wx_2 > -W x_2 \quad \textrm{a.e. in }\, \mathbb{R}^2_+.
\]
If $W < 0$, this implies $\omega^{p-1} \ge (-W/\kappa)x_2$ a.e. in $\mathbb{R}^2_+$, contradicting $\omega \in L^p(\mathbb{R}^2_+)$.
Hence $W \ge 0$, and \eqref{eq: EL} follows.

Finally, let $\tilde{\kappa} > 0$ and $\tilde{W} \ge 0$ also satisfy \eqref{eq: EL}.
Since $\eta_0$ vanishes outside $\{\omega > \delta_0\}$, multiplying the identity $\psi - \tilde{W}x_2 = \tilde{\kappa} \omega^{p-1}$ on $\{ \omega > 0 \}$ by $\eta_0$ and integrating yields
    \[
    \tilde{\kappa}a =\int_{\mathbb{R}^2_+}\psi\eta_0 \dd\bfx =\kappa a,
    \]
where we use \eqref{eq: etacond}.
Thus $\tilde{\kappa} = \kappa$.
By integrating \eqref{eq: EL} against $\omega$, we see that $\kappa + W = 2\calI_{1,1}(p)$, and the same identity holds for
$\tilde{\kappa}$ and $\tilde{W}$.
Consequently, $\tilde{W} = W$.
            }
		\end{proof}

		\subsection{Case \texorpdfstring{$p=\infty$}{p=infinity}}
			We first show that elements of $\calS_{1,1}(\infty)$ are characteristic functions. We then consider arbitrary sets $A\subset \{\omega=1\}$ and $B\subset \{\omega=0\}$ and the function of the form $\omega(\tau)=\omega+\tau \eta$ for $\eta=-\mathbf{1}_{A}+\tau \mathbf{1}_{B}$ belonging to $\calK_{1,1}(\infty)$ for small $\tau\geq 0$ and differentiate $E(\omega(\tau))$ at $\tau =0$.

			\begin{proposition}\label{lem_EL_patch}
				For $\omega \in \calS_{1,1}(\infty)$, there exists a measurable set $\Omega \subsetneq \bbR^2_+$ such that $\omega = \mathbf{1}_{\Omega}$ in $\bbR^2_+$.
			\end{proposition}
			
			\begin{proof}
				We show that the set $\{0 < \omega < 1\}$ has measure zero. Suppose that $\left|\{0<\omega<1\}\right|>0$. Then there exists $\delta_0 > 0$ such that $\left| \{\delta_0 < \omega < 1 - \delta_0\} \right| >0$. We take $h_{*} \in L_c^{\infty}(\bbR^2_+)$ such that 
				$$
				\supp{h_{*}} \subset \{ \delta_0 < \omega < 1 - \delta_0\} \quad \textrm{and} \quad \int_{\bbR^2_+} x_2 h_{*} \, \dd \bfx = 1.
				$$ 
				We take $\delta \in (0,\delta_0)$ and $h \in L_c^{\infty}(\bbR^2_+)$ such that 
				\begin{equation*}
					\left\{
					\begin{aligned}
						h\ge0 \quad &\mbox{ on } \, \{0\le \omega \le \delta\},\\
						h\leq0 \quad &\mbox{ on } \, \{1-\delta \le \omega \le 1 \},
					\end{aligned}
					\right.
				\end{equation*}
				and set $\eta$ by \eqref{eq: eta} so that $\eta$ satisfies \eqref{eq: etacond}. Then for any sufficiently small $\tau> 0$, 
				\[
				0 \le \omega + \tau \eta \le 1 \quad \textrm{and} \quad 
				\int_{\bbR^2_+} x_2(\omega + \tau \eta) \, \dd \bfx=1.
				\]
				By using 
				\begin{align*}
					E(\omega+\tau \eta)=E(\omega)+\left(\int_{\mathbb{R}^{2}_{+}}\psi \eta  \dd \bfx \right)  \tau+E(\eta)\tau^{2},
				\end{align*}
				we observe that 
				\begin{align*}
					0 \ge 
					\int_{\mathbb{R}^{2}_{+}}\psi \eta  \dd \bfx
					=
					\int_{\bbR^2_+}\left( \psi-\left(\int_{\mathbb{R}^{2}_{+}}\psi h_{*}\dd\bfx  \right)x_2\right) h  \dd \bfx.
				\end{align*}
				Since $h$ is arbitrary on $\{\delta<\omega<1-\delta\}$, $\psi-\left(\int_{\mathbb{R}^{2}_{+}}\psi h_{*}\dd\bfx  \right)x_2 = 0$ on $\{0 < \omega < 1 \}$. Thus for $h$ supported on $\{0<\omega<1\}$, $\int_{\mathbb{R}^{2}_{+}}\psi \eta \dd\bfx=0$ and $E(\omega+\tau \eta)=E(\omega)+E(\eta)\tau^{2}>E(\omega)$. This contradicts the maximality of $\omega$. Hence $\omega = \mathbf{1}_{\Omega}$ for some measurable set $\Omega \subset \bbR^2_+$ such that $\int_{\Omega}x_2\dd \bfx=1$.
			\end{proof}

			\begin{lemma}
				\label{prop_representation_patch} 
				For $\omega \in \calS_{1,1}(\infty)$, 
				there exists a unique constant $W=W_{\infty}(\omega) > 0$ such that \eqref{eq: EL} holds.
			\end{lemma}

			\begin{proof}
				The measures of the sets $\{\omega=1\}$ and $\{\omega=0\}$ are positive by Proposition \ref{lem_EL_patch}. For arbitrary $A\subset \{\omega=1\}$ and $B\subset \{\omega=0\}$, we set 
				\begin{align*}
					\eta=-\mathbf{1}_{A}+\tau \mathbf{1}_{B},\quad \tau=\frac{\int_{A}x_2 \dd\bfx}{\int_{B}x_2 \dd\bfx}.
				\end{align*}
				Then, $\int_{\mathbb{R}^{2}_{+}}x_2 \eta \dd\bfx=0$. For sufficiently small $\tau>0$, $\omega+\tau \eta$ satisfies $0\leq \omega+\tau \eta\leq 1$ and $\int_{\mathbb{R}^{2}_{+}}x_2 (\omega+\tau \eta)\dd\bfx=1$. Since $\omega$ is a maximizer, $0\geq \int_{\mathbb{R}^{2}_{+}}\psi \eta \dd\bfx$. By substituting $\eta$ into this inequality, 
				\begin{align*}
					\frac{\fint_{B}\psi  \dd\bfx}{\fint_{B}x_2 \dd\bfx}\leq \frac{\fint_{A}\psi  \dd\bfx}{\fint_{A}x_2 \dd\bfx}.
				\end{align*}
				We choose $A=B_r(\bfx)\cap \{\omega=1\}$ by an open ball $B_r(\bfx)$ with radius $r>0$ centered at $\bfx\in \{\omega=1\}$. Since $\psi$ and $x_2$ are continuous, letting $r\to 0$ implies that the right-hand side converges to $\psi(\bfx)/x_2$. By applying a similar argument for $B$, we obtain 
				\begin{align*}
					\sup\left\{\frac{\psi (\bfy)}{y_2}\Bigg|\ \omega(\bfy)=0,\ \bfy\in \mathbb{R}^{2}_{+} \right\}
					\leq \inf\left\{\frac{\psi (\bfx)}{x_2}\Bigg|\ \omega(\bfx)=1,\ \bfx\in \mathbb{R}^{2}_{+} \right\}=:W.
				\end{align*}
				Since $\psi/x_2$ is continuous, the equality holds. Thus $\psi-Wx_2\geq 0$ on $\{\omega=1\}$ and $\psi-Wx_2\leq 0$ on $\{\omega=0\}$.
				
				Since $0\geq \int_{\mathbb{R}^{2}_{+}}\psi \eta \dd\bfx=\int_{\mathbb{R}^{2}_{+}}(\psi-Wx_2) \eta \dd\bfx$, the maximizer $\omega$ can take arbitrary values on $\{\psi-Wx_2=0\}$. We set $\omega=0$ on $\{\psi-Wx_2=0\}$. Then, $\omega=\mathbf{1}_{\{\psi-Wx_2>0\}}(\bfx)$. Since $\int_{\{\psi-Wx_2>0\}}x_2\dd\bfx$ is decreasing for $W>0$, the constant $W>0$ satisfying $\int_{\{\psi-Wx_2>0\}}x_2\dd\bfx=1$ is unique.
			\end{proof}

			\section{Properties of Maximizers}\label{sec:maximizer_properties}
			
			We show regularity of maximizers $\omega\in \calS_{1,1}(p)$ to the variational problem \eqref{eq_maximization} by studying the Euler--Lagrange equations for the stream function $\psi=\mathcal{G}[\omega]$:
			\begin{equation}
				\left\{
				\begin{aligned}
					-\Delta \psi (\bfx)&= f\left(\psi \left(\bfx \right) - Wx_2  \right)\quad \textrm{ in }\, \bbR^2_+, \\
					\psi&=0 \quad \textrm{ on }\, \partial\bbR^2_+,
				\end{aligned}
				\right.
				\label{eq_ell}
			\end{equation}
			for the function  
			\begin{equation}    
				\begin{aligned}
					f(s) = \begin{cases}
						\kappa^{-\frac{1}{p-1}} s_+^{\frac{1}{p-1}}, & \text{for }\,\, 1<p<\infty, \\[5pt]
						\mathbf{1}_{\{s > 0\}}, & \text{for}\,\, p = \infty.
					\end{cases}
				\end{aligned}
				\label{f}
			\end{equation}
			We show positivity of the Lagrange multiplier by using the identity $W_p=\frac{4p-4}{3p-2}\calI_{1,1}(p)>0$, and deduce that maximizers have compact support in $\overline{\mathbb{R}^{2}_{+}}$. We also derive a different formula for the speed $W_p$ using arbitrary maximizers $\omega\in \calS_{1,1}(p)$ and deduce that maximizers have contiguous contact with the $x_1$-axis. 
			
			\subsection{Regularity of Maximizers}
			
			For $\alpha \in (0,1)$ and a positive integer $m$, we denote by $C^{m,\alpha}(\mathbb{R}^2)$ the space of functions $f$ such that $\partial^\ell_{\mathbf{x}} f$ is uniformly bounded and continuous for $|\ell| \leq m-1$, and $\partial^\ell_{\mathbf{x}} f$ is $\alpha$-H\"older continuous for $|\ell| = m$. The space $C^{m,\alpha}(\mathbb{R}^2)$ is normed with  
			\[
			\|f\|_{C^{m, \alpha}(\mathbb{R}^2)} =  \sum_{|\ell| \leq m} \|\partial_{\mathbf{x}}^\ell f\|_{\infty} + \sum_{|\ell| = m} \sup_{\substack{\mathbf{x} \neq \mathbf{y} \\ \mathbf{x}, \mathbf{y} \in \mathbb{R}^2}} \frac{|\partial_{\bfx}^\ell f(\mathbf{x}) - \partial_{\bfx}^\ell f(\mathbf{y})|}{|\mathbf{x} - \mathbf{y}|^{\alpha}},
			\]  
			where $\ell = (\ell_1, \ell_2)$ is a multi-index with integers $\ell_j$ ($j=1,2$), and  
			\[
			|\ell| = \ell_1 + \ell_2, \quad \partial_{\bfx}^\ell f := \frac{\partial^{|\ell|} f}{\partial x_1^{\ell_1} \partial x_2^{\ell_2}}.
			\]
			We similarly define $C^{m, \alpha}(\overline{\mathbb{R}^2_+})$ and its norm. We denote by $\lceil s \rceil$ the smallest integer greater than or equal to $s \in \bbR$.
			
			\begin{lemma}\label{lem_regularity_1}
				Let $1<p \leq \infty$ and $\omega \in \calS_{1,1}(p)$. Then, $\omega\in L^{\infty}(\bbR^2_+)$ and $\psi =\calG[\omega] \in L^{\infty}(\mathbb{R}^{2}_{+})\cap W^{2,q}_{\mathrm{loc}}(\overline{\mathbb{R}_{+}^{2}})$ for all $q\in (1,\infty)$. In particular, $\psi\in C^{1,\beta}(\overline{\mathbb{R}^{2}_{+}})$ for all $\beta \in (0, 1)$. For $p\in(1,\infty)$, $\omega\in C^{m, \alpha}(\overline{\bbR^2_{+}})$ and 
				\begin{align}
					m = \left\lceil \frac{1}{p-1} - 1 \right\rceil, \quad
					\alpha = \frac{1}{p-1} - m \in (0,1]. \label{malpha}
				\end{align}
			\end{lemma}

			\begin{proof}
				By \eqref{eq250718_2}, $\psi\in L^{\infty}(\mathbb{R}^{2}_{+})$. For $p=\infty$, $\omega =\mathbf{1}_{\{\psi-Wx_2>0\}}\leq 1$. For $p \in (1, \infty)$, $\kappa \omega^{p-1}=(\psi-Wx_2)_{+}\leq \psi$. Since $\psi,\ \omega\in L^{\infty}(\mathbb{R}^{2}_{+})$ and $\psi$ satisfies \eqref{lem_regularity_1}, $\psi\in  W^{2,q}_{\mathrm{loc}}(\overline{\bbR^2_{+}})$ follows by the elliptic regularity theory. 
				
				Suppose that $1<p<\infty$. For $m=0$, $\omega=f(\Psi)\in C^{0,\alpha}(\overline{\mathbb{R}^{2}_{+}})$ for $\Psi=\psi-Wx_2$ by $\psi\in C^{1,\beta}(\overline{\mathbb{R}^{2}_{+}})$. For $m=1$ and $0<\alpha<1$, $-\Delta \Psi=f(\Psi)$ in $\mathbb{R}^{2}_{+}$ and $\Psi(x_1,0)=0$ by \eqref{lem_regularity_1}. By differentiating in the $x
				_1$-variable and applying the elliptic regularity theory, $\partial_{x_1}\Psi\in C^{2,\beta}(\overline{\mathbb{R}^{2}_{+}})$. By $\partial_{x_2}^{2}\Psi=\Delta \Psi-\partial_{x_1}^{2}\Psi=-f(\Psi)-\partial_{x_1}^{2}\Psi$, $\Psi\in C^{3,\alpha}(\overline{\mathbb{R}^{2}_{+}})$. Thus, $\omega=-\Delta \Psi\in C^{1,\alpha}(\overline{\mathbb{R}^{2}_{+}})$. For $\alpha=1$, $\Psi\in C^{3,\beta}(\overline{\mathbb{R}^{2}_{+}})$ and $\omega=f(\Psi)\in C^{1,1}(\overline{\mathbb{R}^{2}_{+}})$. 
				
				For $m\geq 2$, $\omega\in C^{m,\alpha}(\overline{\mathbb{R}^{2}_{+}})$ follows the same argument.
			\end{proof}

			\begin{remark}\label{rmk_reg_odd}
				The regularity stated in Lemma~\ref{lem_regularity_1} is preserved under odd extension across the boundary $\{x_2 = 0\}$.
				This is a consequence of the Euler--Lagrange equation \eqref{eq_ell} and the degeneracy of the nonlinearity $f$ in \eqref{f}, whose derivatives up to order $m = \left\lceil \frac{1}{p-1} - 1 \right\rceil$ vanish at the origin.
			\end{remark}

			\subsection{Identities for Lagrange Multipliers}
			
			We express the Lagrange multipliers by the maximum \eqref{eq_maximization} by using Pohozaev's identity for semilinear elliptic equations, cf. \cite[Theorem 5.1]{QS07}.
			
			\begin{lemma}
				\label{lem_pohozaev}
				Let $1<p< \infty$ and $\omega \in \calS_{1, 1}(p)$.
				The constants $\kappa_p  > 0$ and $W_p> 0$ in \eqref{eq: EL} satisfy 
				\begin{equation}
					\label{eq_pohozaev}
					W_p=\frac{4p-4}{3p-2}\calI_{1,1}(p),\quad \kappa_p =\frac{2p}{3p-2} \calI_{1,1}(p).
				\end{equation}
				For $p=\infty$, $W_{\infty}=\frac{4}{3}\calI_{1,1}(\infty)$.
			\end{lemma}

			\begin{proof}  
				We apply the following identity for arbitrary functions
				\begin{align*}
					\nabla \cdot \left(\nabla \psi (\bfx \cdot \nabla \psi)-\frac{1}{2} \bfx |\nabla \psi|^{2} \right)
					=\Delta \psi (\bfx \cdot \nabla \psi).
				\end{align*}
				For a maximizer $\omega\in \calS_{1,1}({p})$ to \eqref{eq_maximization} and $\psi=\mathcal{G}[\omega]$, the function $\Psi=\psi-Wx_2$ satisfies $-\Delta \Psi=\lambda \Psi_{+}^{\gamma}=\omega$ for $\lambda=\kappa^{-\frac{1}{p-1}}$ and $\gamma=\frac{1}{p-1}\geq 0$ with $\kappa=\kappa_{p}$ and $W=W_{p}$. By using 
				
				\begin{align*}
					\Delta \psi (\bfx \cdot \nabla \psi) = -\lambda \Psi_{+}^{\gamma} \bfx \nabla \Psi-x_2W\omega
					&=-\frac{\lambda}{\gamma+1} \bfx \cdot \nabla \Psi_{+}^{\gamma+1}-x_2W\omega \\
					&=-\frac{\lambda}{\gamma+1}\nabla \cdot (\bfx \,\Psi_{+}^{\gamma+1})
					+\frac{2\lambda}{\gamma+1} \Psi_{+}^{\gamma+1}
					-x_2W\omega \\
					&=-\frac{1}{\gamma+1}\nabla \cdot (\bfx \, \Psi \omega)+\frac{2}{\gamma+1}\psi\omega-\frac{\gamma+3}{\gamma+1}Wx_2\omega,
				\end{align*}\\
				we obtain 
				\begin{align*}
					\nabla \cdot \left(\nabla \psi (\bfx\cdot \nabla \psi)-\frac{1}{2}\bfx |\nabla \psi|^{2} +\frac{1}{\gamma+1}  \bfx  \Psi \omega \right)=\frac{2}{\gamma+1}\psi\omega-\frac{\gamma+3}{\gamma+1}Wx_2\omega.
				\end{align*}\\
				The normal trace of the vector field on the right-hand side on $\{x_2=0\}$ vanishes by the Dirichlet boundary condition of $\psi$. By using $\nabla \psi\in L^{2}(\mathbb{R}^{2}_{+})$ and $\Psi \omega \in L^{1}(\mathbb{R}^{2}_{+})$, the integral of the left-hand side vanishes. Thus 
				\begin{align*}
					W=\frac{4}{\gamma+3}\calI_{1,1}(p).
				\end{align*}
				By $\kappa=2\calI_{1,1}(p)-W$ and $\gamma=1/(p-1)$, we obtain \eqref{eq_pohozaev}.
			\end{proof}
			
			\begin{remark}
				The identities \eqref{eq_pohozaev} imply that the Lagrange multipliers $W_p>0$ and $\kappa_{p}>0$ are invariant for all maximizers $\omega\in \calS_{1,1}(p)$.
				By using the exponents $\alpha$ and $\beta$ of the maximum \eqref{eq: scaling}, the identities \eqref{eq_pohozaev} are expressed as $W_p=\alpha \calI_{1,1}(p)$ and $\kappa_{p}=\beta \calI_{1,1}(p)$. Namely, 
				\begin{align*}
					\frac{\partial}{\partial \mu}\calI_{\mu,\sigma}(p)\Big|_{(\mu,\sigma)=(1,1)}=W_p,\quad \frac{\partial}{\partial \sigma}\calI_{\mu,\sigma}(p)\Big|_{(\mu,\sigma)=(1,1)}=\kappa_p.
				\end{align*}
			\end{remark}

			\subsection{Compactness of  Support}

			We deduce from the positivity of the speed $W_p>0$ in \eqref{eq_pohozaev} that the support of $\omega\in \calS_{1,1}(p)$ is compact. We use the decay $\varphi(x)\to 0$ as $|x|\to\infty$ for $L^{2}$-integrable bounded uniformly continuous functions $\varphi\in \textrm{BUC}(\overline{\mathbb{R}^{2}_{+}})\cap L^{2}(\mathbb{R}^{2}_{+})$.
			
			\begin{lemma}
				\label{lem_cpt_support}
				Let $1<p\leq \infty$ and $\omega \in \calS_{1,1}(p)$. Then, $\supp \omega = \overline{ \{ \bfx \in \bbR^2_+: \omega > 0 \} }$ is compact.
			\end{lemma}
			
			\begin{proof}
				By \eqref{eq_ell}, $\{ \bfx \in \bbR^2_+: \omega > 0 \} = \left\{ \bfx \in \bbR^2_+: \calG[\omega](\bfx) > Wx_2 \right\}$ for some $W = W_p> 0$. By Lemma \ref{lem_regularity_1}, the function 
				\begin{align*}
					\frac{\psi(x_1,x_2)}{x_2}=\int_{0}^{1}\partial_{x_2}\psi(x_1,x_2s) \, \dd s
				\end{align*} 
				is bounded and uniformly continuous in $\overline{\mathbb{R}^{2}_{+}}$. Since $\psi/x_2\in L^{2}(\mathbb{R}^{2}_{+})$ by $\nabla \psi\in L^{2}(\mathbb{R}^{2}_{+})$ and Hardy's inequality \cite[2.7.1]{Mazya}, $\varphi=\psi/x_2\to 0$ as $|\bfx|\to\infty$.
				If there exists a sequence $\{\bfx_n = (x_{n, 1}, x_{n, 2})\}\subset \supp \omega$ such that $|\bfx_n|\to\infty$, $0<W<\psi(\bfx_n)/x_{n,2}\to 0$.
				This is a contradiction. Thus, $\supp \omega$ is compact.
			\end{proof}

			\begin{remark}\label{rmk_W_formula}
				
				Lemma~\ref{lem_cpt_support}, together with the scaling \eqref{two_scaling}, shows that any $\omega\in \calS_{\mu,\sigma}(p)$ has a finite mass, i.e. $\omega \in L^1(\bbR^2_+)$.
                This integrability allows us to introduce the following representation of the traveling speed $W_p = W(p, \mu, \sigma)$ for $1<p\leq \infty$, valid for arbitrary $\omega\in \calS_{\mu, \sigma}(p)$: 
				\begin{align}
					\label{formula_W}
					W_p
					=  \frac{1}{2\pi\|\omega\|_1} \int_{\mathbb{R}^2_+} \int_{\mathbb{R}^2_+} 
					\frac{x_2 + y_2}{|\bfx - \bfy^*|^2}\, \omega(\bfx) \,\omega (\bfy) \,  \dd \bfx \, \dd \bfy.
				\end{align}
				This identity follows from the Euler--Lagrange equation \eqref{eq: EL} and $\int_{\mathbb{R}^{2}_{+}}\omega \partial_{x_2}\Psi \dd \bfx=0$ for $\Psi=\calG[\omega]-W_px_2$. 
			\end{remark}

			\begin{remark}\label{lem_unif}
				For $1<p\leq \infty$, the mass of $\omega\in \calS_{1,1}(p)$ is uniformly bounded from below. Namely, 
				\begin{equation}
					\label{eq_unif_mass}
					\nu(p) := \inf\{||\omega||_{1} : \omega\in \calS_{1,1}(p)\}\geq C_pW_{p}>0.
				\end{equation}
				Indeed, $\sup_{\omega\in \calS_{1,1}(p)}\|\calG[\omega]\|_{\infty}\lesssim 1$ by \eqref{eq250718_2} and for $\omega\in \calS_{1,1}(p)$, 
				\begin{equation}
					\label{unif_mass_low}
					1 = \int_{\bbR^2_+} x_2 \omega \, \dd\bfx
					\leq \left(\sup_{\bfx \in \operatorname{supp} \omega} x_2\right) ||\omega||_{1}
					\leq \frac{\|\calG[\omega]\|_{\infty}}{W_p}||\omega||_{1}
					\lesssim \frac{1}{W_p}\|\omega\|_1.
				\end{equation}
				We will deduce the $L^{1}$-stability of minimal mass Sadovskii vortices from the lower bound \eqref{eq_unif_mass} in \S \ref{sec:patch_stability}.
			\end{remark}

			\begin{remark}\label{rmk_unif_bdd}
				In fact, the Euler--Lagrange equation \eqref{eq: EL} yields further uniform controls on maximizers.
				In particular, one can derive uniform upper bounds on the total mass and the size of the support:
				\begin{align*}
					\sup_{\omega \in \calS_{1,1}(p)} \|\omega\|_1,
					\quad
					\sup_{\omega \in \calS_{1,1}(p)} \big|\{\bfx : \omega(\bfx) > 0\}\big|
					\lesssim_p 1.
				\end{align*}
				The support bound relies on the Steiner symmetry of maximizers (Remark~\ref{rmk_Steiner}), which implies a uniform decay of the associated stream function in the $x_1$ direction.
				These estimates are proved in Lemma~\ref{lem_uniform}.
			\end{remark}
			
			\subsection{Contact of Support with Symmetry Axis}
			
			We show that $\supp {\omega}$ is in contiguous contact with the $x_1$-axis uniformly for all $\omega \in \calS_{1,1}({p})$. For $\omega(\mathbf{x})$ and $\mathbf{x} = (x_1, x_2) \in \mathbb{R}^2_+$, we define the odd extension to $\mathbb{R}^{2}$ by  
			\[
			\tld{\omega}(\mathbf{x}) :=  
			\begin{cases}  
				\omega(\mathbf{x}), & \text{for } \mathbf{x} \in \mathbb{R}^2_+, \\  
				-\omega(\mathbf{x}^*), & \text{for } \mathbf{x}^* \in \mathbb{R}^2_+, \\  
				0, &\text{for}\ \mathbf{x} \in \partial \mathbb{R}^2_+.
			\end{cases}  
			\]
			
			\begin{lemma}
				\label{prop_touching}    
				Let $1<p \leq \infty$. There exists a constant $r = r(p) > 0$ such that $B_{r}(\mathbf{0}) \subset \supp \widetilde{\omega}$ for $\omega \in \calS_{1,1}({p})$ by a suitable translation in the $x_1$-direction.
			\end{lemma}

			\begin{proof}
				By using \eqref{eq: Green}, we set 
				\[
				\psi(\bfx)=\calG[\omega](\bfx) = - \frac{1}{2\pi}\int_{\mathbb{R}^2_+} \log|\bfx-\bfy| \,\omega(\bfy)   \dd \bfy + \frac{1}{2\pi} \int_{\mathbb{R}^2_+} \log|\bfx-\bfy^{*}| \, \omega(\bfy)   \dd \bfy=\psi_1(\bfx)+\psi_2(\bfx).
				\]
				The functions $\psi_1$ and $\psi_2$ satisfy $\psi_1(x_1,x_2)=-\psi_2(x_1,-x_2)$. Since $\omega$ is bounded and has a compact support in $\overline{\mathbb{R}^{2}_{+}}$, $\psi_2\in C^{1}(\overline{\mathbb{R}^{2}_{+}})$ and $\partial_{x_2}\psi_2$ is bounded in $\mathbb{R}^{2}_{+}$.
				Moreover, $\partial_{x_2}\psi_2$ is positive and $\partial_{x_2}\psi_2(x_1,0)$ is decaying as $|x_1|\to\infty$.
				Since $\psi_2$ is harmonic, $\partial_{x_2}\psi_2$ takes a maximum on $\{x_2=0\}$. By translation, we may assume that $\partial_{x_2}\psi_2$ takes a maximum at the origin.
				By $\partial_{x_2}\psi_1(x_1,0)=\partial_{x_2}\psi_2(x_1,0)$, $\partial_{x_2}\psi(\mathbf{0})=2\partial_{x_2}\psi_2(\mathbf{0})=2||\partial_{x_2}\psi_2||_{\infty}$. We multiply $\omega$ by $\partial_{x_2}\psi_2$ and integrate it on $\mathbb{R}^{2}_{+}$.
				By \eqref{formula_W}, 
				\begin{align*}
					W=\frac{1}{2\pi ||\omega||_{1}}\int_{\mathbb{R}^2_+} \int_{\mathbb{R}^2_+} 
					\frac{x_2 + y_2}{|\bfx - \bfy^*|^2}\, \omega(\bfx) \,\omega (\bfy) \,  \dd \bfx \, \dd \bfy
					=\frac{1}{||\omega||_{1}}\int_{\mathbb{R}^2_+}\partial_{x_2}\psi_2 \, \omega \, \dd \bfx
					\leq \frac{1}{2}\partial_{x_2}\psi (\mathbf{0}).
				\end{align*}
				By the continuity of $\partial_{x_2}\psi$, there exists $r>0$ such that 
				\[
				\partial_{x_2}\psi > W \quad \textrm{in} \,\, B_{r}(\mathbf{0}).
				\]
				For any $\bfx=(x_1,x_2)\in B_{r}(\mathbf{0}) \cap \bbR^2_+$, 
				\begin{align*}
					\psi(\bfx) = \int_{0}^{x_2} \partial_{x_2}\psi(x_1,s) \dd s > W x_2,    
				\end{align*}
				and $\bfx \in \{\psi - Wx_2 > 0\}$. We thus conclude that $B_{r}(\mathbf{0}) \cap \bbR^2_+ \subset  \left\{ \psi -Wx_2 > 0 \right\} = \supp \omega$. Since the $C^{1}$-norm of $\psi$ in $\overline{\mathbb{R}^{2}_{+}}$ is uniformly bounded for all $\omega\in \calS_{1,1}(p)$ by Lemma \ref{lem_regularity_1}, the choice of $r=r(p)$ is uniform for all $\omega\in \calS_{1,1}({p})$.
			\end{proof}

			\begin{remark}
				\label{rmk_optimal_regularity}
				The regularity result ${\omega} \in C^{m, \alpha}(\overline{\bbR^2_{+}})$ for $1<p< \infty$ in Lemma \ref{lem_regularity_1} is optimal in the sense that ${\omega} \notin C^{m, \alpha+\varepsilon}(\overline{\bbR^2_{+}})$ for any $\varepsilon > 0$.
				Indeed, for $\omega \in \calS_{1,1}(p)$, $\psi=\mathcal{G}[\omega] \in C^2(\overline{\mathbb{R}^2_+})$ satisfies $\psi(x_1, 0) = 0$ for $x_1 \in \mathbb{R}$ and  
				\[
				\psi(0, x_2)-Wx_2 = (\partial_{x_2}\psi(0,0)-W)x_2 + o(x_2) \quad \textrm{as}\,\, x_2\to0.
				\]
				This means that 
				\[
				\kappa^{p-1} \omega(0, x_2) = \left( (\partial_{x_2}\psi(0,0) - W)x_2+ o(x_2)\right)^{\frac{1}{p - 1}}_+ 
				\ \simeq \ (\partial_{x_2}\psi(0,0) - W)^{\frac{1}{p - 1}} \, x_2^{\frac{1}{p - 1}},
				\]
				and ${\omega}\notin C^{m, \alpha+ \varepsilon}(\overline{\bbR^2_{+}})$ for any $\varepsilon > 0$. 
			\end{remark}

			\section{Compactness of Maximizing Sequences}\label{sec:existence}
			
			We prove existence of maximizers $\omega \in \calS_{\mu, \sigma}(p)$ for all $1<p\leq \infty$ and $0<\mu, \sigma < \infty$. We apply Lion's concentration--compactness lemma \cite{Lions1984}; see also \cite[Lemma 1]{BNL2013}, \cite[Lemma 4.1]{AC2022}. 
			\begin{lemma}
				\label{lem:Lions lemma}
				Let $0<\mu<\infty$.
				Let $\set{\rho_n}\subset L^1(\bbR^2_+)$ satisfy
				\begin{align*}
					\rho_n\geq0 \quad \mbox{ for }n\geq 1,\quad \int_{\bbR^2_+}\rho_n \dd \bfx=\mu_n \to \mu \quad \textrm{as}\,\, n\to\infty.
				\end{align*}
				Then, there exists a subsequence $\set{\rho_{n_k}}$ satisfying one of the following: \\
				
				\noindent(i) \textbf{Compactness.} There exists a sequence $\set{\bfy_k}\subset\overline{\bbR^2_+}$ such that $\rho_{n_k}(\cdot+\bfy_k)$ is tight, i.e. for arbitrary $\varepsilon>0$ there exists $R>0$ such that 
				$$
				\liminf_{k\to\infty}\int_{|\bfx-\bfy_k|<R}\rho_{n_k} \dd \bfx \geq \mu-\varepsilon.$$
				
				\noindent(ii) \textbf{Vanishing.} For each $R>0$, 
				$$
				\lim_{k\to\infty}\left[\sup_{\bfy\in\bbR^2_+}\int_{|\bfx-\bfy|<R}\rho_{n_k} \dd \bfx\right]=0.
				$$
				
				\noindent(iii) \textbf{Dichotomy.} There exists $\bar{\mu}\in (0,\mu)$ such that for arbitrary $\varepsilon>0$ there exist $k_0\geq1$ and $\set{\rho_{k, 1}},\set{\rho_{k, 2}}\subset L^1(\bbR^2_+)$ such that $\supp{\rho_{k, 1}}\cap\supp{\rho_{k, 2}}=\0$ and $0\leq \rho_{k, 1},\rho_{k, 2} \leq \rho_{n_k}$ for $k\geq k_0$, $\dist(\supp{\rho_{k, 1}},\supp{\rho_{k, 2}})\to\infty$ and 
				\begin{align*}
					\limsup_{k\to\infty}\left\{ \nrm{\rho_{n_k}-(\rho_{k, 1}+\rho_{k, 2})}_1+\big|\nrm{\rho_{k, 1}}_1-\bar{\mu} \big|+\big| \nrm{\rho_{k, 2}}_1-(\mu-\bar{\mu})\big|\right\}\leq\varepsilon.
				\end{align*}
			\end{lemma}

			\begin{theorem}\label{thm_existence_new}
				Let $1<p\leq \infty$ and $0<\mu, \sigma <\infty$. For any non-negative sequence $\{\omega_n\}_{n}$ on $\bbR^{2}_{+}$  satisfying 
				\begin{equation}\label{eq_stab_energy}
					\begin{aligned}
						{\liminf_{n \to \infty} E(\omega_n) \ge \calI_{\mu, \sigma}(p), \quad
							\limsup_{n \to \infty}\|x_2 \omega_n\|_1} \le \mu, \quad
						\limsup_{n \to \infty} \|\omega_n\|_p \le \sigma, 
					\end{aligned}
				\end{equation}
				there exists $\omega \in \calS_{\mu, \sigma}(p)$ such that, by a suitable translation for the $x_1$-variable, $\{\omega_{n} \}_n$ subsequently converges to $\omega$ with the norm $||x_2 (\cdot)||_{1}$ and 
				
				\begin{equation}
					\left\{
					\begin{aligned}
						\omega_{n}&\to \omega\quad \textrm{in}\ L^{p}(\mathbb{R}^{2}_{+})\quad \textrm{ for }\ p \in (1, \infty),\\
						\omega_{n} &\rightharpoonup^* \omega \quad \textrm{in}\ L^{\infty}(\mathbb{R}^{2}_{+})\quad \textrm{ for }\ p=\infty.
					\end{aligned}
					\right.
				\end{equation}
			\end{theorem}

			\begin{proof}
				We reduce to the case $\mu=\sigma=1$ by the scaling \eqref{eq: scaling}. By Remark \ref{r: maxmusigma} and choosing a subsequence, we may assume that   
				\begin{equation}
					\label{eq_stab_energy-red}
					{\lim_{n \to \infty} E(\omega_n) = \calI_{1,1}(p), \qquad \lim_{n \to \infty}     \|x_2 \omega_n\|_1 = 1, \qquad \lim_{n \to \infty} \|\omega_n\|_p = 1.}
				\end{equation} 
				By applying Lemma~\ref{lem:Lions lemma} to $\rho_n=x_2\omega_n$, up to a subsequence, one of the following three alternatives -- (i) \textit{Compactness}, (ii) \textit{Vanishing}, or (iii) \textit{Dichotomy} -- occurs. We first exclude the latter two cases.\\

				\noindent \textbf{(ii) Vanishing}. Suppose that vanishing occurs. Namely,  
				\begin{align*}
					\lim_{n\to\infty} \left( \sup_{\bfy\in \mathbb{R}^{2}_{+}}\int_{B_R(\bfy)\cap \mathbb{R}^{2}_{+}}x_2\omega_n \, \dd \bfx \right)  =   0 \quad \text{for each } R>0.
				\end{align*}
				We show that $\lim_{n\to\infty}E(\omg_n)=0$ and derive a contradiction $0<\calI_{1, 1}(p)=\lim_{n\to\infty}E[\omega_n]=0$. For $R > 1$, we decompose the energy as
				\begin{align*}
					2E(\omega_n)  & 
					= \left[ \iint_{\{|\bfx - \bfy|\geq R\}} + \iint_{\{|\bfx - \bfy|< R\}} \right] G(\bfx, \bfy)\omega_n(\bfx)\omega_n(\bfy) \, \dd \bfx\dd \bfy .
				\end{align*}
				By using the bound \eqref{eq: GE}, we estimate the first term as
				\begin{align*}
					\iint_{\{|\bfx-\bfy|\geq R\}} G(\bfx,\bfy)\omega_n(\bfx)\omega_n(\bfy)  \dd \bfx\dd \bfy \lesssim \frac{1}{ R^{2}}.
				\end{align*} 
				We split the second term into
				\begin{align*}
					\left[\iint_{\{{|\bfx - \bfy|< R, \, G\geq R x_2 y_2}\}}
					+\iint_{\{{|\bfx - \bfy|< R,\, G < R x_2 y_2}\}}\right]  G(\bfx, \bfy)\omega_n(\bfx)\omega_n(\bfy) \,  \dd \bfx\dd \bfy  =: \mathbf{I}_{n, R} + \mathbf{II}_{n, R}.
				\end{align*}
				For the points $\bfx$ and $\bfy$ satisfying $R x_2 y_2  \leq  G(\bfx, \bfy)  \leq  x_2 y_2 |\bfx - \bfy|^{-2}$,  we have $|\bfx - \bfy| \leq R^{-1/2}$. We set $0<\alpha=\frac{p-r}{r(p-1)}<1$ for $1<r<p$ and observe that $x_2^{\alpha}\omega_n$ is bounded in $L^{r}(\mathbb{R}^{2}_{+})$ by \eqref{eq: IE} of Proposition \ref{prop: IE}. We choose $r$ so that $\alpha q=\frac{p-r}{2(r-1)(p-1)}<1$ and $|x|^{-2\alpha}$ is integrable near the origin. For the conjugate exponent $r'$ to $r$, we set $1/r'=1/q+1/r-1$ and apply the Green's function estimate \eqref{eq: GE} and Young's convolution inequality to estimate
				\begin{align*}
					\mathbf{I}_{n, R}  
					&\leq  \iint_{\{|\bfx - \bfy|< R^{-1/2}\}} G(\bfx, \bfy) \omega_n(\bfx) \omega_n(\bfy)  \dd \bfx \dd \bfy \\
					&\lesssim 
					\int_{\mathbb{R}^{2}_{+}}\dd\bfx\int_{\mathbb{R}^{2}_{+}}
					\frac{1}{|\bfx-\bfy|^{2\alpha}}\mathbf{1}_{B_{R^{-1/2}}(\mathbf{0})}(\bfx-\bfy)x_2^{\alpha}\omega_n(\bfx)y_2^{\alpha}\omega_n(\bfy)\dd\bfy \\
					&\leq \left\|\frac{1}{|\bfx|^{2\alpha}}\mathbf{1}_{B_{R^{-1/2}}(\mathbf{0})}*(x_2^{\alpha}\omega_n)\right\|_{L^{r'}}||x_2^{\alpha}\omega_n||_{L^{r}} \\
					&\leq \left\|\frac{1}{|\bfx|^{2\alpha}}1_{B_{R^{-1/2}}(\mathbf{0})}\right\|_{L^{q}}||x_2^{\alpha}\omega_n||_{L^{r}}^{2} 
					\lesssim_{p} \frac{1}{R^{\frac{1}{q}-\alpha}}.
				\end{align*}
				We estimate 
				\begin{align*}
					\mathbf{II}_{n, R}  
					&\le R\iint\limits_{\{|\bfx - \bfy|< R\}} x_2 \omega_n(\bfx) \, y_2\omega_n(\bfy) \, \dd \bfx\dd \bfy \le R||x_2\omega_{n}||_{1} \left(\sup_{\bfy\in \mathbb{R}^{2}_{+}}\int_{B_R(\bfy)\cap \mathbb{R}^{2}_{+}}x_2\omega_n(\bfx)  \dd \bfx\right). 
				\end{align*}
				Letting $n\to\infty$ and then $R\to\infty$ imply $\lim_{n\to\infty}E(\omega_n)=0$. Thus, vanishing does not occur.\\
				
				\noindent
				\textbf{(iii) Dichotomy}.
				Suppose that dichotomy occurs for some $\bar{\mu} \in (0, \mu)$. Namely, for an arbitrary $\varepsilon > 0$, there exist $\omega_{1, n}$ and $\omega_{2, n}$ such that $0\leq \omega_{1,n}, \omega_{2,n}\leq \omega_n$ and 
				\[
				\operatorname{supp} (\omega_{1,n}) \cap \operatorname{supp} (\omega_{2,n}) = \emptyset, \qquad d_n  :=  \operatorname{dist} \left( \operatorname{supp} (\omega_{1,n}), \operatorname{supp} (\omega_{2,n}) \right) \to \infty \, \mbox{ as }\, n\to\infty.
				\]
				\begin{align*}
					&\limsup_{n \to \infty} \left(\left\| x_2\tld\omega_{n} \right\|_1 + \left| \mu_{1,n} - \bar{\mu} \right| + \left| \mu_{2,n} - (1 - \bar{\mu}) \right| \right)  \leq  \varepsilon,
				\end{align*} 
				for $\tld{\omg}_{n}=\omg_n - (\omg_{1,n}+\omg_{2,n})$ and $\mu_{i,n} =||x_2 \omega_{i,n}||_{1}$. We may assume that $\sup_{n}\left\| x_2\tld\omega_{n} \right\|_1\leq 2\varepsilon$. By choosing a subsequence, $\mu_{i,n}\to \mu_{i}\neq 0$ for some $\mu_{i}$ satisfying $\mu_1+\mu_2=1$ as $n\to\infty$ and $\varepsilon\to 0$. We set
				\begin{align*}\label{eq_disjoint}
					\sigma_{i,n}=||\omega_{i,n}||_{p},\quad i=1,2.
				\end{align*} 
				The constant $\sigma_{i,n}$ satisfies $(\sigma_{1,n}^{p}+\sigma_{2,n}^{p})^{1/p}\leq \sigma_{n}$ for $1<p<\infty$ and $\sigma_{1,n}\vee \sigma_{2,n}\leq \sigma_{n}$ for $p=\infty$. By $\sigma_n\to1$, we may assume that $\sigma_{i,n}\to \sigma_i$ for some $0\leq \sigma_i\leq 1$. We set $E(\omega_n) = E(\omega_{1, n}) + E(\omega_{2, n}) + \mathbf{I}_{n} + \mathbf{J}_{n}$ for
				\begin{equation*}
					\begin{split}
						\bfI_n = \iint G(\bfx,\bfy) \omega_{1,n}(\bfx)\omega_{2,n}(\bfy) \, \dd \bfx\dd \bfy, \qquad \bfJ_n= \iint G(\bfx,\bfy) (\omega_{n}(\bfx) - \tld{\omg}_{n}(\bfx))  \tld{\omg}_{n}(\bfy) \, \dd \bfx\dd \bfy .
					\end{split}
				\end{equation*}
				By using the bound \eqref{eq: GE}, we estimate 
				\begin{align*}
					\mathbf{I}_{n}  &=  \iint_{\{|\bfx - \bfy|\geq d_n\}} G(\bfx, \bfy)\omega_{1,n}(\bfx)\omega_{2,n}(\bfy)  \dd \bfx\dd \bfy \le  \frac{4}{ d_n^{2}}.
				\end{align*}
				By applying the inequality \eqref{eq_0704_01} and using the uniform bound for $\tld{\omg}_{n}$, 
				\begin{align*}
					\mathbf{J}_{n}  &\lesssim    \|x_2\tld{\omg}_{n}\|_{1}^{\frac{2p-2}{3p-2}} \lesssim \varepsilon^{\frac{2p-2}{3p-2}}.
				\end{align*}
				By $E(\omega_{i,n})\leq \calI_{\mu_{i,n},\sigma_{i,n}}(p)$, 
				\begin{equation*}
					\label{eq_dicho_1}
					\begin{aligned}
						E(\omega_n)\leq \calI_{\mu_{1,n},\sigma_{1,n}}(p)+\calI_{\mu_{2,n},\sigma_{2,n}}(p)+ \frac{4}{ d_n^{2}} + C_{p}\, \varepsilon^{\frac{2p-2}{3p-2}}.
					\end{aligned}
				\end{equation*}
				Letting $n\to\infty$ and $\varepsilon\to 0$ yield $\calI_{1,1}(p)\leq \calI_{\mu_{1},\sigma_{1}}(p)+\calI_{\mu_{2},\sigma_{2}}(p)$. This contradicts the strict superadditivity $\calI_{\mu_1,\sigma_{1}}(p)+\calI_{\mu_2,\sigma_{2}}(p)<\calI_{1,1}(p)$ in \eqref{eq: SSA} of Lemma \ref{lem_superadd}. Thus, dichotomy does not occur. We thus obtain:  \\
				
				\noindent
				\textbf{(i) Compactness}. There exists a sequence $\{\bfy_n = (y_{1, n}, y_{2, n})\}\subset \overline{\mathbb{R}^{2}_{+}}$ such that for arbitrary $\varepsilon>0$, there exists $R = R(\varepsilon)>0$ satisfying
				\begin{align*}
					\liminf_{n \to \infty} \int_{B_R(\bfy_n)\cap \mathbb{R}^{2}_{+}}x_2\omega_n  \dd \bfx  \geq  1 - \varepsilon.
				\end{align*}
				We may assume that $\sup_{n}\int_{\mathbb{R}^{2}_{+}\backslash B_R(\bfy_n)}x_2\omega_n  \dd \bfx\leq 2\varepsilon$.
				
				Suppose that $\limsup_{n\to\infty} y_{2,n}=\infty$. By choosing a subsequence, $y_{2,n} \ge 2R$ for all $n$ and $y_{2,n} \to \infty$. We decompose the energy as 
				\begin{align*}
					2E(\omega_n)  &=  \int_{\mathbb{R}^{2}_{+}} \calG[\omega_n]\, \omega_n \dd \bfx  =  \int_{B_R(\bfy_n)\cap \mathbb{R}^{2}_{+}} \calG[\omega_n]\, \omega_n   \dd \bfx + \int_{\mathbb{R}^{2}_{+} \setminus B_R(\bfy_n)} \calG[\omega_n] \, \omega_n  \dd \bfx.
				\end{align*}
				By applying \eqref{eq250718_2} and \eqref{eq_0704_01},
				\begin{align*}
					\int_{B_R(\bfy_n)\cap \mathbb{R}^{2}_{+}}\calG[\omega_n] \,\omega_n  \dd \bfx 
					&\leq  \left\|\calG[\omega_n]\right\|_{\infty} \int_{B(\bfy_n,R)\cap \mathbb{R}^{2}_{+}} \frac{1}{x_2} \cdot x_2\omega_n  \dd \bfx \lesssim \frac{1}{y_{2,n} - R}, \\
					\int_{\mathbb{R}^{2}_{+} \setminus B_R(\bfy_n)} \calG[\omega_n] \, \omega_n  \dd \bfx
					&\lesssim   
					\left(\int_{\mathbb{R}^{2}_{+} \setminus B_R(\bfy_n)} x_2 \omega_n  \dd \bfx\right)^{\frac{2p-2}{3p-2}}  \lesssim \varepsilon^{\frac{2p-2}{3p-2}}.
				\end{align*}        
				Letting $n\to\infty$ and $\varepsilon\to 0$ imply $\lim_{n\to\infty}E(\omega_n)= 0$. This yields a  contradiction $0<\calI_{1,1}(p)=\lim_{n\to\infty}E(\omega_n)=0$. Thus, $\{y_{2,n}\}$ is bounded.
				
				We may assume that $y_{n}=0$ by translating $\omega_n$ for the $x_1$-variable and replacing $R(\varepsilon)$ by a larger constant. Since $\{\omega_n\}_{n}$ is uniformly bounded in $L^{p}$, by choosing a subsequence,
				\begin{equation*}
					\left\{
					\begin{aligned}
						&\omega_n  \rightharpoonup  \omega \quad \text{in } L^{p}\quad \textrm{for}\ 1<p<\infty, \\
						&\omg_{n} \rightharpoonup^* \omg \quad \text{in } L^{\infty}\quad \textrm{for}\ p=\infty.
					\end{aligned}
					\right.
				\end{equation*}
				By the lower semi-continuity of the $L^{p}$-norm for those convergence, $||\omega||_{p}\leq 1$ for $1<p\leq \infty$ and 
				\begin{equation*}
					\label{impulse_limit}
					\int_{ \bbR^2_+ \setminus B_{\varep}} x_2  \omega   \dd \bfx  \le   2\varepsilon,
				\end{equation*}
				for $B_{\varep}=B_{R(\varepsilon)}(\mathbf{0})$. Since $\varepsilon > 0$ is arbitrary, $\|x_2 \omega\|_1=1$ and $\omega \in \calK_{1,1}(p)$. We set 
				\begin{align*}
					2E(\omega_n)  = \left[ \int_{B_{\varep} \cap \mathbb{R}^{2}_{+}} + \int_{\mathbb{R}^{2}_{+} \setminus B_{\varep}}\right] \,  \mathcal{G}[\omega_n] \omega_n \, \dd \bfx   =: \mathbf{I} + \mathbf{II}.
				\end{align*}
				We decompose $\mathbf{I}$ as
				\begin{align*}
					\mathbf{I}=    \int_{B_{\varep} \cap \mathbb{R}^{2}_{+}} \omega_n(\bfx) \dd \bfx \int_{\mathbb{R}^{2}_{+}} G(\bfx,\bfy) \omega_n(\bfy)  \dd \bfy 
					&=   \int_{B_{\varep} \cap \mathbb{R}^{2}_{+}} \int_{B_{\varep} \cap \mathbb{R}^{2}_{+}} + \int_{B_{\varep} \cap \mathbb{R}^{2}_{+}} \int_{\mathbb{R}^{2}_{+} \setminus B_{\varep}} =: \mathbf{I}_1 + \mathbf{I}_{2}.
				\end{align*}
				Using the symmetry $G(\bfx, \bfy) = G(\bfy, \bfx)$, we see that
				\begin{align*}
					\mathbf{I}_{2} &= \int_{B_{\varep} \cap \mathbb{R}^{2}_{+}} \omega_n(\bfx)  \dd \bfx \int_{\mathbb{R}^{2}_{+} \setminus B_{\varep}} G(\bfx,\bfy) \omega_n(\bfy)  \dd \bfy =  \int_{B_{\varep} \cap \mathbb{R}^{2}_{+}} \omega_n(\bfy)  \dd \bfy \int_{\mathbb{R}^{2}_{+} \setminus B_{\varep}} G(\bfx,\bfy) \omega_n(\bfx)  \dd \bfx \\
					&\leq  \int_{\mathbb{R}^{2}_{+} \setminus B_{\varep}} \omega_n(\bfx)  \dd \bfx \int_{ \mathbb{R}^{2}_{+}} G(\bfx,\bfy) \omega_n(\bfy)  \dd \bfy  =   \int_{\mathbb{R}^{2}_{+} \setminus B_{\varep}} \mathcal{G}[\omega_n](\bfx) \omega_n(\bfx)  \dd \bfx = \mathbf{II}.
				\end{align*}
				Consequently, $2E(\omega_n) = \mathbf{I}_{1} + \mathbf{I}_{2} + \mathbf{II} \le \mathbf{I}_1 + 2\mathbf{II},$
				and we apply \eqref{eq_0704_01} and estimate
				\begin{align*}
					0\le 2E(\omega_n) - \mathbf{I}_{1}
					=  2\mathbf{II} \lesssim   
					\left(\int_{\mathbb{R}^{2}_{+} \setminus B_{\varep}}x_2\omega_n  \dd \bfx\right)^{\frac{2p-2}{3p-2}}\lesssim  \varepsilon^{\frac{2p-2}{3p-2}}.
				\end{align*}
				By applying a similar argument for $E(\omega)$ and combining the estimates for $E(\omega_n)$ and $E(\omega)$,
				\begin{align*}
					2\left| E(\omega_n) - E(\omega) \right| \, &\leq \, \left|\int_{B_{\varepsilon} \cap \mathbb{R}^{2}_{+}} \int_{B_{\varepsilon} \cap \mathbb{R}^{2}_{+}} G(\bfx, \bfy) \left(\omega_n(\bfx)\omega_n(\bfy) - \omega(\bfx)\omega(\bfy)\right)   \dd \bfx \dd \bfy \right| +  C\varepsilon^{\frac{2p-2}{3p-2}}.
				\end{align*}
				Since $G(\bfx,\bfy) \in L^{p'}(B_{\varep} \times B_{\varep})$ for the conjugate exponent $p' \in [1, \infty)$ to $p$, and $\omega_n(\bfx)\omega_n(\bfy) \rightharpoonup \omega(\bfx)\omega(\bfy)$ in $L^{p}(B_{\varep} \times B_{\varep})$ ($\rightharpoonup^*$ when $p = \infty$), we take the limits $n \to \infty$ and $\varepsilon \to 0$ and conclude that $\lim_{n\to \infty}E(\omega_n)
				=  E(\omega)$.

				This convergence implies $E(\omega) = \calI_{1,1}(p)$ and $\omega\in \calS_{1,1}(p)$. By Remark \ref{r: maxmusigma}, $||\omega||_{p}=1$ and $||x_2\omega||_{1}=1$. In particular, 
				\[
				\lim_{n \to \infty} \|\omega_n\|_p = \|\omega\|_p.
				\]
				By the uniform convexity of $L^{p}$ for $1<p<\infty$, $\omega_n \to \omega$ in $L^{p}$. By using  
				\begin{align*}
					\int_{\mathbb{R}^{2}_{+}}x_2|\omega_n - \omega|  \dd \bfx  =
					\int_{B_{\varep}\cap \mathbb{R}^{2}_{+}}x_2|\omega_n - \omega|  \dd \bfx
					+ \int_{\mathbb{R}^{2}_{+}\setminus B_{\varep}}x_2|\omega_n - \omega|  \dd \bfx  \leq  \int_{B_{\varep}\cap \mathbb{R}^{2}_{+}}x_2|\omega_n - \omega|  \dd \bfx+ 2\varepsilon,
				\end{align*}
				we obtain $x_2\omega_n \to x_2\omega$ in $L^1$. For $p=\infty$, $\omega_n  \rightharpoonup^* \omega$ in $L^{\infty}$ and 
				\[
				\sup_{n}\int_{ \bbR^2_+ \setminus B_{\varep}} x_2  \omega_n   \dd \bfx  \le   2\varepsilon   \quad \mbox{ and } \quad			\int_{ \bbR^2_+ \setminus B_{\varep}} x_2  \omega   \dd \bfx  \le   2\varepsilon.
				\]
				Since $\omega\in \calS_{1,1}({\infty})$, $\omega = \mathbf{1}_{\Omega}$ for a bounded set $\Omega\subset \mathbb{R}^{2}_{+}$. We set $B_{\varepsilon}=(B_{\varepsilon}\cap \Omega) \cup (B_{\varepsilon}\cap \Omega^{c})$ so that $|\omega_n-\omega|=1-\omega_n$ on $B_{\varepsilon}\cap \Omega$ and $|\omega_n-\omega|=\omega_n$ on $B_{\varepsilon}\cap \Omega^{c}$. Then, by using 
				\begin{equation*}
					\label{eq_conv_impulse}
					\begin{aligned}
						\int_{\bbR^2_+} x_2 |\omega_n - \omega| \, \dd\bfx &= \int_{B_{\varep}\cap \mathbb{R}^{2}_{+}}x_2|\omega_n - \omega|  \dd \bfx
						+ \int_{\mathbb{R}^{2}_{+}\setminus B_{\varep}}x_2|\omega_n-\omega|  \dd \bfx \\
						& \le \int_{B_{\varep}\cap \Omega }x_2 (\omega - \omega_n)  \dd \bfx + \int_{B_{\varep}\cap \Omega^c} x_2 \omega_n   \dd \bfx + 2\varepsilon
					\end{aligned}
				\end{equation*}
				letting $n\to\infty$ and $\varepsilon\to 0$ imply $x_2\omega_n \to x_2 \omega$ in $L^1$. The proof is now complete.
			\end{proof}

			\section{Existence and Stability of Sadovskii Vortices}
			\label{sec:main}
			
			We deduce the existence and stability of Sadovskii vortices from the variational principle \eqref{eq_maximization} for $1<p\leq \infty$. We also establish the $L^{1}$ stability of minimal-mass Sadovskii vortices and deduce a stability estimate for the horizontal center of mass.

			\subsection{Sadovskii Vortices}\label{sec:finite_p}
			
			We first obtain existence and stability of Sadovskii vortices from that of the set of maximizers $\calS_{\mu,\sigma}(p)$ for $1< p \leq \infty$.
            For $p = \infty$, we obtain stability in the impulse norm { under a smallness condition on the excess vorticity $(\theta_0 - \sigma)_+$.}

			\begin{theorem}[Theorem~\ref{thm_intro_main}]
				\label{thm_main_1}
				Let $1<p \leq \infty$ and $0<\mu,\sigma<\infty$. The following holds:
				
				\begin{enumerate}[itemsep=1em, topsep=1em]
					\item[(i)] \textbf{(Existence)} The odd extension $\tld{\omega}_p$ of each $\omega_p\in \calS_{\mu,\sigma}(p)$ is a Sadovskii vortex travelling with the same speed $W = W(p, \mu, \sigma) > 0$ and 
					\begin{equation}\label{eq_main_supp}
						B_{r}(\mathbf{0}) \subset \supp \tld{\omega}_p,
					\end{equation}
					with the uniform constant $ r = r(p, \mu, \sigma) >0$ by a suitable translation for the $x_1$-variable. 
					\item[(ii)] \textbf{(Regularity)} For $p=\infty$, $\tld{\omega}_{\infty}$ is a Sadovskii vortex patch. For $p \in (1, \infty)$, 
					\begin{equation}\label{eq_main_regularity}
						\tld{\omega}_p \in C^{m, \alpha}(\bbR^2), \quad \textrm{ for } \, m =m_p = \left\lceil \frac{1}{p-1} - 1 \right\rceil,\ \alpha = \alpha_p =\frac{1}{p-1}-m,
					\end{equation}
					and $\tld{\omega}_p \notin C^{m, \alpha + \varepsilon}(\bbR^2)$ for any $\varepsilon > 0$.
					\item[(iii)] \textbf{(Stability)} For each $p \in (1, \infty)$ and $\varepsilon>0$, there exists $\delta = \delta(p, {\mu, \sigma,\, }   \varepsilon)>0$ such that for any initial data $0 \le \theta_0 \in  L_c^{\infty}(\bbR^2_+)$ satisfying  
					\begin{align}
						\inf_{\omega\in \calS_{\mu,\sigma}(p)}\left\{  \left\| \theta_0-\omega \right\|_{p} 
						+\left\| x_2 \left( \theta_0-\omega \right) \right\|_{1} \right\}   \label{eq: initial}
						\leq  \delta,
					\end{align}
					the unique global-in-time solution $\theta(t)$ to 
                    {\eqref{eq: half Euler eq.}}
                    satisfies
					\begin{align} 
						\inf_{\omega\in \calS_{\mu,\sigma}(p)}\left\{   \left\| \theta \left( t \right)-\omega \right\|_{p}
						+\left\|x_2 \left( \theta \left(t \right)-\omega  \right) \right\|_{1}\right\}
						\leq  \varepsilon \quad \mbox{for all } \,\, t \in \bbR.   \label{eq: stability}
					\end{align} 
                   
                    For $p = \infty$, $\varep >0$, and each fixed $q \in (1, \infty)$, there exists $\delta = \delta(\mu, \sigma, \varepsilon, q) > 0$ such that for any initial data $0 \le \theta_0 \in  L_c^{\infty}(\bbR^2_+)$ satisfying
                    \[
                     {\| (\theta_0 - \sigma)_+ \|_q}  +   \inf_{\omega \in \calS_{\mu, \sigma}(\infty)} \| x_2 (\theta_0 - \omega ) \|_1 \le \delta,
                    \]
                    the corresponding solution $\theta(t)$ to 
                    {\eqref{eq: half Euler eq.}}
                    satisfies
                    \[
                    \inf_{\omega \in \calS_{\mu, \sigma}(\infty)} \| x_2(\theta(t) - \omega) \|_1 \le \varepsilon \quad \textrm{for all } \,\, t \in \bbR.
                    \]
                    
				\end{enumerate}
			\end{theorem}
			
			\begin{proof}
				The existence of maximizers $\omega_p\in \calS_{\mu,\sigma}(p)$ follows from Theorem \ref{thm_existence_new}.
                By the Euler--Lagrange equation \eqref{eq: EL} and Lemma~\ref{lem_pohozaev}, together with two-parameter scaling \eqref{two_scaling}, any maximizer $\omega_p \in \calS_{\mu,\sigma}(p)$ is a traveling wave with traveling speed $W = W(p, \mu, \sigma)$.
				The support of $\omega_p$ is compact by Lemma~\ref{lem_cpt_support}, and in contiguous contact with the $x_1$-axis by Lemma \ref{prop_touching}, with touching scale $r = r(p, \mu, \sigma)$ via the scaling \eqref{two_scaling}.
                Thus, $\tld{\omega}_p$ is the Sadovskii vortex.
                The regularity of $\tld{\omega}_p$ follows from Lemma \ref{lem_regularity_1} and Remarks~\ref{rmk_reg_odd} and \ref{rmk_optimal_regularity}.
				
				It remains to show the stability. For $1<p<\infty$, suppose that the claim were false. Then there exist $\varepsilon_0>0$, non-negative $\{\theta_{0, n} \}_{n = 1}^{\infty}\subset L^{\infty}_{c}(\mathbb{R}^{2}_{+})$, and $\{t_n\}_{n = 1}^{\infty}\subset \mathbb{R}$ such that for each $n \ge 1$, 
				\begin{align*}\label{eq241228_2}
					&\inf_{\omega \in \calS_{\mu,\sigma}(p)} \left\{  \left\| \theta_{0, n} - \omega \right\|_{p} + \left\| x_2 (\theta_{0, n} - \omega) \right\|_{1} \right\}   \le  \frac{1}{n},\\
					&\inf_{\omega \in \calS_{\mu,\sigma}(p)}\left\{  \left\| \theta_{n} \left(t_n, \cdot \right) - \omg \left( \cdot  \right) \right\|_{p} + \left\| x_2 \big( \theta_{n} \left(t_n, \cdot \right) - \omega \left( \cdot  \right) \big) \, \right\|_{1}  \right\}  \ge  \varepsilon_0>0,
				\end{align*}
				for the global-in-time weak solution $ \theta_n(t)$ to \eqref{eq: Euler eq.} satisfying $\theta_n(0) = \theta_{0, n}$. We take a sequence $\{\omega_n\}_{n} \subset \calS_{\mu,\sigma}(p)$ such that
				\[
				\| x_2 (\theta_{0, n} - \omega_n) \|_{1} + \| \theta_{0, n} - \omega_n\|_{p}  \to  0,
				\]
				and observe from \eqref{eq_0704_03} that  
				\[
				\lim_{n \to \infty}E(\theta_{0, n}) = \calI_{\mu,\sigma}(p),\quad \limsup_{n\to\infty}||x_2\theta_{0, n}||_{1}\leq \mu, \quad \limsup_{n\to\infty}||\theta_{0, n} ||_{p}\leq \sigma.
				\]
				By the conservation property of the Yudovich solution, $\set{\theta_{n}(t_n)}_{n=1}^\infty$ also satisfies the same property. By applying Theorem~\ref{thm_existence_new} to the sequence $\{ \theta_n(t_n)  \}_{n}$, we obtain a subsequence still denoted by $\{ \theta_n(t_n)  \}_{n}$ such that $\theta_n(t_n)$ converges to some ${\omega} \in \calS_{\mu,\sigma}(p)$ with the norm $||\cdot ||_{p}+||x_2(\cdot)||_{1}$ by a suitable translation for the $x_1$-variable. By letting $n \to \infty$,
				\begin{align*}
					0 &=  \lim_{n \to \infty} \left[	\left\| x_2 \left( \theta_n (t_n) - {\omega} \right) \right\|_{1} +    \left\|  \theta_n (t_n) - {\omega} \right\|_{p} \right] \\
					&\ge  \liminf_{n \to \infty} \left\{ \inf_{\omega \in \calS_{\mu,\sigma}(p)} \left(   \left\| x_2 \left( \theta_{n} \left(t_n\right) - \omega \right)\, \right\|_{1} + \left\| \theta_{n} \left(t_n \right) - \omega \right\|_{p} \right) \right\}  \ge  \varepsilon_0>0.
				\end{align*}
				We obtain a contradiction.
				
                {
                For $p = \infty$, fix $q \in (1, \infty)$ and suppose that the claim were false.
                Then there exist $\varepsilon_0 > 0$, non-negative $\theta_{0,n} \in L^\infty_c(\mathbb{R}^2_+)$, times $t_n \in \bbR$, and $\omega_n \in \calS_{\mu,\sigma}(\infty)$ such that
                \[
                \|x_2(\theta_{0,n} - \omega_n)\|_1 + \|(\theta_{0,n} - \sigma)_+\|_q \to 0,
                \]
                whereas the corresponding solutions satisfy
                \[
                \inf_{\omega \in \calS_{\mu,\sigma}(\infty)} \|x_2(\theta_n(t_n) - \omega)\|_1 \ge \varepsilon_0.
                \]
                Set
                \[
                v_n(t) = \theta_n(t) \wedge \sigma, \quad \eta_n(t) = (\theta_n(t) - \sigma)_+.
                \]
                Conservation of impulse and of the vorticity distribution gives
                \[
                \|x_2\theta_n(t)\|_1 = \|x_2\theta_{0,n}\|_1 \to \mu, \quad \|\eta_n(t)\|_q = \|\eta_n(0)\|_q \to 0.
                \]
                Since $0 \le v_n(t) \le \sigma$ and $0 \le \eta_n(t) \le \theta_n(t)$, \eqref{eq_0704_02} yields, uniformly in $t$,
                \[
                E(v_n(t)) \le C, \quad E(\eta_n(t)) \le C_q \|x_2\theta_{0,n}\|_1^{\frac{4q-4}{3q-2}} \|\eta_n(0)\|_q^{\frac{2q}{3q-2}} \to 0.
                \]
                By \eqref{eq_L2energy} and the Cauchy--Schwarz inequality,
                \[
                0 \le E(\theta_n(t)) - E(v_n(t)) \le 2\sqrt{E(v_n(t)) E(\eta_n(t))} + E(\eta_n(t))\to 0
                \]
                uniformly in $t$.
                Since $0\le \omega_n \le \sigma$,
                \[
                \|x_2(v_n(0) - \omega_n)\|_1 \le \|x_2(\theta_{0,n} - \omega_n)\|_1 \to 0.
                \]
                Then, by applying \eqref{eq_0704_03} with $p=\infty$ we have $E(v_n(0)) \to \calI_{\mu,\sigma}(\infty)$.
                Hence energy conservation and the preceding estimates imply
                \[
                E(v_n(t_n))\to\calI_{\mu,\sigma}(\infty), \quad \limsup_{n\to\infty}\|x_2v_n(t_n)\|_1 \le \mu, \quad \|v_n(t_n)\|_\infty \le \sigma.
                \]
                By Theorem~\ref{thm_existence_new}, after passing to a subsequence and applying horizontal translations,
                \[
                \|x_2(v_n(t_n) - \omega) \|_1 \to 0
                \]
                for some $\omega \in \calS_{\mu,\sigma}(\infty)$.
                We apply the same translations to $\theta_n(t_n)$ and $\eta_n(t_n)$ and use the same notation.
                In particular, $\|x_2 v_n(t_n)\|_1 \to \mu$, and therefore
                \[
                \|x_2\eta_n(t_n)\|_1 = \|x_2\theta_{0,n}\|_1 - \|x_2 v_n(t_n)\|_1 \to 0.
                \]
                Consequently,
                \[
                \|x_2(\theta_n(t_n)-\omega)\|_1 \le \|x_2(v_n(t_n)-\omega)\|_1 + \|x_2\eta_n(t_n)\|_1 \to 0.
                \]
                Since $\calS_{\mu,\sigma}(\infty)$ is invariant under horizontal translations, this contradicts the choice of $t_n$.
                }
			\end{proof}

			\subsection{Minimal-Mass Vortices}\label{sec:patch_stability}
We {refine the stability} of Sadovskii vortex patches by restricting to the minimal-mass maximizers in $\calS_{\mu, \sigma}(\infty)$. 
            More generally, for $\calS_{\mu,\sigma}(p)$ and $1<p\leq \infty$, we consider the minimal mass
			\begin{align*}
				0<\nu(p, \mu, \sigma) := \inf_{\omega\in \calS_{\mu,\sigma}(p) }\|\omega\|_1<\infty \quad \left( \nu \left(p, 1, 1 \right) = \nu \left(p\right) \right),
			\end{align*}
			by using the bound \eqref{eq_unif_mass} and define the set of minimal-mass maximizers:
			\begin{equation}\label{def_min_mass}
			    \calS'_{\mu,\sigma}(p)
			:=\left\{
			\omega\in \calS_{\mu,\sigma}(p):\
			\|\omega\|_1=\nu(p, \mu, \sigma)
			\right\}.
			\end{equation}

			\begin{proposition}
				For any $p \in (1, \infty]$,
				$\calS'_{\mu,\sigma}(p)
				\neq \emptyset
				$.
			\end{proposition}
			\begin{proof}
				For a sequence $\set{\omg_n}\subset \calS_{\mu,\sigma}(p)$ satisfying $\nrm{\omg_n}_1\to\nu(p, \mu, \sigma)$, there exists a subsequence converging to a limit $0\neq \omg\in \calS_{\mu,\sigma}(p)$ with the norm $||x_2(\cdot)||_{1}$ by a suitable translation for the $x_1$-variable by Theorem~\ref{thm_existence_new}. In particular, by choosing a subsequence and Fatou's lemma, $\omega_n\to \omega$ a.e. and  
                \[
                \|\omega\|_1 \le \liminf_{n \to \infty} \|\omega_n \|_1 = \nu(p, \mu, \sigma).
                \]
                Thus, $0\neq \omega\in \calS'_{\mu,\sigma}(p)$.
			\end{proof}

            \begin{theorem}\label{thm_stability of patch}
            For each $\mu, \sigma > 0$ and $p \in (1, \infty]$, the set $\calS'_{\mu, \sigma}(p)$ of minimal-mass maximizers is orbitally stable in the following sense.
            \begin{itemize}[itemsep=0.8em]
            \item[(i)] If $1 < p < \infty$, then for every $\varepsilon > 0$, there exists $\delta = \delta(p, \mu, \sigma, \varepsilon) > 0$ such that, for any non-negative $\theta_0\in L^\infty_c(\bbR^2_+)$ satisfying
            \begin{equation}\label{eq_patch_assumption}
            \inf_{\omega\in\calS'_{\mu,\sigma}(p)} \left\{ \|\theta_0 - \omega\|_{L^1\cap L^p} + \|x_2(\theta_0-\omega)\|_1 \right\}\leq\delta,
            \end{equation}
            the unique global-in-time solution $\theta(t)$ to \eqref{eq: half Euler eq.} with initial data $\theta_0$ satisfies
            \begin{equation}\label{eq_patch_min}
            \inf_{\omega\in\calS'_{\mu,\sigma}(p)} \left\{ \|\theta(t)-\omega\|_{L^1\cap L^p} + \|x_2(\theta(t)-\omega)\|_1 \right\} \le \varepsilon \quad \mbox{for all }\,\, t\in\bbR.
            \end{equation}

             {\item[(ii)] If $p = \infty$, then for every $q \in (1,\infty)$ and $\varepsilon > 0$, there exists $\delta=\delta(q,\mu,\sigma,\varepsilon)>0$ such that, for any non-negative $\theta_0\in L^\infty_c(\bbR^2_+)$ satisfying
            \begin{equation}\label{eq_patch_assumption_infty}
            \inf_{\omega\in\calS'_{\mu,\sigma}(\infty)} \left\{ \|\theta_0-\omega\|_{L^1\cap L^q} + \|x_2(\theta_0-\omega)\|_1 \right\}\leq\delta,
            \end{equation}
            the unique global-in-time solution $\theta(t)$ to \eqref{eq: half Euler eq.} with initial data $\theta_0$ satisfies
            \begin{equation}\label{eq_patch_min_infty}
            \inf_{\omega\in\calS'_{\mu,\sigma}(\infty)} \left\{ \|\theta(t)-\omega\|_{L^1\cap L^q} + \|x_2(\theta(t)-\omega)\|_1 \right\} \le \varepsilon \quad\mbox{for all }\,\, t\in\bbR.
            \end{equation}}
            \end{itemize}
            \end{theorem}

			\begin{proof}
            {We set $q = p$ when $1 < p < \infty$ while fix any $q \in (1, \infty)$ when $p = \infty$.}
				We apply the same contradiction argument as in the proof of Theorem~\ref{thm_main_1} (iii). Suppose that the claim were false. Then, there exists $\varepsilon_0>0$, $\set{t_n}\subset\bbR$, and non-negative $\{\tht_{0,n}\}\subset  L_c^\ift(\mathbb{R}^{2}_{+})$ satisfying 
				\begin{align*}
					&\inf_{\omg \in \calS'_{\mu,\sigma}(p)} \left\{ \nrm{\tht_{0,n}-\omg}_{L^{1}\cap L^{q}} +\nrm{x_2(\tht_{0,n}-\omg)}_{L^{1}}\right\}<\frac{1}{n}, \\
					&\inf_{\omg \in \calS'_{\mu,\sigma}(p)}\left\{\nrm{\tht_n(t_n)-\omg}_{L^{1}\cap L^{q}} + \nrm{x_2(\tht_n(t_n)-\omg )}_{L^{1}}\right\}\geq \varepsilon_0>0,
				\end{align*}
				for the solution $\tht_n(t)$ to 
                \eqref{eq: half Euler eq.}
                satisfying $\theta_n(0) =\tht_{0,n}$. We first consider the case $p\in (1,\infty)$. By using a sequence $\set{\omg_n}\subset \calS'_{\mu,\sigma}(p)$ satisfying
				\begin{equation*}
					\nrm{\tht_{0,n}-\omg_n}_{L^{1}\cap L^{q}} + \nrm{x_2(\tht_{0,n}-\omg_n)}_{L^{1}}\to 0, 
				\end{equation*}
				and \eqref{eq_0704_03}, we observe that the sequence $\{\theta_{0,n}\}$ satisfies 
				\begin{align*}
					\lim_{n\to\infty}E(\theta_{0,n})= \calI_{\mu,\sigma}(p),\quad 
					\limsup_{n\to\infty}||x_2\theta_{0,n}||_{1}\leq \mu,\quad 
					\limsup_{n\to\infty}||\theta_{0,n}||_{p}\leq \sigma.
				\end{align*}
				By conservations, $\theta_{n}(t_n)$ also satisfies the same {bounds} and we apply Theorem~\ref{thm_existence_new} to extract a subsequence converging to a limit $\omega\in \calS_{\mu,\sigma}(p)$ with the norm $||\cdot ||_{p}+||x_2(\cdot)||_{1}$ by a suitable translation for the $x_1$-variable.
				By
				\[
				||\omega||_1\leq \liminf_{n\to\infty}||\theta_n(t_n)||_{1}=\liminf_{n\to\infty}||\theta_{0,n}||_{1}= \nu(p, \mu, \sigma),
				\]
				it follows that $\omega\in \calS'_{\mu,\sigma}(p)$.
				By Lemma \ref{lem_cpt_support}, we take a bounded set $\Omega\subset \overline{\mathbb{R}^{2}_{+}}$ such that $\omega=\omega1_{\Omega}$ and observe that for $q\in [1,\infty)$,  
				\begin{equation}
					\begin{aligned}
						||\theta_n(t_n)-\omega||_{q}^{q}
						&=\int_{\Omega}| \theta_n(t_n)-\omega|^{q} \dd \bfx +\int_{\Omega^{c}} \theta_n(t_n)^{q} \dd \bfx\\
						&=\int_{\Omega}| \theta_n(t_n)-\omega|^{q} \dd \bfx +\int_{\mathbb{R}^{2}_{+}} \theta_n(t_n)^{q} \dd \bfx -\int_{\Omega} \theta_n(t_n)^{q} \dd \bfx \\
						&=\int_{\Omega}| \theta_n(t_n)-\omega|^{q} \dd \bfx +\int_{\mathbb{R}^{2}_{+}} \theta_{0,n}^{q} \dd \bfx - \int_{\Omega} \theta_n(t_n)^{q} \dd \bfx.
					\end{aligned}
					\label{eq: q}
				\end{equation}
				By applying \eqref{eq: q} with exponent $1$, we observe that the first term on the right-hand side vanishes since $\theta_{n}(t_n)\to \omega$ in $L^{p}$.
				By using $||\theta_{0,n}||_1\leq 1/n+\nu(p, \mu, \sigma)$ and $\nu(p, \mu, \sigma)=||\omega||_1=\int_{\Omega}\omega \dd \bfx \leq \liminf_{n\to\infty}\int_{\Omega}\theta_n(t_n) \dd \bfx$, it follows that $\theta_n(t_n)\to \omega$ in $L^{1}$.
				Thus, 
				\begin{align*}
					0= \lim_{n\to\infty}\left\{\nrm{\tht_n(t_n)-\omg}_{L^{1}\cap L^{q}} + \nrm{x_2(\tht_n(t_n)-\omg )}_{L^{1}}\right\}\geq \varepsilon_0>0,
				\end{align*}
				and we obtain a contradiction.
                
                {
                It remains to consider $p = \infty$ and $q \in (1, \infty)$.
                Choose $\omega_n \in \calS'_{\mu,\sigma}(\infty)$ such that
                \[
                \|\theta_{0,n} - \omega_n\|_{L^1 \cap L^q} + \|x_2( \theta_{0,n} - \omega_n)\|_1 \to 0.
                \]
                Since $0 \le \omega_n \le \sigma$,
                \[
                \|(\theta_{0,n}-\sigma)_+\|_q \le \|\theta_{0,n} - \omega_n\|_q \to 0.
                \]
                The compactness argument in the proof of Theorem~\ref{thm_main_1}~(iii) implies, after passing to a subsequence and applying horizontal translations,
                \[
                \|x_2(\theta_n(t_n) - \omega)\|_1 \to 0, \quad \theta_n(t_n) \to \omega \quad\textrm{a.e.},
                \]
                for some $\omega\in\calS_{\mu,\sigma}(\infty)$.
                We keep the same notation for the translated functions.

                Set $\nu = \nu(\infty,\mu,\sigma)$.
                Since $\|\omega_n\|_q^q = \sigma^{q-1}\nu$, the initial $L^q$ convergence and conservation yield
                \[
                \|\theta_n(t_n)\|_q^q = \|\theta_{0,n}\|_q^q \to \sigma^{q-1} \nu.
                \]
                Since $\omega$ is a patch of height $\sigma$, 
                \[
                \sigma^{q-1}\nu \le \sigma^{q-1}\|\omega\|_1 = \|\omega\|_q^q \le \liminf_{n\to\infty}\|\theta_n(t_n)\|_q^q = \sigma^{q-1}\nu.
                \]
                Thus $\omega\in\calS'_{\mu,\sigma}(\infty)$.
                We write $\omega = \sigma\mathbf{1}_{\Omega}$, where $\Omega$ is bounded.
                By conservation property,
                \[
                \|(\theta_n(t_n)-\sigma)_+\|_q = \|(\theta_{0,n}-\sigma)_+\|_q \to 0.
                \]
                Moreover, $(\sigma - \theta_n(t_n))_+ \to 0$ a.e. on $\Omega$ and is bounded by $\sigma$.
                Then
                \[
                \int_{\Omega}|\theta_n(t_n)-\sigma|^q \dd\bfx \le \|(\theta_n(t_n)-\sigma)_+\|_q^q + \int_{\Omega}(\sigma-\theta_n(t_n))_+^q\dd\bfx \to 0.
                \]
                Together with the convergence of $L^q$ norms, \eqref{eq: q} now implies $\theta_n(t_n)\to\omega$ in $L^q$.

                The initial $L^1$ convergence with the conservation property gives 
                \[
                \|\theta_n(t_n)\|_1 = \|\theta_{0,n}\|_1 \to \nu = \|\omega\|_1.
                \]
                By H\"older's inequality,
                \[
                \int_\Omega |\theta_n(t_n) - \omega| \dd\bfx \le |\Omega|^{1-1/q} \left( \int_\Omega |\theta_n(t_n)-\omega|^q\dd\bfx \right)^{1/q} \to 0.
                \]
                Using \eqref{eq: q} with exponent $1$, we obtain
                \[
                \begin{aligned}
                \|\theta_n(t_n) - \omega\|_1 & = \int_\Omega|\theta_n(t_n)-\omega| \dd\bfx + \|\theta_{0,n}\|_1 - \int_\Omega\theta_n(t_n)\,\dd\bfx\\
                & \to  \|\omega\|_1 - \int_\Omega\omega \,\dd\bfx = 0.
                \end{aligned}
                \]
                This yields the required contradiction.
                }
			\end{proof}

            {
            \begin{remark}
            Stability in the norm $\|\cdot\|_1+\|x_2(\cdot)\|_1$ also holds in Theorem~\ref{thm_stability of patch}~(ii) when $0 \le \theta_0 \le \sigma$.
            Indeed, for every $\omega\in\calS'_{\mu,\sigma}(\infty)$,
            \[
            \|\theta_0 - \omega\|_2^2    \le \sigma \|\theta_0-\omega\|_1,
            \]
            so this follows from Theorem~\ref{thm_stability of patch}~(ii) with $q=2$.
            \end{remark}
            }

			\begin{remark}\label{rmk_minmass_p}
				For $p=2$, solutions to the Euler--Lagrange equation \eqref{eq_ell} are unique up to translation for the $x_1$-variable \cite{Bur1996} (see also \cite{AC2022}, \cite{ACJ2025}). Namely, all maximizers $\omega\in \calS_{\mu,\sigma}(2)$ are of the form $\omega_{L}=\lambda (\psi_L-Wx_2)_{+}$ (Chaplygin--Lamb dipole) for $\psi_L=\Psi_{L}+Wx_2$ and 
				\begin{equation*}
					\Psi_{L}(x)=\left\{
					\begin{aligned}
						C_LJ_1(\sqrt{\lambda} r)\sin\theta,\quad r\leq a,\\
						-W\left(r-\frac{a^{2}}{r}\right)\sin\theta,\quad r>a,
					\end{aligned}
					\right. \label{eq: Lamb}
				\end{equation*}
				in the polar coordinates $(r,\theta)$ with the constants
				\begin{align*}
					C_L=-\frac{2W}{\sqrt{\lambda}J_0(c_0)},\quad a=\frac{c_0}{\sqrt{\lambda}},   \label{eq: LambConst}
				\end{align*}
				where $J_{m}(r)$ is the $m$-th order Bessel function of the first kind and $c_0=3.8317\cdots$ is the first zero point of $J_1$, i.e. $J_1(c_0)=0$. Its kinetic energy, impulse, and enstrophy are as follows:
				\begin{align*}
					\frac{1}{2}\int_{\mathbb{R}^{2}_{+}}|\nabla \psi_{L}|^{2}\dd\bfx=\frac{c_0^{2}\pi W^{2}}{\lambda},\quad
					\int_{\mathbb{R}^{2}_{+}}x_2 \omega_{L}\dd\bfx=\frac{c_0^{2}\pi W}{\lambda},\quad
					\int_{\mathbb{R}^{2}_{+}}| \omega_{L}|^{2}\dd\bfx=c_0^{2}\pi W^{2}. 
				\end{align*}
				The Lagrange multipliers $W$ and $\lambda$ are given by $W=\frac{\sigma}{c_0\sqrt{\pi}}$ and $\lambda =\frac{c_0\sqrt{\pi}\sigma}{\mu}$.
				The kinetic energy is $\calI_{\mu,\sigma}(2)=\calI_{1,1}(2)\mu \sigma$ for $\calI_{1,1}(2)=\frac{1}{c_0\sqrt{\pi}}$ and mass is as follows:
				\begin{align*}
					\int_{\mathbb{R}^{2}_{+}} \omega_{L}\dd\bfx=2C_L\int_{0}^{c_0}rJ_1(r)dr=\nu(2, \mu, \sigma).
				\end{align*}
				The set of minimal mass maximizers agrees with the set of maximizers and is characterized as 
				\begin{align*}
					\calS'_{\mu,\sigma}(2)=\calS_{\mu,\sigma}(2)=\{\omega_{L}(\cdot +y \mathbf{e}_1)\ |\ y\in \mathbb{R}\}.
				\end{align*}
				For $p\neq 2$, the set $\calS'_{\mu,\sigma}(p)$ can be strictly smaller than $\calS_{\mu,\sigma}(p)$.
			\end{remark}

			\subsection{Horizontal Displacement}
			
			For a function $\theta(t)$, we set the horizontal center-of-mass  displacement
			\begin{equation}
				P(t) := \frac{1}{ \| \theta(t) \|_1  }\int_{\bbR^2_+}x_1\tht(t)\dd\bfx-\frac{1}{ \| \theta(0) \|_1  }\int_{\bbR^2_+}x_1\tht(0)\dd\bfx.   \label{eq: shift}
			\end{equation}\\
			If $\theta(t)$ is a non-negative solution to the Euler equations \eqref{eq: Euler eq.}, 
			\begin{align*}
				\frac{\dd }{\dd t} P(t) =\frac{1}{ \| \theta_0 \|_1}\int_{\bbR^2_+}u_1(t)\tht(t)\dd\bfx
				&=\frac{1}{ \| \theta_0 \|_1}\int_{\bbR^2_+}\int_{\bbR^2_+}\frac{1}{2\pi}\left[\frac{-x_2+y_2}{|\bfx-\bfy|^2}+\frac{x_2+y_2}{|\bfx-\bfy^*|^2}\right]\tht(t,\bfy)\tht(t,\bfx)\dd\bfy\dd\bfx\\
				&=\frac{1}{ \| \theta_0 \|_1 }\int_{\bbR^2_+}\int_{\bbR^2_+}\frac{1}{2\pi}\frac{x_2+y_2}{|\bfx-\bfy^*|^2} \tht(t,\bfy)\tht(t,\bfx) \dd\bfy\dd\bfx>0. 
			\end{align*}
			In particular, $\frac{\dd }{\dd t} P(t) =W$ for the traveling wave solution $\theta(t)=\omega(x_1-Wt,x_2)$ by the formula \eqref{formula_W}. We deduce the stability of $P(t)$ from the $L^{1}$ stability result (Theorem \ref{thm_stability of patch}).

			\begin{lemma}\label{lem_shift estimate}
				Let $1<p\leq \infty$ and $0 < \mu, \sigma < \infty$.
				Let $W_p = W(p, \mu, \sigma)>0$ be the traveling speed of Sadovskii vortices in $\mathcal{S}_{\mu, \sigma}(p)$.
				For an arbitrary $\varepsilon>0$ and the solution $\theta(t)$ in Theorem \ref{thm_stability of patch}, the horizontal displacement \eqref{eq: shift} satisfies 
				\begin{equation}\label{tech eq: shift estimate} 
					\left| \frac{\dd }{\dd t} P(t) -W_p  \right|\le C(p, \mu, \sigma) \,  \varepsilon (1 + \|\theta_0\|_{\infty}) \quad \mbox{ for all }\, t\in\bbR.
				\end{equation}
			\end{lemma}
			
			\begin{proof}
				We take an arbitrary $\omega\in \calS'_{\mu, \sigma}(p) $ and represent the speed $W_{p} $ by the formula \eqref{formula_W} in Remark~\ref{rmk_W_formula}.
                By subtracting $W_{p}$ from $\frac{\dd}{\dd t} P(t)$, 
				\begin{align*}
					\frac{\dd P(t)}{\dd t}-W_{p}
					&=\frac{1}{2\pi ||\theta_0||_{1}} \int_{\bbR^2_+}\int_{\bbR^2_+} \frac{x_2+y_2}{|\bfx-\bfy^*|^2} \left(\tht(t,\bfy)\tht(t,\bfx)-\omega(\bfy)\omega(\bfx)  \right)\dd\bfy\dd\bfx\\
					&+\frac{1}{2\pi}\left(\frac{1}{||\theta_0||_{1}}-\frac{1}{||\omega||_{1}} \right)\int_{\bbR^2_+}\int_{\bbR^2_+}\frac{x_2+y_2}{|\bfx-\bfy^*|^2} \omega(\bfy)\omega(\bfx) \dd\bfy\dd\bfx.
				\end{align*}
{
Note that
\[
\sup_{\bfx\in\bbR^2_+} \int_{\bbR^2_+} \frac{x_2+y_2}{|\bfx-\bfy^*|^2} |f(\bfy)| \,\dd\bfy \lesssim \|f\|_{L^1\cap L^\infty}.
\]
Together with conservation of the $L^1\cap L^\infty$-norm and
\[
\bigl|\|\theta_0\|_1-\|\omega\|_1\bigr| = \bigl|\|\theta(t)\|_1-\|\omega\|_1\bigr| \le \|\theta(t)-\omega\|_1,
\]
this yields
    \[
    \left|\frac{\dd P(t)}{\dd t}-W_p\right| \lesssim \frac{\|\theta_0\|_{L^1\cap L^\infty} + \|\omega\|_{L^1\cap L^\infty}} {\|\theta_0\|_1} \|\theta(t)-\omega\|_1 .
    \]
}
				For $p=\infty$, $||\omega||_{\infty}=\sigma$.
				For $1<p<\infty$, $||\omega||_{\infty}\leq \kappa^{-\frac{1}{p-1}}||\mathcal{G}[\omega] ||_{\infty}^{\frac{1}{p-1}}$ and $||\mathcal{G}[\omega] ||_{\infty}$ is bounded by a  constant depending on $p$, $\mu$, and $\sigma$ by \eqref{eq250718_2}.
				Thus, $||\omega||_{L^{\infty}\cap L^{1}} \lesssim_{p, \mu, \sigma} 1 $.
				By taking an infimum for $\omega\in \calS_{\mu, \sigma}'(p)$ and applying the stability estimate \eqref{eq_patch_min} for $\delta\leq \varepsilon$, 
				\begin{align*}
					\left| \frac{\dd P(t)}{\dd t} -W_{p}\right|\lesssim \frac{1}{ ||\theta_0||_{1}} \left( 1 + ||\theta_0||_{L^{\infty}\cap L^{1}}  \right) \varepsilon.
				\end{align*} 
				Since $\|\theta_0\|_1 \simeq_{\delta} \|\omega\|_1=\nu(p, \mu, \sigma)$ by the initial assumption, we obtain \eqref{tech eq: shift estimate} with the constant $C=C\left(  p, \mu, \sigma\right)>0$.
			\end{proof}
			
			We conclude this section by proving Theorem~\ref{thm_intro_patch}.
			
			\begin{proof}[\textbf{Proof of Theorem~\ref{thm_intro_patch}}]
                {
                The assertion for general non-negative initial data follows from Theorem~\ref{thm_stability of patch}~(ii) with $\mu = \sigma = 1$.
                For the patch type initial data, take $q=2$ and $\theta_0 = \mathbf{1}_{\Omega_0}$.
                Since
                \[
                \|\mathbf{1}_{\Omega_0} - \mathbf{1}_A\|_2 = |\Omega_0\triangle A|^{1/2},
                \]
                the initial assumption implies \eqref{eq_patch_assumption_infty} in Theorem~\ref{thm_stability of patch}~(ii) after choosing $\delta$ sufficiently small.
                We apply  Theorem~\ref{thm_stability of patch}~(ii) with $\varepsilon/2$ in place of $\varepsilon$.
                }
				Let $\calS'(\infty) = \calS'_{1, 1}(\infty)$, and let $\theta(t) = \mathbf{1}_{\Omega(t)}$ be the solution to \eqref{eq: half Euler eq.} with initial data $\theta(0) = \theta_0$.
				Since
				\[
				\| \theta(t) - \omega\|_1 + \|x_2(\theta(t) - \omega)\|_1 = \int_{\Omega(t) \Delta A} (x_2 + 1) \, \dd \bfx, \qquad  t \in \bbR,
				\]
				for any $\omega = \mathbf{1}_A \in \calS'(\infty)$, Theorem~\ref{thm_intro_patch}-(\textit{i}) follows directly from \eqref{eq_patch_assumption_infty}--\eqref{eq_patch_min_infty} in Theorem~\ref{thm_stability of patch}~(ii).
				Finally, by integrating \eqref{tech eq: shift estimate} of Lemma~\ref{lem_shift estimate} ($p = \infty$) with respect to $t$, we obtain Theorem~\ref{thm_intro_patch}-(\textit{ii}).
			\end{proof}

			\appendix
			
			\section{Shift Function Estimate}\label{app_shift}
			In this section, we present a refined formulation of stability by introducing a shift function that measures the horizontal alignment between a perturbed solution and the family of Sadovskii vortices, leading to the quantitative shift estimate in Proposition~\ref{lem_shift estimate_1}.
			
			\begin{remark}
				For simplicity of exposition, we restrict ourselves to the patch case $p=\infty$.
				The same arguments apply to any $p\in(1,\infty)$ and yield the corresponding results for the minimal-mass class $\calS'_{\mu,\sigma}(p)$.
			\end{remark}

			Throughout this section, we work in the normalized setting $\mu = \sigma = 1$.
			To eliminate the horizontal translation invariance of maximizers, we fix Steiner symmetric representatives (see Remark~\ref{rmk_Steiner}) and introduce the \emph{centered} class
			\[
			\calS_{\mathrm{cen}}(\infty)
			:= \left\{
			\omega \in \calS'_{1,1}(\infty): \omega \,\text{ is Steiner symmetric }
			\right\}.
			\]
			
			We now restate the stability theorem in a form that makes explicit the horizontal translation of the perturbed solution.
			The formulation is obtained by extracting the translation that arises in the contradiction argument in the proof of Theorem~\ref{thm_stability of patch}.
			
			\begin{theorem}[Shift function]\label{thm_stability of patch_shift}
				For {each $q \in (1, \infty)$ and} $\varepsilon>0$, there exists $\dlt=\dlt( \varepsilon, q)>0$ such that the following holds:
				For any initial data $0\leq \tht_0\in L_c^\ift(\bbR^2_+)$ satisfying 
				\begin{equation}\label{eq_patch_assumption_shift}
					\inf_{\omg\in \calS_{\mathrm{cen}}(\ift)} \left\{ \| \theta_0 - \omega\|_{ {L^1 \cap L^q} }  + \| x_2(\tht_0-\omg)\|_1\right\} \le \dlt,
				\end{equation}
				there exists a shift function $\tau:\bbR\to\bbR$ with $\tau(0)=0$ such that the unique global-in-time solution $\tht(t)$ with $\tht(0)=\tht_0$ satisfies 
				\begin{equation}\label{eq_patch_min_shift}
					\inf_{\omg\in \calS_{\mathrm{cen}}(\ift)} \left\{ \| \theta(t) - \omega_{\tau(t)}\|_{{L^1 \cap L^q}  }  + \| x_2(\tht(t)-\omg_{\tau(t)})\|_1\right\} \leq \varepsilon
				\end{equation}
				for all $t \in \bbR$, where $\omega_{\tau(t)}:=\omega(\cdot-(\tau(t),0))$. 
			\end{theorem}

			Although the shift function $\tau(t)$ in \eqref{eq_patch_min_shift} is not uniquely determined and may depend on the choice of closest maximizer, its evolution is rigid: it remains close to the linear motion $W_\infty t$.

			\begin{proposition}[Shift estimate]\label{lem_shift estimate_1}
				There exists $\varepsilon_0>0$ such that, if $\varepsilon < \varepsilon_0$ in Theorem~\ref{thm_stability of patch_shift}, 
				then any shift function satisfying \eqref{eq_patch_min_shift} obeys
				\begin{equation}\label{tech eq: shift estimate_1} 
					|\tau(t)-W_\ift t|<C\,\varepsilon\,(1 + |t|) \quad \textrm{ for all }\, t\in\bbR,
				\end{equation}
				where $C=C\left(\nrm{\tht_0}_\ift, \sup_{\bfx\in\supp\tht_0}|x_1|\right)>0$ and $W_\ift>0$ is the traveling speed of Sadovskii patches in $\calS_{1,1}(\ift)$.
			\end{proposition}
			
			\begin{remark}\label{rmk: shift is optimal}
				The shift estimate \eqref{tech eq: shift estimate_1} is optimal in the sense that the power of $\veps$ on the right-hand side coincides with that in \eqref{eq_patch_min_shift}.
				Indeed, fix any $\bfone_\Omg\in\calS_{\mathrm{cen}}(\ift)$ and, for each $\veps>0$, consider the solution $\theta_\veps(t)$ with initial data 
				\[
				\theta_\veps(0)=(1+\veps)\bfone_\Omega \quad \textrm{ in }\, \bbR^2_+.
				\]
				Since the corresponding traveling speed is $(1+\veps)W_\ift$, the deviation of the shift from $W_{\infty}t$ is of order $\veps  t$.
				On the other hand, the stability estimate \eqref{eq_patch_min_shift} is of order $\varepsilon$.
				Therefore, the bound in \eqref{tech eq: shift estimate_1} cannot be improved in terms of the power of $\varepsilon$.
				
				This optimality stems from our use of the horizontal center of mass, a quantity not exploited in previous shift estimates for the Chaplygin--Lamb dipole \cite{CJ2022, JYZ}, the Kelvin $m$-waves \cite{CJ2022kelvin}, and the Hill's spherical vortex \cite{CJ2023}.
				To the best of the authors' knowledge, such an optimal shift estimate is established here for the first time. 
			\end{remark}
			
			\begin{remark}[Perimeter growth]\label{rmk: perimeter_growth}
				As a consequence of the shift estimate \eqref{tech eq: shift estimate_1}, one can construct perturbed patch-type initial data $\mathbf{1}_{\Omega(0)}$ near any Sadovskii maximizer $\omega\in\calS_{\mathrm{cen}}(\infty)$ whose diameter and perimeter grow linearly in time.
				
				Indeed, starting from a maximizer, one may attach an arbitrarily long and thin backward tail while keeping the perturbation small in the sense of \eqref{eq_patch_min_shift}.
				The shift estimate implies that the main body of the patch propagates with speed close to $W_\infty$, whereas the tip of the thin tail moves much more slowly.
				More precisely, one obtains
				\[
				\sup_{x,y\in\Omega(t)} |x_1-y_1|
				\ge (W_\infty - o(1))\,t,
				\]
				which means linear-in-time growth of both diameter and perimeter.
			\end{remark}

			We now turn to the proof of Proposition~\ref{lem_shift estimate_1}.
			The argument relies on uniform bounds for the set $\calS(p)$ of maximizers.
			Although we focus on the patch case $p=\infty$, these bounds hold for each $p\in(1,\infty]$ and are of independent interest.
			See also Remark~\ref{rmk_unif_bdd} and \cite[Lemma~4.6]{CJS2025}.
			
			\begin{lemma}[Uniform estimates]\label{lem_uniform}
				Let $p \in (1, \infty]$.
				There exists a constant $C = C(p)  > 1$ such that, for any maximizer $\omega \in \calS_{1, 1}(p)$, the following hold:
				\begin{enumerate}
					\item (Mass bound)
					\begin{equation}
						\label{eq_uni_mass}
						C^{-1} \le \|\omega\|_1 \le C .
					\end{equation}
					
					\item (Size of the compact support)
					\begin{equation}
						\label{eq_unif_supp}
						\sup_{\bfx \in \operatorname{supp} \overline{\omega}} |\bfx| \le C,
					\end{equation}
					where $\overline{\omega}(\bfx) := \omega(\bfx - (c,0))$ satisfies the Steiner symmetry condition (Remark~\ref{rmk_Steiner}) for some $c \in \mathbb{R}$.
				\end{enumerate}
			\end{lemma}
			
			\begin{proof}
				We fix $p \in (1, \infty]$ and,
                for notational simplicity, we suppress the dependence on $p$ throughout.
				Let $\omega \in \calS_{1, 1}$ be an arbitrary maximizer with traveling speed $W$ and set $\psi = \calG[\omega]$.
				
				\medskip

				\noindent \textit{(1) Mass bound.}
				Since the lower bound in \eqref{eq_uni_mass} is obtained in \eqref{unif_mass_low} of Remark~\ref{lem_unif}, we only prove the upper bound.
                {
                From the Euler--Lagrange equation \eqref{eq: EL}, we have $\sup_{\bfx \in \supp \omega} x_2 \le \|\psi\|_{\infty}/W$.
                By \eqref{eq250718_2} of Lemma~\ref{lem_stream_bdd_2} and \eqref{eq_pohozaev} of Lemma~\ref{lem_pohozaev}, we may fix sufficiently large $L > 0$ such that $\supp\omg\subset\set{0<x_2<L}$ and}
				 let $A_n := \{ \bfx \in \bbR^2_+:  2^{-(n+1)}L \le  x_2  < 2^{-n}L \}$.
				Then
				\begin{equation}
					\label{eq0923_2}
					\begin{aligned}
						\| \omega\|_1
						= \sum_{n = 0}^{\infty} \int_{ A_n } \omega(\bfx) \, \dd\bfx  
						\le  &\sum_{n = 0}^{\infty} \int_{ A_n } \omega(\bfx) \,  \frac{\psi(\bfx)}{Wx_2}\, \dd\bfx \\
						\lesssim &\sum_{n = 0}^{\infty} 2^n L^{-1} \int_{ A_n } \omega(\bfx) \, \psi(\bfx)\, \dd\bfx \\
						= &\sum_{n = 0}^{\infty} 2^n L^{-1}  \int_{\bbR^2_+}  \left( \int_{A_n }G(\bfx, \bfy) \, \omega(\bfx)  \, \dd\bfx  \right) \,\omega(\bfy) \, \dd\bfy.
					\end{aligned}
				\end{equation}
				Since $\|\omega\|_{\infty} \lesssim 1$, by \eqref{eq241003_1} of Proposition~\ref{lem0801_01},
				\begin{align*}
					\int_{A_n }G(\bfx, \bfy) \, \omega(\bfx)  \, \dd\bfx \lesssim  \int_{0 < x_2  < 2^{-n}L }G(\bfy, \bfx) \, \dd\bfx 
					\lesssim y_2^{1/2} (2^{-n}L)^{3/2}.
				\end{align*}
				Combining this with \eqref{eq0923_2} and applying H\"older's inequality yields
				\[
				\|\omega\|_1 \lesssim \sum_{n = 0}^{\infty} (2^{-n}L)^{1/2} \int_{\bbR^2_+} y_2^{1/2} \omega(\bfy) \, \dd\bfy \lesssim L^{1/2} \|y_2\omega\|_1^{1/2} \|\omega\|_1^{1/2} \lesssim \|\omega\|_1^{1/2}.
				\]
				Thus $\|\omega\|_1 \lesssim 1$, proving \eqref{eq_uni_mass}.
				
				\medskip
				
				\noindent \textit{(2) Size of the compact support.}
				We may assume $c = 0$.
				Since $\sup_{\bfx \in \supp \omega} x_2 \le \|\psi\|_{\infty}/W \lesssim 1$, it suffices to prove
				\begin{equation}
					\label{eq_supp_x_1}
					\sup_{\bfx \in \supp \omega} |x_1| \lesssim 1.
				\end{equation}
				For any $\bfx = (x_1, x_2) \in \bbR^2_+$, the mean value theorem in $x_2$ gives
				\begin{align*}
					\frac{\psi(\bfx)}{x_2} = -\frac{1}{2\pi} \int_{\bbR^2} \frac{ \overline{x}_2 - y_2}{| \overline{\bfx} - \bfy|^2}    \tld\omega(\bfy) \, \dd\bfy \lesssim \int_{\bbR^2_+} \frac{1}{|\overline{\bfx} - \bfy|} \omega(\bfy) \, \dd\bfy
				\end{align*}
				for some $\bar{\bfx} = (x_1, \overline{x}_2) \in \bbR^2_+$ with $\overline{x}_2 \in (0, x_2)$, where $\tld\omg$ is the odd-extension of $\omg$ to $\bbR^2$.
				Let $r > 0$ and, by splitting the integral into $\{| \overline{\bfx} - \bfy|\ge r\}$ and $\{| \overline{\bfx} - \bfy|<r\}$, we estimate
				\begin{align*}
					\frac{\psi(\bfx)}{x_2}  \lesssim   r^{-1} \|\omega\|_1
					+ r^{1/3} \|\omega\|_{L^3 \left([x_1-r, x_1+r]  \times (0,\infty) \right)} \lesssim r^{-1} \|\omega\|_1
					+ r^{1/3} \|\omega\|_{L^1 \left([x_1-r, x_1+r]  \times (0,\infty) \right)}^{1/3},
				\end{align*}
				where we used $\|\omega\|_{\infty} \lesssim 1$ in the last inequality.
				Since $\omega$ is nonincreasing in $x_1 > 0$, for $x_1 > 2r$, 
				\[
				\|\omega\|_{L^1 \left([x_1-r, x_1+r]  \times (0,\infty) \right)} = \int_0^{\infty} \left( \int_{x_1 -r}^{x_1 + r} \omega(s, x_2) \, \dd s \right) \dd x_2 \lesssim \frac{r}{x_1} \|\omega\|_1.
				\]
				Hence, for any $r > 0$ and any $\bfx \in \bbR^2_+$ such that $x_1 > 2r$, we have
				\[
				\frac{\psi(\bfx)}{x_2}  \lesssim r^{-1} \|\omega\|_1 + r^{1/3} \left( \frac{r}{x_1} \|\omega\|_1 \right)^{1/3} \lesssim r^{-1} + r^{2/3}x_1^{-1/3},
				\]
				with implicit constants independent of $x_2$.
				Therefore, we see that there exists $R > 0$ such that $\psi(\bfx)/x_2 < W$ whenever $x_1 > R$.
				Since $\omega$ is even in $x_1$ and $\supp \omega = \{  \psi(\bfx) \ge Wx_2 \}$, this proves \eqref{eq_supp_x_1}.
			\end{proof}

			We are now ready to prove Proposition~\ref{lem_shift estimate_1}.
			
			\begin{proof}[Proof of Proposition~\ref{lem_shift estimate_1}]
				Let $\tau:\bbR\to\bbR$ be a shift function with $\tau(0)=0$ satisfying \eqref{eq_patch_min_shift}.
				For each $t \in \bbR$, there exists $\omega = \omega^t \in \calS_{\mathrm{cen}}(\infty)$ such that
				\begin{equation}
					\label{eq_fix_time}
					\| \theta(t) - \omega_{\tau(t)}\|_1 \leq 2\varepsilon.
				\end{equation}
                Whenever necessary, we may further decrease $\varepsilon$.
				We may assume $t > 0$, since the argument for $t < 0$ is identical.
				Recall that
				\[
				P(t) = \frac{1}{ \| \theta(t) \|_1  }\int_{\bbR^2_+}x_1\tht(t)\dd\bfx - \frac{1}{ \| \theta(0) \|_1  }\int_{\bbR^2_+}x_1\tht(0)\dd\bfx.
				\]
				By \eqref{tech eq: shift estimate} of Lemma~\ref{lem_shift estimate} with $P(0)=0$, we have
				\begin{equation}\label{eq_per_com}
					|P(t)-W_\ift t|\lesssim  (1 + \|\theta_0\|_{\infty})\varepsilon\, t.
				\end{equation}
				Let 
				\begin{equation}\label{eq_per_im}
					Q(t):=  \frac{1}{ \| \theta(t) \|_1  }\int_{\bbR^2_+}x_1\tht(t)\dd\bfx = P(t) + Q(0)
				\end{equation}
				and
				\begin{equation}\label{eq_width}
					R(t):=\sup_{\bfx\in\supp\tht(t)}|x_1|.
				\end{equation}
				
				\medskip
				
				\noindent We first claim that
				\begin{equation}
					\label{eq_per_com2}
					|\tau(t) - Q(t)| \lesssim\varepsilon(1 +  R(t)).
				\end{equation}
				Assuming \eqref{eq_per_com2}, combining it with \eqref{eq_per_com} and \eqref{eq_per_im} yields
                \begin{equation}\label{eq_tau_2}
				\begin{aligned}
					|\tau(t)-W_\ift t| \le   &  |\tau(t)-Q(t)| + |Q(t)-W_\ift t|  \\
					\le &|\tau(t)-Q(t)| + |P(t)-W_\ift t| + |\tau(0) - Q(0)| \\
					\lesssim &\,\varepsilon ( 1 + R(0) + R(t)) + \varepsilon(1+\nrm{\tht_0}_\ift)\, t,
				\end{aligned}
                \end{equation}
				where we used $\tau(0) = 0$ and \eqref{eq_fix_time} with $t=0$.
				Note that $R(0) = \sup_{\bfx \in \supp \theta_0}|x_1|$ and
				$$
				|R(t)-R(0)| \leq t \cdot\sup_{s\in\bbR} \| \nabla^{\perp}(-\Delta_{\bbR^2_+})^{-1}\theta(s)\|_{\infty} \lesssim t\,\nrm{\tht_0}_{L^1\cap L^\ift}
				\lesssim t\,(1+\nrm{\tht_0}_{\ift}).
				$$
				This together with \eqref{eq_tau_2} implies \eqref{tech eq: shift estimate_1}.
				
				\medskip
				
				\noindent It remains to prove \eqref{eq_per_com2}.
				By \eqref{eq_unif_supp} of Lemma~\ref{lem_uniform}, there exists a uniform constant $C_0>0$ such that
				\[
				\supp\omg_{\tau(t)}\subset\calB(t) :=\set{\bfx\in\bbR^2_+:|\bfx-\tau(t)\bfe_1|<C_0}.
				\]
				Then
				$$
				\|\theta(t)\|_1 \cdot Q(t) = \int_{\bbR^2_+} x_1\tht(t)\dd\bfx = \left[ \int_{\calB(t)}  + \int_{\calB(t)^c} \right]x_1\tht(t)\dd\bfx,
				$$
				and we observe that, by definition of $R(t)$ in \eqref{eq_width},
				\begin{equation}
					\label{eq_com_tail}
					\left|\int_{\calB(t)^c} x_1\tht(t)\dd\bfx \right|\leq R(t)\int_{\calB(t)^c}\tht(t)\dd \bfx =  R(t)\int_{\calB(t)^c}|\tht(t)-\omg_{\tau(t)}|\dd\bfx \lesssim \varepsilon \, R(t).
				\end{equation}
				Since $\omega \in \calS_{\mathrm{cen}}(\infty)$ is even in $x_1$, it is clear that
				$$
				\int_{\calB(t)} \tau(t)\omega_{\tau(t)}\, \dd \bfx = \| \omega\|_1 \, \tau(t) = \int_{\bbR^2_+}x_1\omg_{\tau(t)}\dd\bfx=\int_{\calB(t)}x_1\omg_{\tau(t)}\dd\bfx,
				$$
				and then, from \eqref{eq_fix_time} and \eqref{eq_com_tail}, we see that
				\begin{align*}
					\left| \|\theta(t)\|_1 \,Q(t) - \|\theta(t)\|_{L^1(\calB(t))}\,\tau(t) \right| \le &\left| \int_{\calB(t)^c} x_1 \theta(t) \, \dd\bfx \right| + \left| \int_{\calB(t)} (x_1 - \tau(t)) \theta(t) \, \dd\bfx \right| \\
					\lesssim & \,\varepsilon R(t) + \left| \int_{\calB(t)} (x_1 - \tau(t)) \left( \theta(t) - \omega_{\tau(t)} \right) \, \dd\bfx \right| \\
					\lesssim &\, \varepsilon(R(t) + C_0).
				\end{align*}
				Since $\|\theta(t)\|_{L^1(\calB(t))} \ge \|\theta(t)\|_1 - 2\varepsilon = \|\theta(0)\|_1 - 2\varepsilon$, we have
				\begin{align*}
					\left|Q(t) - \tau(t) \right| \le &\left| \left(1 - \frac{\|\theta(t) \|_1}{\| \theta(t) \|_{L^1(\calB(t))}} \right)Q(t) \right| +  \left|   \frac{\|\theta(t) \|_1}{\| \theta(t) \|_{L^1(\calB(t))}} Q(t) - \tau(t)\right| \\
					& \lesssim \, \varepsilon  |Q(t)| + \frac{\varepsilon}{\| \theta(t) \|_{L^1(\calB(t))}} \,(R(t) + C_0) \\
					& \lesssim \, \varepsilon (1 + |Q(t)| + R(t)).
				\end{align*}
				This with the obvious relation $|Q(t)|\leq R(t)$ implies \eqref{eq_per_com2}.
				This completes the proof of Proposition~\ref{lem_shift estimate_1}.
			\end{proof}

			\section*{Acknowledgments}
			{The authors would like to thank the anonymous referees for their constructive comments and valuable suggestions, which have greatly improved the manuscript.}   KA was supported by the JSPS through the Grant in Aid for Scientific Research (C) 24K06800, MEXT Promotion of Distinctive Joint Research Center Program JPMXP0723833165, and Osaka Metropolitan University Strategic Research Promotion Project (Development of International Research Hubs).
            KC and YS were supported by the National Research Foundation of Korea (NRF) grant funded by the Korean government (MSIT) (No. RS-2023-00274499 and No. RS-2022-NR070754).
			IJ was supported by the NRF grant from the Korea government (MSIT), No. 2022R1C1C1011051, RS-2024-00406821, and the Asian Young Scientist Fellowship.
			KW was supported by the NRF grant funded by the Korean government (MSIT) (No. RS-2022-NR070754) and by the Samsung Science and Technology Foundation under Project Number SSTF-BA2002-04.
			
			\bibliographystyle{alpha}
			\bibliography{Sadovskii}

@article{HHR2024,
  archiveprefix = {arXiv},
  author = {Z. Hassainia and T. Hmidi and E. Roulley},
  title = {Desingularization of time-periodic vortex motion in bounded domains via {KAM} tools},
  journal = {\href{https://arxiv.org/abs/2408.16671}{arXiv:2408.16671}},
  year = {2024},
  note = {to appear in J. Eur. Math. Soc.}
}

@article{ACJ2025,
    archiveprefix = {arXiv},
	author = {K. Abe and K. Choi and I.-J. Jeong},
	journal = {\href{https://arxiv.org/abs/2510.00539}{arXiv:2510.00539}},
	title = {{S}tability of {L}amb dipoles for odd-symmetric and non-negative initial disturbances without the finite mass condition},
	year = {2025}}

@article {Abe21_measure,
    AUTHOR = {Abe, K.},
     TITLE = {The vorticity equations in a half plane with measures as
              initial data},
   JOURNAL = {Ann. Inst. H. Poincar\'e{} C Anal. Non Lin\'eaire},
  FJOURNAL = {Annales de l'Institut Henri Poincar\'e{} C. Analyse Non
              Lin\'eaire},
    VOLUME = {38},
      YEAR = {2021},
    NUMBER = {4},
     PAGES = {1055--1094},
      ISSN = {0294-1449,1873-1430},
   MRCLASS = {35Q35 (35K90)},
  MRNUMBER = {4266235},
       DOI = {10.1016/j.anihpc.2020.10.002},
       URL = {https://doi.org/10.1016/j.anihpc.2020.10.002},
}

@article {GallagherGallay,
    AUTHOR = {Gallagher, I. and Gallay, T.},
     TITLE = {Uniqueness for the two-dimensional {N}avier-{S}tokes equation
              with a measure as initial vorticity},
   JOURNAL = {Math. Ann.},
  FJOURNAL = {Mathematische Annalen},
    VOLUME = {332},
      YEAR = {2005},
    NUMBER = {2},
     PAGES = {287--327},
      ISSN = {0025-5831,1432-1807},
   MRCLASS = {35Q30 (35A05 76D03 76D05)},
  MRNUMBER = {2178064},
MRREVIEWER = {Toshiaki\ Hishida},
       DOI = {10.1007/s00208-004-0627-x},
       URL = {https://doi.org/10.1007/s00208-004-0627-x},
}

@article {Delort,
    AUTHOR = {Delort, J.-M.},
     TITLE = {Existence de nappes de tourbillon en dimension deux},
   JOURNAL = {J. Amer. Math. Soc.},
  FJOURNAL = {Journal of the American Mathematical Society},
    VOLUME = {4},
      YEAR = {1991},
    NUMBER = {3},
     PAGES = {553--586},
      ISSN = {0894-0347,1088-6834},
   MRCLASS = {76C05 (35Q30)},
  MRNUMBER = {1102579},
MRREVIEWER = {Denis\ Serre},
       DOI = {10.2307/2939269},
       URL = {https://doi.org/10.2307/2939269},
}

@article{DoGa,
  author = {M. Dolce and T. Gallay},
  title = {The long way of a viscous vortex dipole},
  journal = {Arch. Ration. Mech. Anal.},
  volume = {250},
  number = {2},
  pages = {Paper No. 18},
  year = {2026},
  doi = {10.1007/s00205-026-02169-5}
}

@article{Nielsen,
    author = {Nielsen, A. H. and Rasmussen, J. Juul},
    title = {Formation and temporal evolution of the Lamb-dipole},
    journal = {Physics of Fluids},
    volume = {9},
    number = {4},
    pages = {982-991},
    year = {1997},
    month = {04},
    issn = {1070-6631},
    doi = {10.1063/1.869193},
    url = {https://doi.org/10.1063/1.869193},
    eprint = {https://pubs.aip.org/aip/pof/article-pdf/9/4/982/19090112/982_1_online.pdf},
}

@article{Sipp,
    author = {Sipp, D. and Jacquin, L. and Cosssu, C.},
    title = {Self-adaptation and viscous selection in concentrated two-dimensional vortex dipoles},
    journal = {Physics of Fluids},
    volume = {12},
    number = {2},
    pages = {245-248},
    year = {2000},
    month = {02},
    issn = {1070-6631},
    doi = {10.1063/1.870325},
    url = {https://doi.org/10.1063/1.870325},
    eprint = {https://pubs.aip.org/aip/pof/article-pdf/12/2/245/19138892/245_1_online.pdf},
}

@article{Trieling,
    author = {Trieling, R. R. and van Wesenbeeck, J. M. A. and van Heijst, G. J. F.},
    title = {Dipolar vortices in a strain flow},
    journal = {Physics of Fluids},
    volume = {10},
    number = {1},
    pages = {144-159},
    year = {1998},
    month = {01},
    issn = {1070-6631},
    doi = {10.1063/1.869556},
    url = {https://doi.org/10.1063/1.869556},
    eprint = {https://pubs.aip.org/aip/pof/article-pdf/10/1/144/19267764/144_1_online.pdf},
}

@article{DeRo2009,
    author = {Delbende, I. and Rossi, M.},
    title = {The dynamics of a viscous vortex dipole},
    journal = {Physics of Fluids},
    volume = {21},
    number = {7},
    pages = {073605},
    year = {2009},
    month = {07},
    issn = {1070-6631},
    doi = {10.1063/1.3183966},
    url = {https://doi.org/10.1063/1.3183966},
    eprint = {https://pubs.aip.org/aip/pof/article-pdf/doi/10.1063/1.3183966/15789524/073605_1_online.pdf},
}

@article {Gold,
    AUTHOR = {M.~A.~Goldshtik},
     TITLE = {A mathematical model of discontinuous incompressible fluid flows},
   JOURNAL = {Dokl. Akad. Nauk SSSR}, 
    VOLUME = {147},
      YEAR = {1962},
    NUMBER = {6},
     PAGES = {1310--1313}, 
       URL = {http://mi.mathnet.ru/dan27378},
}

@article{DPMP,
      title={Global in Time Vortex Configurations for the $2${D} {E}uler Equations}, 
      author={J. D{\'a}vila and M. del Pino and M. Musso and S. Parmeshwar},
    journal = {J. Eur. Math. Soc.},
    year = {2026},
    note = {published online first},
    doi = {10.4171/JEMS/1776}
}

@article {CJ2022kelvin,
    AUTHOR = {Choi, K. and Jeong, I.-J.},
     TITLE = {Stability and instability of {K}elvin waves},
   JOURNAL = {Calc. Var. Partial Differential Equations},
  FJOURNAL = {Calculus of Variations and Partial Differential Equations},
    VOLUME = {61},
      YEAR = {2022},
    NUMBER = {6},
     PAGES = {Paper No. 221, 38},
      ISSN = {0944-2669,1432-0835},
   MRCLASS = {35Q31 (37K58)},
  MRNUMBER = {4493276},
       DOI = {10.1007/s00526-022-02334-0},
       URL = {https://doi.org/10.1007/s00526-022-02334-0},
}

@book{Mazya,
	author = {Maz'ya, V.},
	pages = {xxviii+866},
	publisher = {Springer, Heidelberg},
	series = {Grundlehren der Mathematischen Wissenschaften [Fundamental Principles of Mathematical Sciences]},
	title = {Sobolev spaces with applications to elliptic partial differential equations},
	volume = {342},
	year = {2011}}

@article{JYZ,
	archiveprefix = {arXiv},
	author = {I.-J. Jeong and Y. Yao and T. Zhou},
	journal = {\href{https://arxiv.org/abs/2507.15739}{arXiv:2507.15739}},
	primaryclass = {math.AP},
	title = {Superlinear gradient growth for $2${D} {E}uler equation without boundary},
	year = {2025}}

@book{QS07,
	address = {Basel},
	author = {P. Quittner and P. Souplet},
	doi = {10.1007/978-3-7643-8442-5},
	isbn = {978-3-7643-8441-8},
	publisher = {Birkh{\"a}user},
	series = {Birkh{\"a}user Advanced Texts: Basler Lehrb{\"u}cher},
	title = {Superlinear Parabolic Problems: Blow-up, Global Existence and Steady States},
	url = {https://doi.org/10.1007/978-3-7643-8442-5},
	year = {2007}}

@article{MZ23,
	author = {Masmoudi, N. and Zhao, W.},
	doi = {10.4007/annals.2024.199.3.3},
	fjournal = {Annals of Mathematics. Second Series},
	issn = {0003-486X,1939-8980},
	journal = {Ann. of Math. (2)},
	mrclass = {35Q31 (76E05)},
	mrnumber = {4740211},
	number = {3},
	pages = {1093--1175},
	title = {Nonlinear inviscid damping for a class of monotone shear flows in a finite channel},
	url = {https://doi.org/10.4007/annals.2024.199.3.3},
	volume = {199},
	year = {2024}}

@book{Pra,
	author = {Prandtl, L. and United States. National Advisory Committee for Aeronautics},
	publisher = {National Advisory Committee for Aeronautics},
	series = {Technical memorandum},
	title = {Motion of Fluids with Very Little Viscosity},
	url = {https://books.google.co.kr/books?id=VM3pSgAACAAJ},
	year = {1928}}

@article{Bat,
	author = {Batchelor, G. K.},
	doi = {10.1017/S0022112056000238},
	fjournal = {Journal of Fluid Mechanics},
	issn = {0022-1120,1469-7645},
	journal = {J. Fluid Mech.},
	mrclass = {76.0X},
	mrnumber = {84310},
	mrreviewer = {Y.\ H.\ Kuo},
	pages = {388--398},
	title = {A proposal concerning laminar wakes behind bluff bodies at large {R}eynolds number},
	url = {https://doi.org/10.1017/S0022112056000238},
	volume = {1},
	year = {1956}}

@book{McL,
	author = {McLachlan, R. I.},
	mrclass = {99-05},
	mrnumber = {2638706},
	note = {Thesis (Ph.D.)--California Institute of Technology},
	pages = {119},
	publisher = {ProQuest LLC, Ann Arbor, MI},
	title = {Separated viscous flows via multigrid},
	url = {http://gateway.proquest.com/openurl?url_ver=Z39.88-2004&rft_val_fmt=info:ofi/fmt:kev:mtx:dissertation&res_dat=xri:pqdiss&rft_dat=xri:pqdiss:9031470},
	year = {1990}}

@article{Lions84b,
	author = {Lions, P.-L.},
	journal = {Ann. Inst. H. Poincar\'{e} Anal. Non Lin\'{e}aire},
	pages = {223--283},
	title = {The concentration-compactness principle in the calculus of variations. {T}he locally compact case. {II}},
	volume = {1},
	year = {(1984)}}

@inproceedings{Ben1976,
	address = {Berlin, Heidelberg},
	author = {Benjamin, T. B.},
	booktitle = {Applications of Methods of Functional Analysis to Problems in Mechanics},
	editor = {Germain, Paul and Nayroles, Bernard},
	isbn = {978-3-540-38165-5},
	pages = {8--29},
	publisher = {Springer Berlin Heidelberg},
	title = {The alliance of practical and analytical insights into the nonlinear problems of fluid mechanics},
	year = {1976}}

@article{Cha1903,
	author = {Chaplygin, S. A.},
	doi = {10.1134/S1560354707020074},
	fjournal = {Regular and Chaotic Dynamics. International Scientific Journal},
	issn = {1560-3547,1468-4845},
	journal = {Regul. Chaotic Dyn.},
	mrclass = {76B47 (37N10)},
	mrnumber = {2350307},
	number = {2},
	pages = {219--232},
	title = {One case of vortex motion in fluid},
	url = {https://doi.org/10.1134/S1560354707020074},
	volume = {12},
	year = {2007}}

@book{Lamb1906,
	author = {Lamb, H.},
	edition = {sixth},
	isbn = {0-521-45868-4},
	mrclass = {76-01},
	mrnumber = {1317348},
	note = {With a foreword by R. A. Caflisch [Russel E. Caflisch]},
	pages = {xxvi+738},
	publisher = {Cambridge University Press, Cambridge},
	series = {Cambridge Mathematical Library},
	title = {Hydrodynamics},
	year = {1993}}

@article{AC2022,
	author = {Abe, K. and Choi, K.},
	doi = {10.1007/s00205-022-01782-4},
	fjournal = {Archive for Rational Mechanics and Analysis},
	issn = {0003-9527,1432-0673},
	journal = {Arch. Ration. Mech. Anal.},
	mrclass = {35Q31 (76B47)},
	mrnumber = {4419609},
	number = {3},
	pages = {877--917},
	title = {Stability of {L}amb dipoles},
	url = {https://doi.org/10.1007/s00205-022-01782-4},
	volume = {244},
	year = {2022}}

@article{Arn1966,
	author = {Arnold, V.},
	fjournal = {Universit\'e{} de Grenoble. Annales de l'Institut Fourier},
	issn = {0373-0956,1777-5310},
	journal = {Ann. Inst. Fourier (Grenoble)},
	mrclass = {57.50 (57.55)},
	mrnumber = {202082},
	mrreviewer = {R.\ Hermann},
	pages = {319--361},
	title = {Sur la g\'eom\'etrie diff\'erentielle des groupes de {L}ie de dimension infinie et ses applications \`a{} l'hydrodynamique des fluides parfaits},
	url = {http://www.numdam.org/item?id=AIF_1966__16_1_319_0},
	volume = {16},
	year = {1966}}

@article{BCK2024,
author = {E. Bru{\`e} and M. Colombo and A. Kumar},
title = {{Flexibility of two-dimensional Euler flows with integrable vorticity}},
volume = {175},
journal = {Duke Mathematical Journal},
number = {9},
publisher = {Duke University Press},
pages = {1593 -- 1647},
year = {2026},
doi = {10.1215/00127094-2025-0056},
URL = {https://doi.org/10.1215/00127094-2025-0056}
}

@article{Bur1987,
	author = {Burton, G. R.},
	doi = {10.1007/BF01450739},
	fjournal = {Mathematische Annalen},
	issn = {0025-5831,1432-1807},
	journal = {Math. Ann.},
	mrclass = {49A50 (30C25 35J60 35R35)},
	mrnumber = {870963},
	mrreviewer = {J.\ E.\ Rubio},
	number = {2},
	pages = {225--253},
	title = {Rearrangements of functions, maximization of convex functionals, and vortex rings},
	url = {https://doi.org/10.1007/BF01450739},
	volume = {276},
	year = {1987}}

@article{Bur1988,
	author = {Burton, G. R.},
	doi = {10.1017/S0308210500014669},
	fjournal = {Proceedings of the Royal Society of Edinburgh. Section A. Mathematics},
	issn = {0308-2105,1473-7124},
	journal = {Proc. Roy. Soc. Edinburgh Sect. A},
	mrclass = {35Q20 (35J60 76C05)},
	mrnumber = {943803},
	mrreviewer = {Charles\ J.\ Amick},
	number = {3-4},
	pages = {269--290},
	title = {Steady symmetric vortex pairs and rearrangements},
	url = {https://doi.org/10.1017/S0308210500014669},
	volume = {108},
	year = {1988}}

@article{Bur1996,
	author = {Burton, G. R.},
	doi = {10.1098/rspa.1996.0125},
	fjournal = {Proceedings of the Royal Society. London. Series A. Mathematical, Physical and Engineering Sciences},
	issn = {0962-8444,2053-9169},
	journal = {Proc. Roy. Soc. London Ser. A},
	mrclass = {35Q35 (76C05 76M30)},
	mrnumber = {1421744},
	mrreviewer = {Alexander\ Yurjevich\ Chebotarev},
	number = {1953},
	pages = {2343--2350},
	title = {Uniqueness for the circular vortex-pair in a uniform flow},
	url = {https://doi.org/10.1098/rspa.1996.0125},
	volume = {452},
	year = {1996}}

@article{Bur2005,
	author = {Burton, G. R.},
	doi = {10.1007/s00021-004-0126-6},
	fjournal = {Journal of Mathematical Fluid Mechanics},
	issn = {1422-6928,1422-6952},
	journal = {J. Math. Fluid Mech.},
	mrclass = {76B47 (49N15 76M30)},
	mrnumber = {2126130},
	pages = {S68--S80},
	title = {Isoperimetric properties of {L}amb's circular vortex-pair},
	url = {https://doi.org/10.1007/s00021-004-0126-6},
	volume = {7},
	year = {2005}}

@article{BNL2013,
	author = {Burton, G. R. and Nussenzveig Lopes, H. J. and Lopes Filho, M. C.},
	doi = {10.1007/s00220-013-1806-y},
	fjournal = {Communications in Mathematical Physics},
	issn = {0010-3616,1432-0916},
	journal = {Comm. Math. Phys.},
	mrclass = {35Q30 (35B35 49K10)},
	mrnumber = {3117517},
	mrreviewer = {John\ Albert},
	number = {2},
	pages = {445--463},
	title = {Nonlinear stability for steady vortex pairs},
	url = {https://doi.org/10.1007/s00220-013-1806-y},
	volume = {324},
	year = {2013}}

@article{BCM2025,
	author = {Butt\`a, P. and Cavallaro, G. and Marchioro, C.},
	doi = {10.1137/24M1642391},
	fjournal = {SIAM Journal on Mathematical Analysis},
	issn = {0036-1410,1095-7154},
	journal = {SIAM J. Math. Anal.},
	mrclass = {76B47 (37N10)},
	mrnumber = {4854830},
	number = {1},
	pages = {789--824},
	title = {Leapfrogging vortex rings as scaling limit of {E}uler equations},
	url = {https://doi.org/10.1137/24M1642391},
	volume = {57},
	year = {2025}}

@article{CLZ2021,
	author = {Cao, D. and Lai, S. and Zhan, W.},
	doi = {10.1007/s00526-021-02068-5},
	fjournal = {Calculus of Variations and Partial Differential Equations},
	issn = {0944-2669,1432-0835},
	journal = {Calc. Var. Partial Differential Equations},
	mrclass = {35J60 (35Q31 76B47)},
	mrnumber = {4295232},
	number = {5},
	pages = {Paper No. 190, 16},
	title = {Traveling vortex pairs for 2{D} incompressible {E}uler equations},
	url = {https://doi.org/10.1007/s00526-021-02068-5},
	volume = {60},
	year = {2021}}

@article{CLW2024,
	author = {Cao, D. and Li, S. and Wang, G.},
	date = {2025/10/10},
	doi = {10.1007/s00526-025-03154-8},
	id = {Cao2025},
	isbn = {1432-0835},
	journal = {Calculus of Variations and Partial Differential Equations},
	number = {9},
	pages = {279},
	title = {Desingularization of vortices for the incompressible Euler equation on a sphere},
	url = {https://doi.org/10.1007/s00526-025-03154-8},
	volume = {64},
	year = {2025}}

@article{CQZZ2025,
	author = {Cao, D. and Qin, G. and Zhan, W. and Zou, C.},
	doi = {10.1007/s40818-024-00191-y},
	fjournal = {Annals of PDE. Journal Dedicated to the Analysis of Problems from Physical Sciences},
	issn = {2524-5317,2199-2576},
	journal = {Ann. PDE},
	mrclass = {35J65 (35J20 35Q31 76B03 76B47)},
	mrnumber = {4842908},
	number = {1},
	pages = {Paper No. 1, 55},
	title = {Uniqueness and stability of traveling vortex pairs for the incompressible {E}uler equation},
	url = {https://doi.org/10.1007/s40818-024-00191-y},
	volume = {11},
	year = {2025}}

@article{CWZ2020,
	author = {Cao, D. and Wang, G. and Zhan, W.},
	doi = {10.1137/19M1292151},
	fjournal = {SIAM Journal on Mathematical Analysis},
	issn = {0036-1410,1095-7154},
	journal = {SIAM J. Math. Anal.},
	mrclass = {76B47 (35Q31 76B03)},
	mrnumber = {4169262},
	mrreviewer = {Xinyu\ He},
	number = {6},
	pages = {5363--5388},
	title = {Desingularization of vortices for two-dimensional steady {E}uler flows via the vorticity method},
	url = {https://doi.org/10.1137/19M1292151},
	volume = {52},
	year = {2020}}

@article{Car2002,
	author = {Carley, M.},
	doi = {10.1006/jcph.2002.7107},
	fjournal = {Journal of Computational Physics},
	issn = {0021-9991,1090-2716},
	journal = {J. Comput. Phys.},
	mrclass = {76M23 (76B47)},
	mrnumber = {1917860},
	number = {2},
	pages = {616--641},
	title = {A triangulated vortex method for the axisymmetric {E}uler equations},
	url = {https://doi.org/10.1006/jcph.2002.7107},
	volume = {180},
	year = {2002}}

@Article{Childress.2008,
  author     = {Childress, S.},
  journal    = {Phys. D},
  title      = {Growth of anti-parallel vorticity in {E}uler flows},
  year       = {2008},
  issn       = {0167-2789},
  number     = {14-17},
  pages      = {1921--1925},
  volume     = {237},
  doi        = {10.1016/j.physd.2008.02.028},
  fjournal   = {Physica D. Nonlinear Phenomena},
  mrclass    = {76B47},
  mrnumber   = {2449775},
  mrreviewer = {Bartosz Protas},
  url        = {https://doi.org/10.1016/j.physd.2008.02.028},
}

@Article{Childress.2007,
  author  = {Childress, S.},
  journal = {AML reports 05-07 and 06-07, courant institute of mathematical sciences},
  title   = {Models of vorticity growth in {E}uler flows {I}, {A}xisymmetric flow without swirl and {I}{I}, {A}lmost 2-{D} dynamics},
  year    = {2007},
}

@article{Child1966,
	author = {Childress, S.},
	doi = {10.1063/1.1761786},
	eprint = {https://pubs.aip.org/aip/pfl/article-pdf/9/5/860/12572657/860\_1\_online.pdf},
	issn = {0031-9171},
	journal = {The Physics of Fluids},
	month = {05},
	number = {5},
	pages = {860-872},
	title = {Solutions of {E}uler's Equations containing Finite Eddies},
	url = {https://doi.org/10.1063/1.1761786},
	volume = {9},
	year = {1966}}

@article{CGV2016,
	author = {Childress, S. and Gilbert, A. D. and Valiant, P.},
	doi = {10.1017/jfm.2016.573},
	fjournal = {Journal of Fluid Mechanics},
	issn = {0022-1120,1469-7645},
	journal = {J. Fluid Mech.},
	mrclass = {76B47},
	mrnumber = {3569025},
	pages = {1--30},
	title = {Eroding dipoles and vorticity growth for {E}uler flows in {$\Bbb{R}^3$}: axisymmetric flow without swirl},
	url = {https://doi.org/10.1017/jfm.2016.573},
	volume = {805},
	year = {2016}}

@article{CG2018,
	author = {Childress, S. and Gilbert, A. D.},
	doi = {10.1088/1873-7005/aa9880},
	fjournal = {Fluid Dynamics Research. An International Journal},
	issn = {0169-5983,1873-7005},
	journal = {Fluid Dyn. Res.},
	mrclass = {76B47},
	mrnumber = {3757290},
	number = {1},
	pages = {011418, 40},
	title = {Eroding dipoles and vorticity growth for {E}uler flows in {$\Bbb R^3$}: the hairpin geometry as a model for finite-time blowup},
	url = {https://doi.org/10.1088/1873-7005/aa9880},
	volume = {50},
	year = {2018}}

@article{Choi2024,
	author = {Choi, K.},
	doi = {10.1002/cpa.22134},
	fjournal = {Communications on Pure and Applied Mathematics},
	issn = {0010-3640,1097-0312},
	journal = {Comm. Pure Appl. Math.},
	mrclass = {76B47},
	mrnumber = {4666623},
	mrreviewer = {Weicheng\ Zhan},
	number = {1},
	pages = {52--138},
	title = {Stability of {H}ill's spherical vortex},
	url = {https://doi.org/10.1002/cpa.22134},
	volume = {77},
	year = {2024}}

@article{CJ2021,
	author = {Choi, K. and Jeong, I.-J.},
	doi = {10.1353/ajm.2025.a971091.},
	fjournal = {American Journal of Mathematics},
	journal = {Amer. J. Math.},
	number = {1},
	title = {On vortex stretching for anti-parallel axisymmetric flows},
	volume = {147},
	year = {2025}}

@article{CJ2022,
	author = {Choi, K. and Jeong, I.-J.},
	doi = {10.1016/j.nonrwa.2021.103470},
	fjournal = {Nonlinear Analysis. Real World Applications. An International Multidisciplinary Journal},
	issn = {1468-1218,1878-5719},
	journal = {Nonlinear Anal. Real World Appl.},
	mrclass = {35Q31 (76B47)},
	mrnumber = {4350517},
	pages = {Paper No. 103470, 20},
	title = {Infinite growth in vorticity gradient of compactly supported planar vorticity near {L}amb dipole},
	url = {https://doi.org/10.1016/j.nonrwa.2021.103470},
	volume = {65},
	year = {2022}}

@article{CJ2023,
	author = {Choi, K. and Jeong, I.-J.},
	doi = {10.1080/03605302.2022.2139721},
	fjournal = {Communications in Partial Differential Equations},
	issn = {0360-5302,1532-4133},
	journal = {Comm. Partial Differential Equations},
	mrclass = {76B47 (35Q31)},
	mrnumber = {4555167},
	mrreviewer = {Da-Wen\ Deng},
	number = {1},
	pages = {54--85},
	title = {Filamentation near {H}ill's vortex},
	url = {https://doi.org/10.1080/03605302.2022.2139721},
	volume = {48},
	year = {2023}}

@article{CJS2025_atmos,
    archiveprefix = {arXiv},
	author = {Choi, K. and Jeong, I.-J. and Sim, Y.-J.},
    journal = {\href{https://arxiv.org/abs/2512.10412}{arXiv:2512.10412}},
	title = {Vortex atmospheres of traveling vortices: rigorous definition, existence, and topological classification},
    year = {2025}
	 }

@article{CJS2025,
	author = {Choi, K. and Jeong, I.-J. and Sim, Y.-J.},
	doi = {10.1007/s40818-025-00212-4},
	fjournal = {Annals of PDE. Journal Dedicated to the Analysis of Problems from Physical Sciences},
	issn = {2524-5317,2199-2576},
	journal = {Ann. PDE},
	mrclass = {99-06},
	mrnumber = {4927140},
	number = {2},
	pages = {Paper No. 18},
	title = {On existence of {S}adovskii vortex patch: a touching pair of symmetric counter-rotating uniform vortices},
	url = {https://doi.org/10.1007/s40818-025-00212-4},
	volume = {11},
	year = {2025}}

@article{CJY2024,
    archiveprefix = {arXiv},
	author = {K. Choi and I.-J. Jeong and Y. Yao},
	journal = {\href{https://arxiv.org/abs/2409.19822}{arXiv:2409.19822}},
	title = {Stability of vortex quadrupoles with odd-odd symmetry},
	year = {2024}}

@article{DHLM2025,
    archiveprefix = {arXiv},
	author = {M. Donati and L. E. Hientzsch and C. Lacave and E. Miot},
	journal = {\href{https://arxiv.org/abs/2503.21604}{arXiv:2503.21604}},
	title = {On the dynamics of leapfrogging vortex rings},
	year = {2025}}

@article{FT1981,
	author = {A. Friedman and B. Turkington},
	doi = {10.1090/S0002-9947-1981-0628444-6},
	issn = {0002-9947},
	journal = {Transactions of the American Mathematical Society},
	language = {English},
	number = {1},
	pages = {1--37},
	title = {Vortex rings: Existence and asymptotic estimates},
	volume = {268},
	year = {1981}}

@article{GMT2023,
	author = {Gustafson, S. and Miller, E. and Tsai, T.-P.},
	date = {2026/01/29},
	doi = {10.1007/s00021-026-01001-0},
	id = {Gustafson2026},
	isbn = {1422-6952},
	journal = {Journal of Mathematical Fluid Mechanics},
	number = {1},
	pages = {19},
	title = {Growth Rates for Anti-Parallel Vortex Tube Euler Flows in Three and Higher Dimensions},
	url = {https://doi.org/10.1007/s00021-026-01001-0},
	volume = {28},
	year = {2026}}

@article{HMW2020,
	author = {Hassainia, Z. and Masmoudi, N. and Wheeler, M. H.},
	doi = {10.1002/cpa.21855},
	fjournal = {Communications on Pure and Applied Mathematics},
	issn = {0010-3640,1097-0312},
	journal = {Comm. Pure Appl. Math.},
	mrclass = {35Q31 (76B47)},
	mrnumber = {4156612},
	mrreviewer = {Piotr\ Biler},
	number = {9},
	pages = {1933--1980},
	title = {Global bifurcation of rotating vortex patches},
	url = {https://doi.org/10.1002/cpa.21855},
	volume = {73},
	year = {2020}}

@article{HM2017,
	author = {Hmidi, T. and Mateu, J.},
	doi = {10.1007/s00220-016-2784-7},
	fjournal = {Communications in Mathematical Physics},
	issn = {0010-3616,1432-0916},
	journal = {Comm. Math. Phys.},
	mrclass = {76B47 (35B65 35Q31 35Q35 35Q83 76B03 76E09)},
	mrnumber = {3607460},
	mrreviewer = {Peter\ Bernard\ Weichman},
	number = {2},
	pages = {699--747},
	title = {Existence of corotating and counter-rotating vortex pairs for active scalar equations},
	url = {https://doi.org/10.1007/s00220-016-2784-7},
	volume = {350},
	year = {2017}}

@article{HT2025,
	author = {Huang, D. and Tong, J.},
	doi = {10.1007/s00205-025-02113-z},
	fjournal = {Archive for Rational Mechanics and Analysis},
	issn = {0003-9527,1432-0673},
	journal = {Arch. Ration. Mech. Anal.},
	mrclass = {99-06},
	mrnumber = {4933909},
	number = {4},
	pages = {Paper No. 46},
	title = {Steady contiguous vortex-patch dipole solutions of the 2{D} incompressible {E}uler equation},
	url = {https://doi.org/10.1007/s00205-025-02113-z},
	volume = {249},
	year = {2025}}

@incollection{LDW2016,
	author = {Leweke, T. and Le Diz\`es, S. and Williamson, C. H. K.},
	booktitle = {Annual review of fluid mechanics. {V}ol. 48},
	isbn = {978-0-8243-0748-6},
	mrclass = {76D17 (76B47)},
	mrnumber = {3726395},
	pages = {507--541},
	publisher = {Annual Reviews, Palo Alto, CA},
	series = {Annu. Rev. Fluid Mech.},
	title = {Dynamics and instabilities of vortex pairs},
	volume = {48},
	year = {2016}}

@article{LJ2025,
	author = {Lim, D. and Jeong, I.-J.},
	doi = {10.1007/s00205-025-02103-1},
	fjournal = {Archive for Rational Mechanics and Analysis},
	issn = {0003-9527,1432-0673},
	journal = {Arch. Ration. Mech. Anal.},
	mrclass = {76B47 (35Q35)},
	mrnumber = {4903691},
	number = {3},
	pages = {Paper No. 32},
	title = {On the optimal rate of vortex stretching for axisymmetric {E}uler flows without swirl},
	url = {https://doi.org/10.1007/s00205-025-02103-1},
	volume = {249},
	year = {2025}}

@article{FW2012,
	author = {Luzzatto-Fegiz, P. and Williamson, C. H. K.},
	doi = {10.1017/jfm.2012.255},
	fjournal = {Journal of Fluid Mechanics},
	issn = {0022-1120,1469-7645},
	journal = {J. Fluid Mech.},
	mrclass = {76E07},
	mrnumber = {2971551},
	pages = {323--350},
	title = {Determining the stability of steady two-dimensional flows through imperfect velocity-impulse diagrams},
	url = {https://doi.org/10.1017/jfm.2012.255},
	volume = {706},
	year = {2012}}

@article{Nor1975,
	author = {Norbury, J.},
	doi = {10.1002/cpa.3160280602},
	fjournal = {Communications on Pure and Applied Mathematics},
	issn = {0010-3640,1097-0312},
	journal = {Comm. Pure Appl. Math.},
	mrclass = {35J05},
	mrnumber = {399645},
	mrreviewer = {Clarence\ M.\ Ablow},
	number = {6},
	pages = {679--700},
	title = {Steady planar vortex pairs in an ideal fluid},
	url = {https://doi.org/10.1002/cpa.3160280602},
	volume = {28},
	year = {1975}}

@article{Sad1971,
	author = {V. S. Sadovskii},
	doi = {https://doi.org/10.1016/0021-8928(71)90070-0},
	issn = {0021-8928},
	journal = {Journal of Applied Mathematics and Mechanics},
	number = {5},
	pages = {729-735},
	title = {Vortex regions in a potential stream with a jump of {B}ernoulli's constant at the boundary: PMM vol. 35, no. 5, 1971, pp. 773--779},
	url = {https://www.sciencedirect.com/science/article/pii/0021892871900700},
	volume = {35},
	year = {1971}}

@article{SLF2008,
	author = {Shariff, K. and Leonard, A. and Ferziger, J. H.},
	doi = {10.1016/j.jcp.2007.10.005},
	fjournal = {Journal of Computational Physics},
	issn = {0021-9991,1090-2716},
	journal = {J. Comput. Phys.},
	mrclass = {76M23 (76B47)},
	mrnumber = {2463198},
	number = {21},
	pages = {9044--9062},
	title = {A contour dynamics algorithm for axisymmetric flow},
	url = {https://doi.org/10.1016/j.jcp.2007.10.005},
	volume = {227},
	year = {2008}}

@article{SS2010,
	author = {Smets, D. and Van Schaftingen, J.},
	doi = {10.1007/s00205-010-0293-y},
	fjournal = {Archive for Rational Mechanics and Analysis},
	issn = {0003-9527,1432-0673},
	journal = {Arch. Ration. Mech. Anal.},
	mrclass = {35Q31 (35B25 35J91 76B03 76B47)},
	mrnumber = {2729322},
	mrreviewer = {Stefano\ Bianchini},
	number = {3},
	pages = {869--925},
	title = {Desingularization of vortices for the {E}uler equation},
	url = {https://doi.org/10.1007/s00205-010-0293-y},
	volume = {198},
	year = {2010}}

@article{Tur1983,
	author = {Turkington, B.},
	doi = {10.1080/03605308308820293},
	fjournal = {Communications in Partial Differential Equations},
	issn = {0360-5302,1532-4133},
	journal = {Comm. Partial Differential Equations},
	mrclass = {35Q20 (76C05)},
	mrnumber = {702729},
	mrreviewer = {Wei\ Ming\ Ni},
	number = {9},
	pages = {999--1030, 1031--1071},
	title = {On steady vortex flow in two dimensions. {I}, {II}},
	url = {https://doi.org/10.1080/03605308308820293},
	volume = {8},
	year = {1983}}

@article{Wan1986,
	author = {Wan, Y.-H.},
	fjournal = {Communications in Mathematical Physics},
	issn = {0010-3616,1432-0916},
	journal = {Comm. Math. Phys.},
	mrclass = {76C05 (58E99 76E30)},
	mrnumber = {861881},
	mrreviewer = {Jacob\ Burbea},
	number = {1},
	pages = {1--20},
	title = {The stability of rotating vortex patches},
	url = {http://projecteuclid.org/euclid.cmp/1104115929},
	volume = {107},
	year = {1986}}

@article{Wan1988,
	author = {Wan, Y.-H.},
	doi = {10.1016/0167-2789(88)90056-5},
	fjournal = {Physica D. Nonlinear Phenomena},
	issn = {0167-2789,1872-8022},
	journal = {Phys. D},
	mrclass = {76C05 (58F10)},
	mrnumber = {969034},
	mrreviewer = {Boris\ A.\ Khesin},
	number = {2},
	pages = {277--295},
	title = {Desingularizations of systems of point vortices},
	url = {https://doi.org/10.1016/0167-2789(88)90056-5},
	volume = {32},
	year = {1988}}

@article{Wang2024,
	author = {Wang, G.},
	doi = {10.1090/tran/9105},
	fjournal = {Transactions of the American Mathematical Society},
	issn = {0002-9947,1088-6850},
	journal = {Trans. Amer. Math. Soc.},
	mrclass = {35Q35 (35C07 76B47 76E30)},
	mrnumber = {4744767},
	mrreviewer = {V\'aclav\ M\'acha},
	number = {4},
	pages = {2635--2661},
	title = {On concentrated traveling vortex pairs with prescribed impulse},
	url = {https://doi.org/10.1090/tran/9105},
	volume = {377},
	year = {2024}}

@article{Yang1991,
	author = {Yang, J. F.},
	doi = {10.1142/S021820259100023X},
	fjournal = {Mathematical Models and Methods in Applied Sciences},
	issn = {0218-2025,1793-6314},
	journal = {Math. Models Methods Appl. Sci.},
	mrclass = {35Q35 (35B40 35R35 76C05)},
	mrnumber = {1144183},
	mrreviewer = {A.\ J.\ Meir},
	number = {4},
	pages = {461--475},
	title = {Existence and asymptotic behavior in planar vortex theory},
	url = {https://doi.org/10.1142/S021820259100023X},
	volume = {1},
	year = {1991}}

@article{Lions1984,
	author = {Lions, P.-L.},
	fjournal = {Annales de l'Institut Henri Poincar\'e. Analyse Non Lin\'eaire},
	issn = {0294-1449},
	journal = {Ann. Inst. H. Poincar\'e{} Anal. Non Lin\'eaire},
	mrclass = {49A50 (49A22)},
	mrnumber = {778970},
	mrreviewer = {Gianfranco\ Bottaro},
	number = {2},
	pages = {109--145},
	title = {The concentration-compactness principle in the calculus of variations. {T}he locally compact case. {I}},
	url = {http://www.numdam.org/item?id=AIHPC_1984__1_2_109_0},
	volume = {1},
	year = {1984}}

@article{Pierrehumbert1980,
	author = {Pierrehumbert, R. T.},
	doi = {10.1017/S0022112080000559},
	journal = {Journal of Fluid Mechanics},
	number = {1},
	pages = {129--144},
	title = {A family of steady, translating vortex pairs with distributed vorticity},
	volume = {99},
	year = {1980}}

@article{MV1994,
	author = {Meleshko, V. V. and van Heijst, G. J. F.},
	doi = {10.1017/S0022112094004428},
	fjournal = {Journal of Fluid Mechanics},
	issn = {0022-1120,1469-7645},
	journal = {J. Fluid Mech.},
	mrclass = {76C05 (01A55 01A60 76-03)},
	mrnumber = {1289112},
	mrreviewer = {Oscar\ F.\ Orellana},
	pages = {157--182},
	title = {On {C}haplygin's investigations of two-dimensional vortex structures in an inviscid fluid},
	url = {https://doi.org/10.1017/S0022112094004428},
	volume = {272},
	year = {1994}}

@article{MST1988,
	author = {Moore, D. W. and Saffman, P. G. and Tanveer, S.},
	doi = {10.1063/1.866718},
	eprint = {https://pubs.aip.org/aip/pfl/article-pdf/31/5/978/12738100/978\_1\_online.pdf},
	issn = {0031-9171},
	journal = {The Physics of Fluids},
	month = {05},
	number = {5},
	pages = {978-990},
	title = {The calculation of some {B}atchelor flows: The {S}adovskii vortex and rotational corner flow},
	url = {https://doi.org/10.1063/1.866718},
	volume = {31},
	year = {1988}}

@article{AZ2021,
    archiveprefix = {arXiv},
	author = {Y. A. Antipov and A. Y. Zemlyanova},
	journal = {\href{https://arxiv.org/abs/2010.08118}{arXiv:2010.08118}},
	title = {Sadovskii vortex in a wedge and the associated {R}iemann-{H}ilbert problem on a torus},
	year = {2021}}

@article {Chernyshenko.1988,
    AUTHOR = {Chernyshenko, S. I.},
     TITLE = {Asymptotics of a steady separated flow around a body with
              large {R}eynolds numbers},
   JOURNAL = {Prikl. Mat. Mekh.},
  FJOURNAL = {Akademiya Nauk SSSR. Otdelenie Tekhnicheskikh Nauk. Institut
              Mekhaniki. Prikladnaya Matematika i Mekhanika},
    VOLUME = {52},
      YEAR = {1988},
    NUMBER = {6},
     PAGES = {958--966},
      ISSN = {0032-8235},
   MRCLASS = {76D15},
  MRNUMBER = {985037},
MRREVIEWER = {Frederick\ A.\ Howes},
       DOI = {10.1016/0021-8928(88)90010-X},
       URL = {https://doi.org/10.1016/0021-8928(88)90010-X},
}

@article{Chernyshenko1993,
	author = {Chernyshenko, S. I.},
	doi = {10.1017/S002211209300151X},
	journal = {Journal of Fluid Mechanics},
	pages = {423--431},
	title = {Stratified {S}adovskii flow in a channel},
	volume = {250},
	year = {1993}}

@article{Freilich2016,
	author = {Freilich, D. V.},
	journal = {UC San Diego},
	title = {The {S}adovskii vortex},
	year = {2016}}

@article{FL2017,
	author = {Freilich, D. V. and Llewellyn Smith, S. G.},
	doi = {10.1017/jfm.2017.401},
	journal = {Journal of Fluid Mechanics},
	pages = {479--501},
	title = {The {S}adovskii vortex in strain},
	volume = {825},
	year = {2017}}

@article{ST1982,
	author = {Saffman, P. G. and Tanveer, S.},
	doi = {10.1063/1.863679},
	fjournal = {The Physics of Fluids},
	issn = {0031-9171},
	journal = {Phys. Fluids},
	mrclass = {76C05},
	mrnumber = {680023},
	number = {11},
	pages = {1929--1930},
	title = {The touching pair of equal and opposite uniform vortices},
	url = {https://doi.org/10.1063/1.863679},
	volume = {25},
	year = {1982}}

@article{GS2024,
	author = {Gallay, T. and {\v{S}}ver\'ak, V.},
	doi = {10.2140/apde.2024.17.681},
	fjournal = {Analysis \& PDE},
	issn = {2157-5045,1948-206X},
	journal = {Anal. PDE},
	mrclass = {35Q30 (35Q31)},
	mrnumber = {4713111},
	mrreviewer = {Isabelle\ Gruais},
	number = {2},
	pages = {681--722},
	title = {Arnold's variational principle and its application to the stability of planar vortices},
	url = {https://doi.org/10.2140/apde.2024.17.681},
	volume = {17},
	year = {2024}}

@article{BCM2019,
	author = {Bedrossian, J. and Coti Zelati, M. and Vicol, V.},
	doi = {10.1007/s40818-019-0061-8},
	fjournal = {Annals of PDE. Journal Dedicated to the Analysis of Problems from Physical Sciences},
	issn = {2524-5317,2199-2576},
	journal = {Ann. PDE},
	mrclass = {76B47 (35B40 35P25 35Q31 76B03)},
	mrnumber = {3987441},
	mrreviewer = {Xinyu\ He},
	number = {1},
	pages = {Paper No. 4, 192},
	title = {Vortex axisymmetrization, inviscid damping, and vorticity depletion in the linearized 2{D} {E}uler equations},
	url = {https://doi.org/10.1007/s40818-019-0061-8},
	volume = {5},
	year = {2019}}

@article{BM2015,
	author = {Bedrossian, J. and Masmoudi, N.},
	doi = {10.1007/s10240-015-0070-4},
	fjournal = {Publications Math\'ematiques. Institut de Hautes \'Etudes Scientifiques},
	issn = {0073-8301,1618-1913},
	journal = {Publ. Math. Inst. Hautes \'Etudes Sci.},
	mrclass = {35Q31 (35B35 35Q35)},
	mrnumber = {3415068},
	mrreviewer = {Matthew\ Paddick},
	pages = {195--300},
	title = {Inviscid damping and the asymptotic stability of planar shear flows in the 2{D} {E}uler equations},
	url = {https://doi.org/10.1007/s10240-015-0070-4},
	volume = {122},
	year = {2015}}

@article{IJ2022,
	author = {Ionescu, A. D. and Jia, H.},
	doi = {10.1007/s00205-022-01815-y},
	fjournal = {Archive for Rational Mechanics and Analysis},
	issn = {0003-9527,1432-0673},
	journal = {Arch. Ration. Mech. Anal.},
	mrclass = {35Q31},
	mrnumber = {4487510},
	number = {1},
	pages = {61--137},
	title = {Linear vortex symmetrization: the spectral density function},
	url = {https://doi.org/10.1007/s00205-022-01815-y},
	volume = {246},
	year = {2022}}

@article{IJ2023,
	author = {Ionescu, A. D. and Jia, H.},
	doi = {10.4310/acta.2023.v230.n2.a2},
	fjournal = {Acta Mathematica},
	issn = {0001-5962,1871-2509},
	journal = {Acta Math.},
	mrclass = {76E30 (35Q35 76B03)},
	mrnumber = {4628607},
	mrreviewer = {Yanguang\ (Charles)\ Li},
	number = {2},
	pages = {321--399},
	title = {Non-linear inviscid damping near monotonic shear flows},
	url = {https://doi.org/10.4310/acta.2023.v230.n2.a2},
	volume = {230},
	year = {2023}}

@article{AF1986,
	author = {Amick, C. J. and Fraenkel, L. E.},
	doi = {10.1007/BF00251252},
	fjournal = {Archive for Rational Mechanics and Analysis},
	issn = {0003-9527},
	journal = {Arch. Rational Mech. Anal.},
	mrclass = {35Q20 (76C05)},
	mrnumber = {816615},
	mrreviewer = {Charles\ J.\ Amick},
	number = {2},
	pages = {91--119},
	title = {The uniqueness of {H}ill's spherical vortex},
	url = {https://doi.org/10.1007/BF00251252},
	volume = {92},
	year = {1986}}
			
		\end{document}